\documentclass[12pt]{article}
\usepackage{amssymb,amsmath,amsthm,epsfig}
\usepackage{latexsym, enumerate}
\usepackage{eepic}
\usepackage{epic}
\usepackage{graphicx}
\usepackage{color}
\usepackage{ifpdf}
\usepackage{subfigure}

\usepackage{dsfont}
\usepackage{multirow}
\usepackage{makecell}
\usepackage{algorithm}
\usepackage{bm}

\theoremstyle{plain}

\theoremstyle{definition}

\theoremstyle{remark}

\begin{document}
\title{ \large\bf Identification of the degradation coefficient for an anomalous diffusion process in hydrology}

\author{
Guang-Hui Zheng\thanks{Corresponding author. College of Mathematics and Econometrics, Hunan University, Changsha 410082, Hunan Province, China. Email: zhenggh2012@hnu.edu.cn}
\and
Ming-Hui Ding\thanks{College of Mathematics and Econometrics, Hunan University, Changsha 410082, Hunan Province, China. Email: minghuiding@hnu.edu.cn}
}

\date{}
\maketitle

\begin{center}{\bf ABSTRACT}
\end{center}\smallskip
In hydrology, the degradation coefficient is one of the key parameters to describe the water quality change and to determine the water carrying capacity. This paper is devoted to identify the degradation coefficient in an anomalous diffusion process by using the average flux data at the accessible part of boundary. The main challenges in inverse degradation coefficient problems (IDCP) is the average flux measurement data only provide very limited information and cause the severe ill-posedness of IDCP. Firstly, we prove the average flux measurement data can uniquely determine the degradation coefficient. The existence and uniqueness of weak solution for the direct problem are established, and the Lipschitz continuity of the corresponding forward operator is also obtained. Secondly, to overcome the ill-posedness, we combine the variational regularization method with Laplace approximations (LA) to solve the IDCP. This hybrid method is essentially the combination of deterministic regularization method and stochastic method. Thus, it is able to calculate the minimizer (MAP point) more rapidly and accurately, but also enables captures the statistics information and quantifying the uncertainty of the solution. Furthermore, the existence, stability and convergence of the minimizer of the variational problem are proved. The convergence rate estimate between the LA posterior distribution and the actual posterior distribution in the sense of Hellinger distance is given, and the skewness are introduced for characterizing the symmetry or slope of LA solution, especially the relationship with the symmetry of the measurement data. Finally, the one-dimensional and two-dimensional numerical examples are presented to confirm the efficiency and robustness of the proposed method.

\smallskip
{\bf keywords:} fractional diffusion equation, average flux data, degradation coefficient, Laplace approximations, Hellinger distance, skewness.

\section{Introduction}

In recent years, fractional calculus and fractional differential equations have been more and more extensively used in many
scientific fields. For example, physical, chemical, biology, engineering, medicine, hydrology, finance and so on, refer to \cite{R. Metzler2000, I. Podlubny1999, A.A.Kilbas2006, S.G. Samko1993, V. V. Uchaikin2013, D. Brockmann2006, E. Scalas2000, D. Benson2000, M. Hall2008, B. Henry2008}.

As is known to all, in hydrology, the normal solute diffusion obeys Darcy¡¯s law:
\begin{align}\label{dl}
q=-\kappa(x,u)\nabla u,
\end{align}
and mass conservation law:
\begin{align}\label{mcl}
\frac{\partial u}{\partial t}+\nabla\cdot q=f(x,t,u),
\end{align}
where $q$ is diffusion flux, $\kappa$ is diffusion coefficient, $u$ is concentration of solute and $f$ denotes some source or sink.
By substituting (\ref{dl}) into (\ref{mcl}), the following classical diffusion equation can be derived
\begin{align}\label{de}
\frac{\partial u}{\partial t}-\nabla\cdot(\kappa(x,u)\nabla u)=f(x,t,u).
\end{align}

It is well known that the classical diffusion equation can describe the normal diffusion quite well, and in probability theory, it corresponds to
Brownian motion. However, there are an increasing number of so-called anomalous diffusion arises in real world, especially the diffusion phenomena
occurred in some media with memory and hereditary properties \cite{R. Metzler2000, I. Podlubny1999, A.A.Kilbas2006, S.G. Samko1993, V. V. Uchaikin2013}. The anomalous diffusion is not consistent with
the classical mass conservation law, but satisfies what is called time fractional mass conservation law \cite{S.Y. Lukashchuk2015,G. Wang2015}, i.e.
\begin{align}\label{tfmcl}
_0D_{t}^\alpha u+\nabla\cdot q=f(x,t,u),
\end{align}
where $\alpha\in(0,1)$ ,and $_{0}D_t^{\alpha}u$ is the Caputo fractional derivative defined by \cite{R. Metzler2000, I. Podlubny1999, A.A.Kilbas2006, S.G. Samko1993}
\begin{eqnarray*}
\begin{split}
_{0}D_t^{\alpha}u=\frac{1}{\Gamma(1-\alpha)}\int_0^t(t-s)^{-\alpha}\frac{\partial u(x,s)}{\partial s}ds.
\end{split}
\end{eqnarray*}
Similarly, by combining Darcy¡¯s law (\ref{dl}), one can obtain the time fractional diffusion equation (TFDE)
\begin{align}\label{tfde}
_{0}D_t^{\alpha}u-\nabla\cdot(\kappa(x,u)\nabla u)=f(x,t,u).
\end{align}
The time fractional diffusion equation, compared with the classical Brownian motion, is closely related to fractional Brownian motion and gradually accepted as
an important tool for describing anomalous diffusion. Particularly in hydrology, the TFDE models sticking and trapping between
mobile periods for contaminant particles in a porous medium \cite{R. Schumer2003} or a river flow \cite{P. Chakraborty2009}. About the direct problems for TFDE, i.e., initial
value problem and initial boundary value problem, which have been studied extensively in the past few years \cite{K. Sakamoto2011, S.D. Eidelman2004, Y. Luchko2010, R. Gorenflo2015, Y. Lin and C. Xu2007, Y. N. Zhang2014, F. Zeng2013, K. Mustapha2016, Q. Xu2013, B. Jin2015}. However, in some practical problems,
the boundary data on the whole boundary cannot be obtained. We only know the noisy data on a part of the boundary or at some interior points of the concerned
domain, which will lead to some inverse problems, i.e., fractional inverse diffusion problems. Recently, there are also rapidly growing publications on the time
fractional inverse diffusion problems, such as a tutorial review \cite{B. JinRundell2015}, inverse initial boundary value problems \cite{B. JinRundell2015, J.J. Liu2010, D.A. Murio2008, J. Liu 2015, G.H. Zheng2012, W. Rundell 2013}, inverse source problems \cite{T. Wei2016(1), Y. Zhang2011, Liu2016(1), Wei2014(1), N.H. Tuan2017},
and inverse coefficient problems \cite{J. Cheng2009, G.S. Li2013, Z.D. Zhang2016, B.T. Jin2012, L. MillerYamamoto2013, V.K. Tuan2011, Z. Li2016, T. Wei2018, L.sun2017, SunYanWei2018, T. Wei2016, Weizhang2016, YamamotoZhang2012}.

In hydrology and water resources management, the pollutant degradation coefficient is one of the key parameters to describe the water quality change and to determine
the water carrying capacity \cite{Wang2006,Huang2017}. So, identification of degradation coefficient is very important for water quality evaluation, monitoring and protection.
In anomalous diffusion, identification of degradation coefficient is often correspond to the inverse reaction coefficient problems (IRCP) for TFDE.

As for the inverse reaction (or degradation) coefficient problems for TFDE, there are only limited papers studying this topic. For example, Tuan \cite{V.K. Tuan2011} gave the uniqueness of
IRCP by only using finite measurement data on the boundary. Jin and Rundell \cite{B.T. Jin2012} established the uniqueness in determining the reaction
coefficient from the direct flux measurements for one dimensional TFDE, and an algorithm of the quasi-Newton type is proposed for the numerical reconstruction. In \cite{YamamotoZhang2012}, Yamamoto and Zhang obtained conditional stability in recovering the reaction coefficient in a one dimensional TFDE with one half order Caputo derivative by a Carleman estimate. By means of integral transform method, Miller and Yamamoto \cite{L. MillerYamamoto2013} proved the uniqueness for the IRCP for TFDE from the internal measurement data. Li et al. \cite{G.S. Li2013} suggested an
optimal perturbation algorithm for the simultaneous numerical recovery of the diffusion coefficient and fractional order in a one dimensional TFDE.
In \cite{L.sun2017}, Sun and Wei proved the uniqueness in identifying the reaction coefficient for one dimensional TFDE, and used the conjugate gradient method to solve it numerically. The IRCP for TFDE from the final time data measurement was discussed by Jin and Rundell in \cite{B. JinRundell2015}.

In this paper, based on pollutant degradation model in hydrology and compared with the above one dimensional case which was studied more extensively, we
consider the inverse degradation coefficient problem for TFDE in higher dimensions. Let $\Omega$ be a bounded domain in $\mathds{R}^d$ $(d\geq1)$ with smooth boundary $\partial\Omega$. The pollutant degradation model can be described by TFDE as follows
\begin{eqnarray}
\label{mix-eq1}
\begin{cases}
\begin{split}
_{0}D_t^{\alpha}u-\sum\limits_{k,l=1}^{d}\frac{\partial}{\partial x_k}(a_{kl}(x)\frac{\partial u}{\partial x_l})+q(x)u&=\phi(x)v(t),\ \ \ \text{in}\ \ D=\Omega\times(0, T),\\
u&=0,\ \ \ \text{on}\ \ \partial\Omega\times(0, T),\\
u&=0,\ \ \ \text{in}\ \ \Omega\times\{0\},\\
\end{split}
\end{cases}
\end{eqnarray}
where the degradation coefficient
$q\in C({\overline\Omega})$, $q\geq 0$, the diffusion coefficient Matrix $A=(a_{kl})_{d\times d}\in C^{1}({\overline\Omega})^{{d\times d}}$ and satisfy $a_{kl}=a_{lk}$, $1\leq k,l\leq d$,
\begin{eqnarray*}
\begin{split}
\sum\limits_{k,l=1}^{d} a_{kl}(x)\xi_{k}\xi_{l}\geq \lambda_0\sum\limits_{k=1}^{d}\xi_{k}^{2},\ x\in\overline\Omega,\ \xi\in R^d,\ \lambda_0>0.
\end{split}
\end{eqnarray*}
The input source formed with separated variables $\phi(x)v(t)$, where $v(t)$ is the time-varying strength of source, and $\phi(x)$ denotes the space-position information.
For example, the usual modelling of point source has the form $\phi(x)=\delta(x-a_0)$, here $\delta(x)$ is Dirac function and $a_0$ is the source location.

Throughout this paper, The solutions to system (\ref{mix-eq1}) will be denoted by $u_j^i(x,t;q)$ in order to indicate its dependence on the degradation coefficient $q$, and correspond to the input sources $\phi_j(x)v_i(t)$, $i=1,2$; $j=1,2,\cdots$. Hereafter, $C$ refers to a generic constant which may differ at different occurrences. Then, our aim is to solve the following inverse problem.

\begin{figure}[H]
\centering
  \includegraphics[width=2.5in, height=2in]{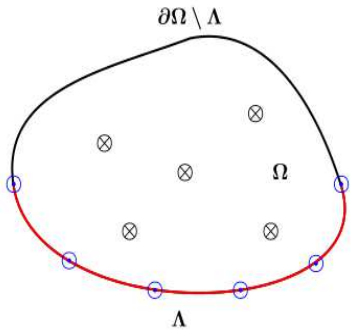}
  \caption{ The physical domain: $\Omega$,\ \ accessible boundary: $\Lambda$,\ \
  inaccessible boundary: $\partial\Omega\setminus\Lambda$,\ \ input source locations: $\otimes$,\ \ measurement locations: $\odot$. }
  \label{schematic of IDCPs}
\end{figure}

\emph{Inverse degradation coefficient problem (IDCP) for TFDE}: Given the input source $\phi_j(x)v_i(t)$, $i=1,2$; $j=1,2,\cdots$, such that the measurement data set
\begin{eqnarray}
\label{mea}
\begin{split}
\int_{\Lambda\times(0,T)}\frac{\partial u_j^i(x,t;q)}{\partial\nu_{A}}h(s,t)dsdt
\end{split}
\end{eqnarray}
determine the degradation coefficient $q$ (see Figure $\ref{schematic of IDCPs}$ for a schematic illustration). Here $\Lambda\subseteq\partial\Omega$ is an accessible part of the boundary, $h$ is a nonzero nonnegative function, and
\begin{eqnarray*}
\begin{split}
\frac{\partial u_j^i}{\partial\nu_{A}}=\sum\limits_{k,l=1}^{d} a_{kl}(x)\frac{\partial u_j^i}{\partial x_k}\nu_l(x),
\end{split}
\end{eqnarray*}
where $\nu_l(x)$ is the $l$th component of the outward unit normal vector $\nu(x)$.
The measurement data set (\ref{mea}) is a weighted integral data on accessible part of the boundary, the weight function $h$ can be interpreted as a characterization of measure instrument.
Compared with the direct flux measurement (or Neumann measurement data) on boundary, i.e. $\frac{\partial u_j^i(x,t;q)}{\partial\nu_{A}}\big|_{\Lambda\times(0,T)}$ (see Jin and Rundell \cite{B.T. Jin2012}), it is a average flux measurement, which is rather easier to measure as a practical matter \cite{Dostert2009,McEnroe2009}. However, because of the average effect, the amounts of average flux measurement data is much less than the amounts of direct flux measurement data. The similar situation also arises in the popular topic called inverse scattering problem with phaseless data \cite{Bao2013, Ammari2016, Zhang2017}. For IDCP, owing to the reduced measurement data, the average flux measurement reconstruction is much more severely ill-posed than the direct flux measurement reconstruction. Moreover, the high nonlinearity is also the inherent difficulty for IDCP. In this paper, we will prove that the limited average flux measurement is enough to identify the degradation coefficient uniquely. Then we combine the variational regularization method with Laplace approximations to solve the IDCP. This hybrid method is essentially the organic combination of deterministic method and stochastic method.

Our contributions of this work are fourfold:
\begin{itemize}
\item{We establish the uniqueness of IDCP, i.e., the average flux measurement data (\ref{mea}) can uniquely determine the degradation coefficient $q$.}
\item{Due to the severe ill-posedness of IDCP, the reconstruction accuracy of degradation coefficient are very sensitive to the data noise, we combine the $L^r$ ($r>1$) variational regularization method and Laplace approximations theory to overcome this difficulty. Especially, it can effectively capture the statistics information and quantifying the uncertainty of the solution.}
\item{We analyze the Hellinger distance between the exact posterior measure and its approximation given by Laplace approximations and prove the second order convergence rate at the MAP point.}
\item{Since the average flux measurement only provide very limited information, such as the symmetry of degradation coefficient. We first use a term called "skewness" from statistics to characterize the symmetry of solution, and show the relation between the symmetry of solution and the symmetry of data. Furthermore, we recover numerically the skewness of degradation coefficient by improving the symmetry of measurement data.}
\end{itemize}

The paper is organized as follows. In Section 2, we present the uniqueness for the IDCP. In Section 3, we consider the well-posedness of weak solution of the direct problem, and prove the  Lipschitz continuity of the corresponding forward operator. Furthermore, the variational regularization scheme with $L^r$ ($r>1$) penalty term is introduced to overcome the ill-posedness of IDCP. The existence, stability, convergence theorem of the minimizer of the variational problem are established. Moreover, the first-order Fr\'{e}chet derivative and second-order Fr\'{e}chet derivative of forward operator are obtained by using adjoint method, and the conjugate gradient method is presented for calculating the minimizer (MAP point). In Section 4, the  Bayesian theory and Laplace approximations method are recalled. The convergence rate estimate between the LA posterior distribution and the actual posterior distribution in the sense of Hellinger distance is given, and the confidence region and skewness are introduced for characterizing the accuracy and reliability of LA. Several one-dimensional and two-dimensional numerical examples are shown in Section 5. Finally, the conclusions are given in Section 6.

\section{The uniqueness of the IDCP}
In this section, the classical solutions of the problem (\ref{mix-eq1}) which belong to the space $C(\overline{D})\cap W^1_t(0,T)\cap C^2_x(\Omega)$ are considered \cite{Y. Luchko2010}, here $W^1_t(0,T)$ denotes the space of functions $f\in C^1(0,T]$ such that $f'\in L^1(0,T)$. Motivated by the ideas in \cite{Badia1999}, we show that the average flux measurement data can determine the degradation coefficient uniquely.\\
\\
\label{anvs}
\textbf{Theorem 1.1.} Let $\{\phi_j\}_{j=1}^{\infty}\in C(\Omega)$ be a complete set in $L^2(\Omega)$, $v\in C^1(0,T)$ and $h\in C_0(\Lambda\times(0,T))$ be given nonzero nonnegative functions, and $v$ satisfies $v(0)=0$. Assume $p(x)$, $q(x)\in C(\overline\Omega)$, $p$, $q\geq 0$ on $\Omega$. Let $u_j^i(x,t;p)$, $u_j^i(x,t;q)$ be the classical solutions of problem (\ref{mix-eq1}) corresponding to the input sources $\phi_j(x)v_i(t)$ ($i=1,2$; $j=1,2,\cdots$) with the degradation coefficients $p$ and $q$ respectively. If we choose $v_1=v$, $v_2={_{0}D_t^\alpha v}$ such that
\begin{eqnarray}
\label{condtion 1}
\begin{split}
\int_{\Lambda\times(0,T)}\frac{\partial u_j^1(x,t;p)}{\partial\nu_{A}}h(x,t)dxdt&=\int_{\Lambda\times(0,T)}\frac{\partial u_j^1(x,t;q)}{\partial\nu_{A}}h(x,t)dxdt;\\
\int_{\Lambda\times(0,T)}\frac{\partial u_j^2(x,t;p)}{\partial\nu_{A}}h(x,t)dxdt&=\int_{\Lambda\times(0,T)}\frac{\partial u_j^2(x,t;q)}{\partial\nu_{A}}h(x,t)dxdt,
\end{split}
\end{eqnarray}
then $q=p$ in $\Omega$.\\
\\
\textbf{Proof}. Notice that $h\in C_0(\Lambda\times(0,T))$ be a given nonzero nonnegative function. Then we set $h_0=h$ on $(\Lambda\times(0,T))$ and $h_0=0$ on $(\partial\Omega\setminus\Lambda)\times(0,T)$, and introduce the function $ w(x,t;q)$ as the solution of the following adjoint problem
\begin{eqnarray}
\label{mix-eq-weak1}
\begin{cases}
\begin{split}
_{t}D_T^{\alpha}w-\sum\limits_{k,l=1}^{d}\frac{\partial}{\partial x_k}(a_{kl}(x)\frac{\partial w}{\partial x_l})+q(x)w&=0,\ \ \ \text{in}\ \ D=\Omega\times(0, T),\\
w&=h_0(x,t),\ \ \ \text{on}\ \ \partial\Omega\times(0, T),\\
w&=0,\ \ \ \text{in}\ \ \Omega\times\{T\}.\\
\end{split}
\end{cases}
\end{eqnarray}
In fact, by using the transform formula $\tilde{w}(x,t)=w(x,T-t)$, (\ref{mix-eq-weak1}) becomes [48]
\begin{eqnarray}
\label{mix-eq-weak2}
\begin{cases}
\begin{split}
_{0}D_t^{\alpha}\tilde{w}-\sum\limits_{k,l=1}^{d}\frac{\partial}{\partial x_k}(a_{kl}(x)\frac{\partial \tilde{w}}{\partial x_l})+q(x)\tilde{w}&=0,\ \ \ \text{in}\ \ D,\\
\tilde{w}&=h_0(x,T-t),\ \ \ \text{on}\ \ \partial\Omega\times(0, T),\\
\tilde{w}&=0,\ \ \ \text{in}\ \ \Omega\times\{0\}.\\
\end{split}
\end{cases}
\end{eqnarray}
Since we choose $v_1=v$, from  $(\ref{mix-eq1})$, $(\ref{mix-eq-weak1})$ and the Green's formula we compute
\begin{eqnarray*}
\label{Qresidual-eq}
\begin{split}
&\int_{0}^T \int_{\Omega}\phi_{j}(x)v(t)w(x,t;q)dxdt\\= & \int_{0}^T\int_{\Omega}[_{0}D_t^{\alpha}u_j^1(x,t;q)-\sum\limits_{k,l=1}^{d}\frac{\partial}{\partial x_k}(a_{kl}(x)\frac{\partial u_j^1(x,t;q)}{\partial x_l})+q(x)u_j^1(x,t;q)]w(x,t;q)dxdt\\
=&\int_0^T \int_{\Omega} ({_{t}D_T^{\alpha}w)u_j^1}dxdt-\int_0^T\int_{\Omega}\sum\limits_{k,l=1}^{d}\frac{\partial}{\partial x_k}(a_{kl}(x)\frac{\partial u_j^1}{\partial x_l})wdxdt+\int_0^T \int_{\Omega}q(x)u_j^1 wdxdt\\
=&\int_0^T \int_{\Omega} ({_{t}D_T^{\alpha}w)u_j^1}dxdt-\int_0^T\int_{\Omega}\sum\limits_{k,l=1}^{d}\frac{\partial}{\partial x_k}(a_{kl}(x)\frac{\partial w}{\partial x_l})u_j^1dxdt+\int_0^T \int_{\Omega}q(x)u_j^1 wdxdt\\&+\int_0^T\int_{\partial\Omega} (\sum\limits_{k,l=1}^{d}a_{kl}(x)\frac{\partial w}{\partial x_l}\nu^k)u_j^1dxdt-\int_0^T\int_{\partial\Omega} (\sum\limits_{k,l=1}^{d}a_{kl}(x)\frac{\partial u_j^1}{\partial x_l}\nu^k)wdxdt\\=&\int_{0}^T\int_{\Omega}[_{t}D_T^{\alpha}w-\sum\limits_{k,l=1}^{d}\frac{\partial}{\partial x_k}(a_{kl}(x)\frac{\partial w}{\partial x_l})+q(x)w]u_j^1dxdt-\int_0^T\int_{\partial\Omega} (\sum\limits_{k,l=1}^{d}a_{kl}(x)\frac{\partial u_j^1}{\partial x_l}\nu^k)wdxdt\\
\end{split}
\end{eqnarray*}
\begin{flalign*}
\text{\ \ $=-\int_0^T\int_{\Lambda}\frac{\partial u_j^1(x,t;q)}{\partial \nu_{A}}h(x,t)dxdt$},&&
\end{flalign*}
where the second equality uses the following fractional integration by parts formula (see \cite{T. Wei2016} Lemma 2.1):
\begin{eqnarray*}
\begin{split}
\int_0^T{_{0}D_t^{\alpha}}u(t)w(t)dt=\int_0^T u(t){_{t}D_T^{\alpha}}w(t)dt,\ \ \ \text{for}\ u,\ w\in AC[0,T], \text{and}\ w(T)=0,
\end{split}
\end{eqnarray*}
here $AC[0,T]$ is the space of functions which are absolutely continuous on $[0,T]$. Similarly, by setting $v_2(t)={_{0}D_t^{\alpha}v}$, it follows that
\begin{eqnarray*}
\begin{split}
\int_{0}^T\int_{\Omega} \phi_{j}(x)_{0}D_t^{\alpha}v(t)w(x,t;q)dxdt=-\int_0^T\int_{\Lambda}\frac{\partial u_j^2(x,t;q)}{\partial \nu_{A}}h(x,t)dxdt.
\end{split}
\end{eqnarray*}
For $p(x)$ and the corresponding function $w(x,t;p)$ given by $(\ref{mix-eq-weak1})$ ,we see that
\begin{eqnarray*}
\begin{split}
\int_{0}^T\int_{\Omega} \phi_{j}(x)v(t)w(x,t;p)dxdt&=-\int_0^T\int_{\Lambda}\frac{\partial u_j^1(x,t;p)}{\partial \nu_{A}}h(x,t)dxdt;\\
\int_{0}^T \int_{\Omega} \phi_{j}(x)_{0}D_t^{\alpha}v(t)w(x,t;p)dxdt&=-\int_0^T\int_{\Lambda}\frac{\partial u_j^2(x,t;p)}{\partial \nu_{A}}h(x,t)dxdt.
\end{split}
\end{eqnarray*}
Then we see from  $(\ref{condtion 1})$ that
\begin{eqnarray*}
\begin{split}
\int_{0}^T\int_{\Omega} \phi_{j}(x)v(t)w(x,t;q)dxdt=&\int_{0}^T \int_{\Omega} \phi_{j}(x)v(t)w(x,t;p)dxdt;\\
\int_{0}^T\int_{\Omega} \phi_{j}(x){_{0}D_t^{\alpha}}v(t)w(x,t;q)dxdt=&\int_{0}^T\int_{\Omega} \phi_{j}(x){_{0}D_t^{\alpha}}v(t)w(x,t;p)dxdt.
\end{split}
\end{eqnarray*}
By the completeness of $\{\phi_j(x)\}_{j=1}^{\infty}$, we obtain
\begin{eqnarray}
\label{eq3}
\begin{split}
\int_{0}^T v(t)w(x,t;q)dt=&\int_{0}^T v(t)w(x,t;p)dt;\\
\int_{0}^T {_{0}D_t^{\alpha}}v(t)w(x,t;q)dt=&\int_{0}^T {_{0}D_t^{\alpha}}v(t)w(x,t;p)dt.\\
\end{split}
\end{eqnarray}
Multiplying equation (\ref{mix-eq-weak1}) by $v$, integrating by parts over $(0,T)$, we find
\begin{eqnarray}
\label{eq2}
\begin{split}
\int_{0}^T{_{t}D_T^{\alpha}}w(x,t;q)vdt-\int_0^T\sum\limits_{k,l=1}^{d}\frac{\partial}{\partial x_k}(a_{kl}(x)\frac{\partial w(x,t;q)}{\partial x_l})vdt+\int_0^T q(x)w(x,t;q)vdt=0;\\
\int_{0}^T{_{t}D_T^{\alpha}}w(x,t;p)vdt-\int_0^T\sum\limits_{k,l=1}^{d}\frac{\partial}{\partial x_k}(a_{kl}(x)\frac{\partial w(x,t;p)}{\partial x_l})vdt+\int_0^T p(x)w(x,t;p)vdt=0.
\end{split}
\end{eqnarray}
The two expressions of $(\ref{eq2})$ are subtracted from each other, and using (\ref{eq3}) we have
\begin{eqnarray}
\label{eq8}
\begin{split}
\int_{0}^T{_{0}D_t^{\alpha}}v[w(x,t;q)-w(x,t;p)]dt=(q-p)\int_0^T w(x,t;p)v(t)dt,
\end{split}
\end{eqnarray}
we finally combine $(\ref{eq3})$, $(\ref{eq8})$ to discover
\begin{eqnarray*}
\begin{split}
(q-p)\int_0^T w(x,t;p)v(t)dt=0.
\end{split}
\end{eqnarray*}
The maximum principle for time factional diffusion equation (\ref{mix-eq-weak1}) or (\ref{mix-eq-weak2}) (see \cite{Y. Luchko2010}, \cite{Al-Refai2015}) can be applied to deduce that $w(x,t;p)>0$, then $q=p$ in $\Omega$.\\
\\
\textbf{Remark 1.2.}
Obviously, when the average flux measurement (\ref{condtion 1}) is replaced by the direct flux measurement and the other assumptions of Theorem 1.1 are satisfied, the uniqueness result still holds, i.e., if
\begin{eqnarray}
\label{condtion 2}
\begin{split}
\frac{\partial u_j^1(x,t;p)}{\partial\nu_{A}}h(x,t)dxdt&=\frac{\partial u_j^1(x,t;q)}{\partial\nu_{A}}h(x,t)dxdt;\\
\frac{\partial u_j^2(x,t;p)}{\partial\nu_{A}}h(x,t)dxdt&=\frac{\partial u_j^2(x,t;q)}{\partial\nu_{A}}h(x,t)dxdt,\ \ j=1,2,\cdots,
\end{split}
\end{eqnarray}
then $q=p$ in $\Omega$. In \cite{B.T. Jin2012}, the direct flux measurement was used by Jin and Rundell for coefficient identification in one dimensional TFDE (see Theorem 3.1 (a) in \cite{B.T. Jin2012}). However,
they assume that the reaction coefficient must be known beforehand in the neighborhood of right boundary $x=1$.\\
\\
\textbf{Remark 1.3.}
Let $\epsilon\in(0, T)$, $v\in H^1(0,T)$ satisfying $v=0$ in $(0, \epsilon]$, $v>0$ in $(\epsilon, T)$, and the other assumptions of Theorem 1.1 are satisfied. By choosing $v_1=v$, $v_2={_{0}D_t^\alpha v}$, we see the uniqueness also holds. The proof is the same as in Theorem 1.1.

\section{Variational regularization method}
\subsection{Weak solution of TFDE}
In order to obtain the continuity of forward map and provide feasibility for numerical computation (such as finite element method), we study the weak solution of problem (\ref{mix-eq1}) with $q(x)\in L^\infty(\Omega)$. The existence, uniqueness and stability of weak solution are given.
Let $_{0}C^{\infty}(0,T)$ denote the space of infinitely differentiable functions on $(0,T)$ with compact support in $(0,T]$. The Sobolev space $_{0}H^s(0,T)$ is the closure of $_{0}C^{\infty}(0,T)$  with respect to the norm $\parallel \cdot\parallel_{H^s(0,T)}$, where $\parallel\cdot\parallel_{H^{s}(0,T)}$ denotes the norm in the usual fractional Sobolev spaces $H^{s}(0,T)$ \cite{ford2011, Lions 1972Non-homogeneous}.

We recall $D=(\Omega\times(0,T))$ and define a Hilbert space
\begin{eqnarray}
\begin{split}
B^s(D):={_{0}H^s}(0,T;L^2(\Omega))\cap L^2(0,T;H_0^1(\Omega)),\ \ s\in(0,1)
\end{split}
\end{eqnarray}
equipped with the norm
\begin{eqnarray}
\begin{split}
\parallel u\parallel_{B^s(D)}=\left(\parallel u\parallel_{{H^s}(0,T;L^2(\Omega))}^2+\parallel u\parallel_{L^2(0,T;H_0^1(\Omega))}^2\right)^{\frac{1}{2}}.
\end{split}
\end{eqnarray}
\textbf{Definition 3.1.1. (Weak solution)} We call that $u\in B^{\frac{\alpha}{2}}(D)$ is a weak solution of the initial boundary value problem $(\ref{mix-eq1})$ provide
\begin{eqnarray}
\begin{split}
A(u,w)=F(w),\ \ \forall w\in B^{\frac{\alpha}{2}}(D),
\end{split}
\end{eqnarray}
where the bilinear form $A(\cdot,\cdot)$ and $F(\cdot)$ are defined by
\begin{eqnarray}
\begin{split}
A(u,w):=&({_{0}D_t^{\frac{\alpha}{2}}}u,{_{t}D_T^{\frac{\alpha}{2}}}w)_{L^2(D)}
 +\int_{D}\sum\limits_{k,l=1}^{d}a_{kl}(x)\frac{\partial u_j}{\partial x_{k}}\frac{\partial w}{\partial x_{l}}dxdt\\&+(qu,w)_{L^2(D)};\\
 F(w):=&(\phi v,w)_{L^2(D)}.
\end{split}
\end{eqnarray}
\textbf{Theorem 3.1.2.} If $q\in L^{\infty}(\Omega)$, $q>0$, and $a_{kl}(x)$, $\phi$, $v$ satisfy the same assumption in Theorem 1.1. Then the problem $(\ref{mix-eq1})$ exists a unique solution $u$ in $B^{\frac{\alpha}{2}}(D)$ and the solution $u$ satisfies
\begin{eqnarray}
\label{eq4}
\begin{split}
\parallel u\parallel_{B^{\frac{\alpha}{2}}(D)}\leq C \parallel \phi v\parallel_{L^2(D)},
\end{split}
\end{eqnarray}
where $C$ is a constant independent of $u$.

By using the Lax-Milgram theorem, the proof of Theorem 2.1 is standard (refer to \cite{ford2011, Y. Lin and C. Xu2007, T. Wei2016}). Hence, we omit it.

\subsection{The forward operator}
In this section, we solve numerically the reaction coefficient $q(x)$ by problem $(\ref{mix-eq1})$.  The inverse coefficient problem is formulated into a variational problem by using the Tikhonov regularization. Then the existence, stability and convergence of minimizer for the variational problem are provided.\\
Define a forward or solution operator
\begin{eqnarray*}
\begin{split}
\mathcal{F}:q(x)\in Q\mapsto\Phi\in R^{2\times N},
\end{split}
\end{eqnarray*}
where $Q=\{{q(x)\in L^{\infty}(\Omega)}\mid q_{\textrm{min}}\leq q(x)\leq q_{\textrm{max}}\}$, and
\begin{eqnarray}
\begin{split}
\Phi=(\varphi_{ij})_{2\times N},\ \ \varphi_{ij}=\int_{\partial\Lambda\times(0,T)}\frac{\partial u_j^i(q;x,t)}{\partial\nu_A}h(x,t)dxdt,\ \ i=1,2;\ j=1,2,...,N,
\end{split}
\end{eqnarray}
and the norm of $R^{2\times N}$ is defined by $\parallel\Phi\parallel_{s}=(\sum\limits_{i=1}^{2}\sum\limits_{j=1}^{N}\mid\varphi_{ij}\mid^s)^{\frac{1}{s}},\ s\geq1$.
To get the continuity of the forward operator, we restrict $h\in L^2(0,T;H_0^{\frac{1}{2}}(\Lambda))$.\\
\\
\textbf{Theorem 3.2.1.} If $h\in L^2(0,T;H_0^{\frac{1}{2}}(\Lambda))$, and the other assumption is the same as Theorem 1.1. Then the nonlinear forward map $\mathcal{F}:Q\mapsto R^{2\times N}$ is the Lipschitz continuous.\\
\\
\textbf{Proof}. Setting $u_j^i=u_{j}^i(x,t;q),\ \tilde{u}_j^i=u_{j}^i(x,t;\tilde{q})$, it follows that
\begin{eqnarray*}
\label{mix-eq-weak}
\begin{cases}
\begin{split}
_{0}D_t^{\alpha}u_j^i(x,t)-\sum\limits_{k,l=1}^{d}\frac{\partial}{\partial x_k}(a_{kl}(x)\frac{\partial u_j^i(x,t)}{\partial x_l})+q(x)u_j^i(x,t)&=\phi_j(x)v_i(t),\ x\in\Omega,0\leq t\leq T;\\
u_j^i(x,t)&=0,\ x\in\partial\Omega,0\leq t\leq T;\\
u_j^i(x,0)&=0,\ x\in\Omega,\\
\end{split}
\end{cases}
\end{eqnarray*}
and
\begin{eqnarray*}
\label{mix-eq-weak}
\begin{cases}
\begin{split}
_{0}D_t^{\alpha}\tilde{u}_j^i(x,t)-\sum\limits_{k,l=1}^{d}\frac{\partial}{\partial x_k}(a_{kl}(x)\frac{\partial \tilde{u}_j^i(x,t)}{\partial x_l})+\tilde{q}(x)\tilde{u}_j^i(x,t)&=\phi_j(x)v_i(t),\ x\in\Omega,0\leq t\leq T;\\
\tilde{u}_j^i(x,t)&=0,\ x\in\partial\Omega,0\leq t\leq T;\\
\tilde{u}_j^i(x,0)&=0,\ x\in\Omega.\\
\end{split}
\end{cases}
\end{eqnarray*}
Let $z_{j}^i=u_j^i-\tilde{u}_j^i$, then $z_j^i$ satisfies the following equation
\begin{eqnarray*}
\label{mix-eq-weak}
\begin{cases}
\begin{split}
_{0}D_t^{\alpha}z_j^i(x,t)-\sum\limits_{k,l=1}^{d}\frac{\partial}{\partial x_k}(a_{kl}(x)\frac{\partial z_j^i(x,t)}{\partial x_l})+q(x)z_j^i(x,t)&=(\tilde{q}(x)-q(x))\tilde{u}_j^i,\ x\in\Omega,0\leq t\leq T;\\
z_j^i(x,t)&=0,\ x\in\partial\Omega,0\leq t\leq T;\\
z_j^i(x,0)&=0,\ x\in\Omega.\\
\end{split}
\end{cases}
\end{eqnarray*}
Since $v_i(t)\in H^1[0,T];\ \phi_j\in L^2(\Omega)$ and $v_i(t)\phi_j(x)\in H^1(0,T;L^2(\Omega))\subseteq L^2(\Omega\times(0,T))$, by using Theorem 3.1.2, we have
\begin{eqnarray}
\label{eq5}
\begin{split}
\parallel \tilde{u}_j^i\parallel_{B^{\frac{\alpha}{2}}(D)}\leq C\parallel\phi_j(x)v_i(t)\parallel_{L^2(D)}.
\end{split}
\end{eqnarray}
Similarly from $(\ref{eq5})$
\begin{eqnarray}
\label{seneq}
\begin{split}
\parallel z_j^i\parallel_{B^{\frac{\alpha}{2}}(D)}\leq& C\parallel(\tilde{q}-q)\tilde{u}_j^i\parallel_{L^2(D)}\\
\leq&C\parallel\tilde{q}-q\parallel_{L^{\infty}(\Omega)}\parallel \tilde{u}_j^i\parallel_{B^{\frac{\alpha}{2}}(D)}\\
\leq&C\parallel\tilde{q}-q\parallel_{L^{\infty}(\Omega)}\parallel\phi_j(x)v_i(t)\parallel_{L^2(D)}\\
\leq&C\parallel\tilde{q}-q\parallel_{L^{\infty}(\Omega)}.
\end{split}
\end{eqnarray}
Furthermore, by trace theorem, it implies
\begin{eqnarray*}
\begin{split}
\left\| \frac{\partial z_j^i}{\partial \nu_A}\right\|_{L^{{2}}(0,T;H^{-\frac{1}{2}}(\Lambda))}\leq& C\parallel z_j^i\parallel_{B^{\frac{\alpha}{2}}(\Omega\times(0,T))}\leq&C\parallel\tilde{q}-q\parallel_{L^{\infty}(\Omega)},
\end{split}
\end{eqnarray*}
and then
\begin{eqnarray*}
\begin{split}
\mid\tilde{\varphi}_{ij}-\varphi_{ij}\mid=& \Big|\left\langle\frac{\partial u_j^i}{\partial\nu_A}-\frac{\partial \tilde{u}_j^i}{\partial \nu_A},h\right\rangle \Big|\\
\leq&\left\|\frac{\partial z_j^i}{\partial \nu_A}\right\|_{L^2(0,T;H^{-\frac{1}{2}}(\Lambda))}\parallel h\parallel_{L^2(0,T;H^{\frac{1}{2}}(\Lambda))}\\
\leq&C\parallel\tilde{q}-q\parallel_{L^{\infty}(\Omega)}.
\end{split}
\end{eqnarray*}
Therefore, the Lipschitz continuity of $\mathcal{F}$ is proved.

\subsection{Variational regularization method with $L^r$ ($r>1$) penalty term}
In order to overcome the ill-posedness of the problem, we apply the variational regularization method with $L^r$ ($r>1$) penalty term to deal with it. Then the corresponding variational functional is defined as follows
\begin{eqnarray}
\label{eq9}
\begin{split}
J(q)=\frac{1}{s}\parallel\mathcal{F}(q)-\Phi^{\delta}\parallel_{s}^s+\frac{\mu}{r}\parallel q\parallel_{L^r(\Omega)}^{r},\ \ r>1,\ s\geq 1,
\end{split}
\end{eqnarray}
where $\Phi^\delta=(\varphi_{ij}^{\delta})_{2\times N}$ denotes the measure data matrix and satisfies $\parallel\Phi^{\delta}-\Phi\parallel_{s}\leq\delta$, here $\delta$ is the noise level. $\mu$ is the regularization parameter. Next, we will prove the existence, stability and convergence of minimizer of the variational functional (\ref{eq9}). Although the proof is similar to \cite{L.sun2017, hofmann2007}, for the completeness, we give some details.\\
\\
\textbf{Theorem 3.3.1.} There exists a minimizer $q_\mu^\delta\in Q$ for variation functional $J(q)$.\\
\\
\textbf{Proof}. Because of the nonnegativity of $J(q)$, there is a minimizing sequence $\{q_k\}$ in $Q$ such that
\begin{eqnarray*}
\begin{split}
J_0=\inf\limits_{q\in Q} J(q),\ \lim\limits_{k\rightarrow \infty}J(q_k)=J_0.
\end{split}
\end{eqnarray*}
Since $\{q_k\}\subseteq Q$ is bounded in $L^{\infty}(\Omega)$, there exists a subsequence, which is again denoted by $\{q_k\}$, such that $q_k\stackrel{\ast}\rightharpoonup q_0$ in $L^{\infty}(\Omega)$. Moreover, by using the reflexivity of $L^r(\Omega)$ and the density of $L^r(\Omega)$ $(r>1)$ in $L^1(\Omega)$ \cite{brezis2011}, it follows that $q_k\rightharpoonup q_0$ in $L^r(\Omega)$. Owing to the closed convexity of $Q$, by using Mazur Theorem, we find $q_k\rightarrow q_0$, and then $q_0\in Q$. From the weak lower semicontinuity of norm, and the Lipschitz continuity of $\mathcal{F}$, we have
\begin{eqnarray*}
\begin{split}
J(q_0)=&\frac{1}{s}\parallel\mathcal{F}(q_0)-\Phi^{\delta}\parallel_{s}^s+\frac{\mu}{r}\parallel q_0\parallel_{L^r(\Omega)}^{r}\\
\leq&\liminf\limits_{k\rightarrow\infty}(\frac{1}{s}\parallel\mathcal{F}(q_k)-\Phi^{\delta}\parallel_{s}^s+\frac{\mu}{r}\parallel q_k\parallel_{L^r(\Omega)}^{r})\\
=&\liminf\limits_{k\rightarrow\infty} J(q_k)\\
=&J_0.
\end{split}
\end{eqnarray*}
Therefore, $q_0$ is a minimizer of $J(q)$.\\
\\
\textbf{Theorem 3.3.2.} Assume that $\{\Phi_k\}$ is a sequences which satisfy $\Phi_k\rightarrow \Phi^\delta,\ k\rightarrow\infty$ in $R^{2\times N}$, and $\{q_k\}$ is a minimizer of $J(q)$ with $\Phi^\delta$ replaced by $\Phi_k$. Then the minimizers of $J(q)$ are stable with respect to the measurement data $\Phi^\delta$.\\
\\
\textbf{Proof}. From the definition of $\{q_k\}$, we find
\begin{eqnarray}
\label{eq6}
\begin{split}
\frac{1}{s}\parallel\mathcal{F}(q_k)-\Phi_{k}\parallel_{s}^s+\frac{\mu}{r}\parallel q_k\parallel_{L^r(\Omega)}^{r}
\leq\frac{1}{s}\parallel\mathcal{F}(q)-\Phi_{k}\parallel_{s}^s+\frac{\mu}{r}\parallel q\parallel_{L^r(\Omega)}^{r},\ \ \forall q\in Q.
\end{split}
\end{eqnarray}
Since $\{q_k\}\subseteq Q$ is bounded in $L^{\infty}(\Omega)$, there has a subsequence, still denoted by $\{q_k\}$, such that $q_k\stackrel{\ast}\rightharpoonup q_0$. Similar to the proof of Theorem 3.3.1, it implies $q_k\rightarrow q_0$. Based on continuity of $\mathcal{F}$ and weak lower semicontiunity of norm, by $(\ref{eq6})$ it follows that
\begin{eqnarray*}
\begin{split}
&\frac{1}{s}\parallel\mathcal{F}(q_0)-\Phi^\delta\parallel_{s}^s+\frac{\mu}{r}\parallel q_0\parallel_{L^r(\Omega)}^{r}\\
&\leq\liminf\limits_{k\rightarrow\infty}(\frac{1}{s}\parallel\mathcal{F}(q_k)-\Phi_{k}\parallel_{s}^s+\frac{\mu}{r}\parallel q_k\parallel_{L^r(\Omega)}^{r})\\
&\leq\limsup\limits_{k\rightarrow\infty}(\frac{1}{s}\parallel\mathcal{F}(q_k)-\Phi_{k}\parallel_{s}^s+\frac{\mu}{r}\parallel q_k\parallel_{L^r(\Omega)}^{r})\\
&\leq\lim\limits_{k\rightarrow\infty}(\frac{1}{s}\parallel\mathcal{F}(q)-\Phi_{k}\parallel_{s}^s+\frac{\mu}{r}\parallel q\parallel_{L^r(\Omega)}^{r})\\
&=\frac{1}{s}\parallel\mathcal{F}(q)-\Phi^\delta\parallel_{s}^s+\frac{\mu}{r}\parallel q\parallel_{L^r(\Omega)}^{r},\\
\end{split}
\end{eqnarray*}
for all $q\in Q$. This deduce that $q_0$ is a minimizer of $J(q)$. Furthermore, we set $q=q_0$, and find that
\begin{eqnarray*}
\begin{split}
\frac{1}{s}\parallel\mathcal{F}(q_0)-\Phi^\delta\parallel_{s}^s+\frac{\mu}{r}\parallel q_0\parallel_{L^r(\Omega)}^{r}
=\lim\limits_{k\rightarrow\infty}(\frac{1}{s}\parallel\mathcal{F}(q_k)-\Phi_{k}\parallel_{s}^s+\frac{\mu}{r}\parallel q_k\parallel_{L^r(\Omega)}^{r}).
\end{split}
\end{eqnarray*}
Thus, the stability result holds.\\
\\
\textbf{Definition 3.3.3.} $q^\dagger\in Q$ is called an $L^r$-minimizing solution if\\
\begin{eqnarray}
\begin{split}
\parallel q^\dagger\parallel_{L^r(\Omega)}=\min\limits_{\mathcal{F}(q)=\Phi}\parallel q\parallel_{L^r(\Omega)}.
\end{split}
\end{eqnarray}
\textbf{Theorem 3.3.4.} Suppose that the noise level sequence $\{\delta_k\}$ convergence monotonically to $0$, and the corresponding measurement data $\Phi^{\delta_k}$
satisfy $\parallel \Phi^{\delta_k}-\Phi\parallel_s\leq\delta_k$. Moreover, assume that the regularization parameter $\mu(\delta)$ satisfies $\mu(\delta)\rightarrow 0$, and $\frac{\delta^s}{\mu(\delta)}\rightarrow 0$, (as $\delta\rightarrow 0$), and $\alpha(\delta)$ is also monotonically increasing. $\{q_{\mu(\delta_k)}^{\delta_k}\}$ is a minimizer of $J(q)$ with $\Phi^\delta$ replaced by $\Phi^{\delta_k}$. Then $\{q_{\mu(\delta_k)}^{\delta_k}\}$ has a subsequence which convergent to an $L^r$-minimizing solution, and satisfies that $\parallel q^\dagger\parallel_{L^r(\Omega)}=\lim\limits_{k\rightarrow \infty}\parallel q_{\mu(\delta_k)}^{\delta_k}\parallel_{L^r(\Omega)}$.\\
\\
\textbf{Proof}. From the definition of $q_{\mu(\delta_k)}^{\delta_k}$, we get
\begin{eqnarray*}
\begin{split}
&\frac{1}{s}\parallel\mathcal{F}(q_{\mu(\delta_k)}^{\delta_k})-\Phi^{\delta_k}\parallel_{s}^s+\frac{\mu(\delta_k)}{r}\parallel q_{\mu(\delta_k)}^{\delta_k}\parallel_{L^r(\Omega)}^{r}\\
&\leq\frac{1}{s}\parallel\mathcal{F}(q^\dagger)-\Phi_{\delta_k}\parallel_{s}^s+\frac{\mu(\delta_k)}{r}\parallel q^\dagger\parallel_{L^r(\Omega)}^{r}\\
&\leq\frac{1}{s}\delta_k^s+\frac{\mu(\delta_k)}{r}\parallel q^\dagger\parallel_{L^r(\Omega)}^r.
\end{split}
\end{eqnarray*}
Then, take $k\rightarrow \infty$, we obtain
\begin{eqnarray}
\begin{split}
\lim\limits_{k\rightarrow\infty}\mathcal{F}(q_{\mu(\delta_k)}^{\delta_k})=\Phi.
\end{split}
\end{eqnarray}
Moreover, notice that
\begin{eqnarray*}
\begin{split}
\parallel q_{\mu(\delta_k)}^{\delta_k}\parallel_{L^r(\Omega)}^r
\leq \frac{r\delta_k^s}{s\mu(\delta_k)}+\parallel q^\dagger\parallel_{L^r(\Omega)}^r,
\end{split}
\end{eqnarray*}
and $\frac{\delta^s}{\mu(\delta)}\rightarrow 0,\ (k\rightarrow\infty)$, it implies that
\begin{eqnarray}
\label{eq7}
\begin{split}
\limsup\limits_{k\rightarrow\infty}\parallel q_{\mu(\delta_k)}^{\delta_k} \parallel_{L^r(\Omega)}^r\leq\parallel q^\dagger\parallel_{L^r(\Omega)}^r.
\end{split}
\end{eqnarray}
Since $\{q_{\mu(\delta_k)}^{\delta_k}\}\subseteq Q$ is bounded, there has a subsequence, which denoted again by $\{q_{\mu(\delta_k)}^{\delta_k}\}$ satisfying $q_{\mu(\delta_k)}^{\delta_k}\stackrel{\ast}\rightharpoonup q_0\in L^\infty(\Omega)$. Likewise, by Mazur Theorem, we also get $q_{\mu(\delta_k)}^{\delta_k}\rightarrow q_0\in L^\infty(\Omega)$. Hence, by (\ref{eq7}), a similar argument implies that
\begin{eqnarray}
\begin{split}
\parallel q_{0}\parallel_{L^r(\Omega)}^r&\leq\liminf\limits_{k\rightarrow\infty}\parallel q_{\mu(\delta_k)}^{\delta_k} \parallel_{L^r(\Omega)}^r\leq\limsup\limits_{k\rightarrow\infty}\parallel q_{\mu(\delta_k)}^{\delta_k}\parallel_{L^r(\Omega)}^r\\
&\leq\parallel q^\dagger\parallel_{L^r(\Omega)}^r\leq\parallel q\parallel_{L^r(\Omega)}^r,
\end{split}
\end{eqnarray}
for all $q\in Q$ satisfying $\mathcal{F}(q)=\Phi$. Setting $q=q_0$ deduce that $\parallel q^\dagger\parallel_{L^r(\Omega)}=\parallel q_0\parallel_{L^r(\Omega)}$. That is, $q_{0}$ is an $L^r$-minimizing solution, and satisfies $\parallel q^\dagger\parallel_{L^r(\Omega)}=\parallel q_0\parallel_{L^r(\Omega)}=\lim\limits_{k\rightarrow\infty}\parallel q_{\mu(\delta_k)}^{\delta_k}\parallel_{L^r(\Omega)}$.\\

\subsection{The Fr\'{e}chet derivative of variation functional}
For simplicity, we only focus on the case $s=r=2$, then the variation functional becomes
\begin{eqnarray}
\label{3.29}
\begin{split}
J(q)=\frac{1}{2}\parallel\mathcal{F}(q)-\Phi^{\delta}\parallel_{2}^2+\frac{\mu}{2}\parallel q\parallel_{L^2(\Omega)}^{2}.
\end{split}
\end{eqnarray}
In order to find the minimizer of variation functional (\ref{3.29}), the efficient evaluation of the Fr\'{e}chet derivative is critical, here we adopt the adjoint method. First, we claim that the following asymptotic expansion formula of solution holds.\\
\\
\textbf{Theorem 3.4.1.} The solution $u_j^i(q)$ is differentiable in the sense that: for any direction $\delta q\in L^\infty(\Omega)$, we have that
\begin{eqnarray}
\label{}
\begin{split}
u_j^i(q+\delta q)=u_j^i(q)+\vartheta_j^i(q)[\delta q]+o(\|\delta q\|_{L^\infty(\Omega)}),
\end{split}
\end{eqnarray}
where $\vartheta_j^i(q)[\delta q]$ satisfies the \textbf{sensitive equation}:
\begin{eqnarray}
\label{3.31}
\begin{cases}
\begin{split}
_{0}D_t^{\alpha}\vartheta_j^i(q)[\delta q]&-\sum\limits_{k,l=1}^{d}\frac{\partial}{\partial x_k}\left(a_{kl}(x)\frac{\partial \vartheta_{j}^i(q)[\delta q]}{\partial x_l}\right)+q\cdot\vartheta_{j}^i(q)[\delta q]=-\delta q\cdot u_j^i(q),\ \ \ \text{in}\ \ D,\\
\vartheta_j^i(q)[\delta q]&=0,\ \ \ \text{on}\ \ \partial\Omega\times(0, T),\\
\vartheta_j^i(q)[\delta q]&=0,\ \ \ \text{in}\ \ \Omega\times\{0\}.
\end{split}
\end{cases}
\end{eqnarray}
\textbf{Proof}. Setting $\tilde{\vartheta}_j^i=u_j^i(q+\delta q)-u_j^i(q)-\vartheta_j^i(q)[\delta q]$. It is easy to verify that $\tilde{\vartheta}_j^i$ satisfies
\begin{eqnarray}
\label{}
\begin{cases}
\begin{split}
_{0}D_t^{\alpha}\tilde{\vartheta}_j^i&-\sum\limits_{k,l=1}^{d}\frac{\partial}{\partial x_k}\left(a_{kl}(x)\frac{\partial \tilde{\vartheta}_{j}^i}{\partial x_l}\right)+(q+\delta q)\tilde{\vartheta}_{j}^i=-\delta q\cdot \vartheta_j^i(q)[\delta q],\ \ \ \text{in}\ \ D,\\
\tilde{\vartheta}_j^i&=0,\ \ \ \text{on}\ \ \partial\Omega\times(0, T),\\
\tilde{\vartheta}_j^i&=0,\ \ \ \text{in}\ \ \Omega\times\{0\}.
\end{split}
\end{cases}
\end{eqnarray}
By Theorem 3.1.2, it deduces that $\|\tilde{\vartheta}_j^i\|_{B^{\frac{\alpha}{2}}(\Omega)}\leq\|\delta q\|_{L^\infty(\Omega)}\|\vartheta_j^i(q)[\delta q]\|_{L^2(\Omega)}$. Similarly, by using Theorem 3.1.2 to (\ref{3.31}), we find  $\|\vartheta_j^i(q)[\delta q]\|_{B^{\frac{\alpha}{2}}(\Omega)}\leq\|\delta q\|_{L^\infty(\Omega)}\|u_j^i(q)\|_{L^2(\Omega)}$, and then the claim now holds.

The following theorem shows that the solution $u_j^i(q)$ is twice Fr\'{e}chet differentiable. Since the proof is analogous to Theorem 3.4.1 and is omitted.\\
\\
\textbf{Theorem 3.4.2.} $\vartheta_j^i(q)[\delta q]$ is differentiable in the sense that: for the direction $\tilde{\delta} q\in L^\infty(\Omega)$, we have that
\begin{eqnarray*}
\label{}
\begin{split}
\vartheta_j^i(q+\tilde{\delta} q)[\delta q]=\vartheta_j^i(q)[\delta q]+\zeta_j^i(q)[\delta q,\tilde{\delta} q]+o(\|\tilde{\delta} q\|_{L^\infty(\Omega)}),\ \ \ \text{as}\ \ \delta q,\tilde{\delta} q\rightarrow0\ \ \text{in}\ L^\infty(\Omega),
\end{split}
\end{eqnarray*}
where $\zeta_j^i(q)[\delta q,\tilde{\delta} q]$ satisfies the \textbf{second-order sensitive equation}:
\begin{eqnarray}
\label{3.33}
\begin{cases}
\begin{split}
_{0}D_t^{\alpha}\zeta_j^i(q)[\delta q,\tilde{\delta} q]&-\sum\limits_{k,l=1}^{d}\frac{\partial}{\partial x_k}\left(a_{kl}(x)\frac{\partial \zeta_{j}^i(q)[\delta q,\tilde{\delta} q]}{\partial x_l}\right)+q\cdot\zeta_{j}^i(q)[\delta q,\tilde{\delta} q]\\
&\ \ \ \ \ \ \ \ \ \ \ \ \ \ \ \ \ \ \ \ \ =-\tilde{\delta} q\cdot \vartheta_j^i(q)[\delta q]-\delta q\cdot\vartheta_j^i(q)[\tilde{\delta} q],\ \ \ \text{in}\ \ D,\\
\zeta_j^i(q)[\delta q,\tilde{\delta} q]&=0,\ \ \ \text{on}\ \ \partial\Omega\times(0, T),\\
\zeta_j^i(q)[\delta q,\tilde{\delta} q]&=0,\ \ \ \text{in}\ \ \Omega\times\{0\}.
\end{split}
\end{cases}
\end{eqnarray}
\\
\\
\textbf{Theorem 3.4.3.} For the direction $\delta q\in L^\infty(\Omega)$, the Fr\'{e}chet derivative of variation functional $J(q)$ at $q\in Q$ is given by
\begin{eqnarray}
\label{}
\begin{split}
J'(q)[\delta q]&=\sum_{i=1}^2\sum_{j=1}^{N}\int_0^T\int_\Omega\delta q(x)\cdot u_{j}^i(x,t)\varpi_{j}^i(x,t)dxdt+\mu\int_\Omega\delta q(x)\cdot q(x)dx,
\end{split}
\end{eqnarray}
where $\varpi_j^i$ satisfies the \textbf{adjoint equation}:
\begin{eqnarray}
\label{3.35}
\begin{cases}
\begin{split}
_{t}D_T^{\alpha}\varpi_{j}^i&-\sum\limits_{k,l=1}^{d}\frac{\partial}{\partial x_k}\left(a_{kl}(x)\frac{\partial \varpi_{j}^i}{\partial x_l}\right)+q(x)\varpi_{j}^i=0,\ \ \ \text{in}\ \ D,\\
\varpi_j^i&=h\left(\int_{\Lambda\times(0,T)}\frac{\partial u_{j}^i(q)}{\partial\nu_A}h(x,t)dxdt-\varphi_{ij}^\delta\right),\ \ \ \text{on}\ \ \Lambda\times(0, T),\\
\varpi_j^i&=0,\ \ \ \text{on}\ \ (\partial\Omega\setminus\Lambda)\times(0, T),\\
\varpi_j^i&=0,\ \ \ \text{in}\ \ \Omega\times\{0\}.
\end{split}
\end{cases}
\end{eqnarray}
\textbf{Proof}. In fact, we only need to prove that the Fr\'{e}chet derivative of $J_1(q)=\frac{1}{2}\parallel\mathcal{F}(q)-\Phi^{\delta}\parallel_{2}^2$ is given by
\begin{eqnarray}
\label{3.36}
\begin{split}
J_1'(q)[\delta q]=\sum_{i=1}^2\sum_{j=1}^{N}\int_0^T\int_\Omega\delta q(x)\cdot u_{j}^i(x,t)\varpi_{j}^i(x,t)dxdt.
\end{split}
\end{eqnarray}
Notice that $J_1(q)$ can be written as the Lagrangian multiplier formula,
\begin{eqnarray*}
\label{}
\begin{split}
J_1(q)=\frac{1}{2}\parallel\mathcal{F}(q)-\Phi^{\delta}\parallel_{2}^2+\sum_{i=1}^2\sum_{j=1}^N\left( _{0}D_t^{\alpha}u_j^i-\sum\limits_{k,l=1}^{d}\frac{\partial}{\partial x_k}\left(a_{kl}(x)\frac{\partial u_j^i}{\partial x_l}\right)+q(x)u_j^i-\phi_jv_i,\ \lambda\right)_{L^2(D)},
\end{split}
\end{eqnarray*}
for any multiplier function $\lambda\in L^2(D)$.
Then, from Theorem 3.4.2 and integration by parts, we have that
\begin{align*}
\label{}
J_1&(q+\delta q)-J_1(q)=\frac{1}{2}\left(\left\|\mathcal{F}(q+\delta q)-\Phi^{\delta}\right\|_2^2-\left\|\mathcal{F}(q)-\Phi^{\delta}\right\|_2^2\right)\\
&\ +\left( _{0}D_t^{\alpha}u_j^i(q+\delta q)-\sum\limits_{k,l=1}^{d}\frac{\partial}{\partial x_k}\left(a_{kl}(x)\frac{\partial u_j^i(q+\delta q)}{\partial x_l}\right)+(q+\delta q)u_j^i(q+\delta q)-\phi_jv_i,\ \lambda\right)_{L^2(D)}\\
&\ \ \ -\left( _{0}D_t^{\alpha}u_j^i(q)-\sum\limits_{k,l=1}^{d}\frac{\partial}{\partial x_k}\left(a_{kl}(x)\frac{\partial u_j^i(q)}{\partial x_l}\right)+qu_j^i-\phi_jv_i,\ \lambda\right)_{L^2(D)}\\
=&\sum_{i=1}^2\sum_{j=1}^N\bigg[\left(\int_{\Lambda\times(0,T)}\frac{\partial u_{j}^i(q)}{\partial\nu_A}h(x,t)dxdt-\varphi_{ij}^\delta\right)\cdot\int_{\Lambda\times(0,T)}\frac{\partial \vartheta_{j}^i}{\partial\nu_A}h(x,t)dxdt\\
&\ \ \ +\left( _{0}D_t^{\alpha}\vartheta_j^i-\sum\limits_{k,l=1}^{d}\frac{\partial}{\partial x_k}\left(a_{kl}(x)\frac{\partial \vartheta_{j}^i}{\partial x_l}\right)+q\vartheta_{j}^i+\delta q\cdot u_j^i(q),\ \lambda\right)_{L^2(D)}\bigg]+o(\|\delta q\|_{L^\infty(\Omega)})\\
=&\sum_{i=1}^2\sum_{j=1}^N\bigg[\left(\int_{\Lambda\times(0,T)}\frac{\partial u_{j}^i(q)}{\partial\nu_A}h(x,t)dxdt-\varphi_{ij}^\delta\right)\cdot\int_{\Lambda\times(0,T)}\frac{\partial \vartheta_{j}^i}{\partial\nu_A}h(x,t)dxdt\\
&\ \ \ +\left( _{t}D_T^{\alpha}\lambda-\sum\limits_{k,l=1}^{d}\frac{\partial}{\partial x_k}\left(a_{kl}(x)\frac{\partial \lambda}{\partial x_l}\right)+q\lambda,\ \vartheta_j^i\right)_{L^2(D)}+\left(\delta q\cdot u_j^i(q),\ \lambda\right)_{L^2(D)}\\
&\ \ \ -\int_{\partial\Omega\times(0,T)}\frac{\partial \vartheta_{j}^i}{\partial\nu_A}\lambda(x,t)dxdt\bigg]+o(\|\delta q\|_{L^\infty(\Omega)}).
\end{align*}
By choosing $\lambda=\varpi_{j}^i$, and $\varpi_{j}^i$ satisfies (\ref{3.35}), then we obtain (\ref{3.36}).

\subsection{The conjugate gradient method}
The conjugate gradient method (CGM) combined with an appropriate stopping rule can serve as a regularization method, and it has been applied to various inverse problems \cite{L.sun2017,T. Wei2016, jin2007}. However, the CGM is a deterministic regularization method which yield only a point estimate of the solution, without quantifying the associated uncertainties of measurement data noise. Here, we use CGM to calculate a maximum a posteriori (MAP) estimator of the solution, and then in conjunction with the Laplace approximations to sample the associated posterior distribution. The CGM we utilized is presented in Algorithm 1.

\begin{algorithm}[H]
\caption{The conjugate gradient method for solving the variational problem.}
      1:~~Choose $q_0$, and set $k=0$;\\
      ~2:~~Solve the direct problem $(\ref{mix-eq1})$ with $q=q_k$, $\phi=\phi_j$, $v=v_i$, and determine the residual
\begin{eqnarray*}
\label{}
\begin{split}
r_{ij}=\mathcal{F}(q_k)-\varphi_{ij}^\delta;
\end{split}
\end{eqnarray*}
      3:~~Solve the adjoint equation $(\ref{3.35})$ and determine the gradient $J'(q_k)$ by $(\ref{3.36})$;\\
      ~4:~~Calculate the conjugate coefficient $\gamma_k$ by
\begin{eqnarray*}
\label{Qresidual-eq1}
\begin{split}
\gamma_k=\frac{\|J'(q_k)\|_{L^2(\Omega)}^2}{\|J'(q_{k-1})\|_{L^2(\Omega)}^2},\ \gamma_0=0,
\end{split}
\end{eqnarray*}
       and the decent direction $d_k$ by
\begin{eqnarray*}
\label{Qresidual-eq2}
\begin{split}
d_k=-J'(q_k)+\gamma_kd_{k-1},\ d_0=-J'(q_0);
\end{split}
\end{eqnarray*}
      5:~~Solve the sensitivity equation $(\ref{3.31})$ for $\vartheta_j^i(q_k)$ with $\delta q=d_k$;\\
      ~6:~~Update the coefficient $q_k$ by
\begin{eqnarray*}
\label{Qresidual-eq3}
\begin{split}
q_{k+1}=q_k+\beta_kd_k,\ k=0,1,2,\cdots,
\end{split}
\end{eqnarray*}
where $\beta_k$ is given by
\begin{eqnarray*}
\label{Qresidual-eq}
\begin{split}
\beta_k=&-\frac{\sum_{i=1}^2\sum_{j=1}^{N}\left(\mathcal{F}(q_k)-\varphi_{ij}^\delta\right)\int_{\Lambda\times(0,T)}\frac{\partial\vartheta_{j}^i(x,t)}{\partial \nu_A}h(x,t)dxdt+\mu\int_{\Omega}q_kd_kdx}{\sum_{i=1}^2\sum_{j=1}^{N}\left(\int_{\Lambda\times(0,T)}\frac{\partial\vartheta_{j}^i(x,t)}{\partial \nu_A}h(x,t)dxdt\right)^2+\mu\int_{\Omega}d_k^2dx};
\end{split}
\end{eqnarray*}
      ~7:~~Increase k by one and go to step (2), repeat the above procedure until a stopping criterion \\ $~~~~~$is satisfied.\\
   \label{algorithm-pgddf}
\end{algorithm}

\section{Bayesian theory and Laplace approximation}

\subsection{Convergence of Laplace approximation}
From classical Bayesian theory, and notice that the observation error $\eta$ is an independent and identically distributed (i.i.d) Gauss random vector with mean zero and the covariance matrix  $B=\delta^{2}I$, here $I$ is the unit matrix, we can write the minimization functional $(\ref{3.29})$ as \cite{stuart2010}
\begin{eqnarray}
\label{Qresidual-eq}
\begin{split}
J(q)&\propto\frac{1}{2\delta^2}\parallel\mathcal{F}(q)-\Phi^\delta\parallel_2^2+\frac{\mu}{2\delta^2}\parallel q\parallel_{L^2(\Omega)}^2\\
&=\frac{1}{2}\parallel\mathcal{F}(q)-\Phi^\delta\parallel_B^2+\frac{\mu}{2\delta^2}\parallel q\parallel_{L^2(\Omega)}^2\\
&:= J_B(q),
\end{split}
\end{eqnarray}
where $\parallel\cdot\parallel_B$ is a covariance weighted norm on $R^{2N\times1}$ given by $\parallel\cdot\parallel_B=\parallel\delta^{-1}I\cdot\parallel_2$, and $\Phi^\delta=(\varphi_{11}^\delta,\cdots,\varphi_{1N}^\delta,\varphi_{21}^\delta,\cdots,\varphi_{2N}^\delta)^T\in R^{2N\times1}$ is a stretched vector version of the one defined in $(\ref{3.29})$. Furthermore, the minimizer of $J_B$ defines the maximum a posteriori estimator
\begin{eqnarray}
\label{eq20}
\begin{split}
q_{\mathrm{MAP}}=\arg\min\limits_{q\in Q} J_B(q).
\end{split}
\end{eqnarray}
The Laplace approximation (LA) in essence is a linearization around the MAP point $q_{\mathrm{MAP}}$ for sampling the posterior
distribution of the solution (refer to \cite{iglesias2010}). It consists of approximating the posterior measure(or distribution) by $\omega\thickapprox N(q_{\mathrm{MAP}},C_{\mathrm{MAP}})$, here $C_{\mathrm{MAP}}=(J_B''(q_{\mathrm{MAP}}))^{-1}$ is the inverse of Hessian of $J_B(q_{\mathrm{MAP}})$. Next, motivated by \cite{cotter2010}, we analyze the Hellinger distance between the exact posterior measure and its approximation given by $N(q_{\mathrm{MAP}},C_{\mathrm{MAP}})$ and obtain the convergence as $q\rightarrow q_{\mathrm{MAP}}$. In \cite{wacher2017}, Wacher gave a bound of the Hellinger distance between the posterior and its LA, but it seems not to deduce the convergence. Moreover, the estimation of other distances (such as Kullback-Leibler divergence) between the posterior distribution and its approximation were applied in analysing the stochastic surrogate models (e.g., \cite{yan2015, jiang2017}). \\
\\
\textbf{Lemma 5.1.1.} For every $\varepsilon >0$ and $q\in Q$, there exists $M\in R$ such that, the forward map $\mathcal{F}:Q\mapsto R^{2N\times1}$ satisfies,
\begin{eqnarray}
\label{Qresidual-eq}
\begin{split}
\parallel\mathcal{F}(q)\parallel_B\leq\exp(\varepsilon\parallel q\parallel_{L^{\infty}(\Omega)}^2+M).
\end{split}
\end{eqnarray}
\textbf{Proof}. From theorem 3.1.2, and trace theorem, it follows that
\begin{eqnarray*}
\label{Qresidual-eq}
\begin{split}
\mid\varphi_{ij}\mid&=\mid<\frac{\partial u_j^i}{\partial\nu_A},h>\mid\\
&\leq\parallel\frac{\partial u_j^i}{\partial\nu_A}\parallel_{L^2(0,T;H^{-\frac{1}{2}}(\Lambda)}\parallel h\parallel_{L^2(0,T;H^{\frac{1}{2}}(\Lambda))}\\
&\leq C\parallel u_j^i\parallel_{B^{\frac{\alpha}{2}}(\Omega\times(0,T))}\parallel h\parallel_{L^2(0,T;H^{\frac{1}{2}}(\Lambda))}\\
&\leq C\parallel\phi_j(x)v_i(t)\parallel_{H^1(0,T;L^2(\Omega))}\parallel h\parallel_{L^2(0,T;H^{\frac{1}{2}}(\Lambda))}\\
&\leq C<\infty.
\end{split}
\end{eqnarray*}
Thus
\begin{eqnarray*}
\label{Qresidual-eq}
\begin{split}
\parallel\mathcal{F}(q)\parallel_B&=\delta^{-1}(\sum\limits_{i=1}^2\sum\limits_{j=1}^N\mid\varphi_{ij}\mid^2)^{\frac{1}{2}}
&\leq\frac{2NC}{\delta}
&:=C_{N,\delta}<\infty,
\end{split}
\end{eqnarray*}
and then
\begin{eqnarray*}
\label{Qresidual-eq}
\begin{split}
\parallel\mathcal{F}(q)\parallel_B&\leq e^{\log C_{N,\delta}}e^{\varepsilon\parallel q\parallel_{L^{\infty}(\Omega)}^2}\\
&\leq e^{(\varepsilon\parallel q\parallel_{L^{\infty}(\Omega)}^2+M)},
\end{split}
\end{eqnarray*}
here, $M=\log C_{N,\delta}$.

From the Laplace approximation, the Taylor expansion of $J_B(q)$ at MAP point $q_{\mathrm{MAP}}$ is given by
\begin{eqnarray}
\label{eq16}
\begin{split}
J_B(q)&=J_B(q_{\mathrm{MAP}})+\frac{1}{2}J_B''(q_{\mathrm{MAP}})(q-q_{\mathrm{MAP}})^2
+o(\parallel q-q_{\mathrm{MAP}}\parallel_{L^{\infty}(\Omega)}^2)\\
:&=\tilde{J_B}(q)+o(\parallel q-q_{\mathrm{MAP}}\parallel_{L^{\infty}(\Omega)}^2).
\end{split}
\end{eqnarray}

In the infinite dimensional version of Beyesian theory. The posterior measure is absolutely continuous with respect to the prior, and the Rodon-Nikodym derivative between them is determined by the data likelihood, i.e.
\begin{eqnarray}
\label{Qresidual-eq}
\begin{split}
\frac{dw}{dw_0}(q)=\frac{1}{Z}e^{-J_B(q)};\ \ Z=\int_{L^{\infty}(\Omega)}e^{-J_B(q)}dw_0(q).
\end{split}
\end{eqnarray}
Similarly, the Laplace approximation version of Bayesian formula is defined as
\begin{eqnarray}
\label{Qresidual-eq}
\begin{split}
\frac{d\tilde{w}}{dw_0}(q)=\frac{1}{\tilde{Z}}e^{-\tilde{J}_B(q)};\ \ \tilde{Z}=\int_{L^{\infty}(\Omega)}e^{-\tilde{J}_B(q)}dw_0(q).
\end{split}
\end{eqnarray}\\
\\
\textbf{Lemma 5.1.2.} The forward map $\mathcal{F}:Q\rightarrow R^{2N\times1}$ satisfies:\\
(i) for every $\varepsilon>0$ and $r>0$, there exsits $M=M(\varepsilon,r)\in R$ such that, for all $q\in Q$ and $\parallel\Phi^\delta\parallel_B<r$,
\begin{eqnarray}
\label{Qresidual-eq}
\begin{split}
\frac{1}{2}\parallel\mathcal{F}(q)-\Phi^\delta\parallel_B^2\geq M-\varepsilon\parallel q\parallel_{L^{\infty}(\Omega)}^2;
\end{split}
\end{eqnarray}
(ii) for every $r>0$, there exsits a $L=L(r)>0$ such that for all $q\in Q$ and $\Phi^\delta\in R^{2N\times1}$ with $\max\{\parallel q\parallel_{L^{\infty}(\Omega)},\parallel\Phi^\delta\parallel_B\}<r$,
\begin{eqnarray}
\label{Qresidual-eq}
\begin{split}
\frac{1}{2}\parallel\mathcal{F}(q)-\Phi^\delta\parallel_B^2\leq L(r).
\end{split}
\end{eqnarray}

The Lemma 5.1.2 is the direct results of Lemma 5.1.1, by using Lemma 2.1 in \cite{cotter2010}, so we omit the proof.\\
\\
\textbf{Theorem 5.1.3.} Assume that $L(r)$ and $M$ are defined in Lemma 5.1.2, then the measure $w$ and its Laplace approximation measure $\tilde{w}$ are close with respect to the Hellinger distance, i.e., there is a constant $C$, such that
\begin{eqnarray}
\label{5.53}
\begin{split}
d_{Hell}(w,\tilde{w})\leq C\sqrt{e^{L(r)+\frac{\mu}{2\delta^2}q_{\mathrm{max}}^2\mid\Omega\mid}+e^{-M-\frac{\mu}{\delta^2}q_{\mathrm{min}}\mid \Omega\mid}}\cdot o(\parallel q-q_{\mathrm{MAP}}\parallel_{L^{\infty}(\Omega)}^2),
\end{split}
\end{eqnarray}
and
\begin{eqnarray*}
\label{}
\begin{split}
d_{Hell}(w,\tilde{w})\rightarrow0,\ \ \ \text{as}\ \ q\rightarrow q_{\mathrm{MAP}},\ \ \text{in}\ L^\infty(\Omega),
\end{split}
\end{eqnarray*}
where the Hellinger distance is defined by
\begin{eqnarray*}
\label{}
\begin{split}
d_{Hell}(w,\tilde{w})=\left(\frac{1}{2}\int_{L^{\infty}(\Omega)}\left(\sqrt{\frac{dw}{dw_0}}-\sqrt{\frac{d\tilde{w}}{dw_0}}\right)^2dw_0\right)^{\frac{1}{2}}.
\end{split}
\end{eqnarray*}
\\
\textbf{Proof.} By using Lemma 5.1.2 (ii), we obtain
\begin{eqnarray}
\label{eq18}
\begin{split}
\mid Z\mid &=\int_{L^{\infty}(\Omega)}e^{-\frac{1}{2}\parallel\mathcal{F}(q)-\Phi^\delta\parallel_B^2-\frac{\mu}{2\delta^2}\parallel q\parallel_{L^{2}(\Omega)}^2} dw_0(q)\\
&\geq\int_{L^{\infty}(\Omega)}e^{-L(r)-\frac{\mu}{2\delta^2}q_{\mathrm{max}}^2\mid \Omega\mid}dw_0(q)\\
&\geq e^{-L(r)-\frac{\mu}{2\delta^2}q_{\mathrm{max}}^2\mid \Omega\mid}\cdot\mid Q_r\mid_\infty\\
&>0,
\end{split}
\end{eqnarray}
where $\mid \Omega\mid=\int_{\Omega}1dx$, $\mid Q_r\mid_\infty=\int_{L^{\infty}(\Omega)}1dw_0(q)$, and
\begin{eqnarray}
\label{eq17}
\begin{split}
\mid \tilde{Z}\mid&=\int_{L^{\infty}(\Omega)}e^{-\tilde{J}_{B(q)}} dw_0(q)\\
&=\int_{L^{\infty}(\Omega)}e^{-\frac{1}{2}\parallel\mathcal{F}(q_{\mathrm{MAP}})-\Phi^{\delta}\parallel_B^2-\frac{\mu}{2\delta^2}\parallel q_{\mathrm{MAP}}\parallel_{L^{2}(\Omega)}^2-\frac{1}{2}J''_B(q_{\mathrm{MAP}})(q-q_{\mathrm{MAP}})^2} dw_0(q)\\
&=e^{-\frac{1}{2}\parallel\mathcal{F}(q_{\mathrm{MAP}})-\Phi^{\delta}\parallel_B^2-\frac{\mu}{2\delta^2}\parallel q_{\mathrm{MAP}}\parallel_{L^{2}(\Omega)}^2}\int_{L^{\infty}(\Omega)}e^{-\frac{1}{2}J''_B(q_{\mathrm{MAP}})(q-q_{\mathrm{MAP}})^2}dw_0(q).
\end{split}
\end{eqnarray}
From the second-order necessary condition for minimizer, lemma 5.1.2 (i) and Fernique theorem (see \cite{cotter2010}), we see that
\begin{eqnarray*}
\label{Qresidual-eq}
\begin{split}
\mid\tilde{Z}\mid&=C\int_{L^{\infty}(\Omega)}e^{-\frac{1}{2}\parallel\mathcal{F}(q_{\mathrm{MAP}})-\Phi^{\delta}\parallel_B^2-\frac{\mu}{2\delta^2}\parallel q_{\mathrm{MAP}}\parallel^2_{L^2(\Omega)}}dw_0(q)\\&\leq C\int_{L^\infty(\Omega)}e^{\varepsilon\parallel q_{\mathrm{MAP}}\parallel_{L^\infty(\Omega)}^2-M-\frac{\mu}{2\delta^2}q_{\mathrm{min}}^2\mid \Omega\mid}dw_0(q)\\
&\leq Ce^{-M-\frac{\mu}{2\delta^2} q_{\mathrm{min}}^2}\mid Qr\mid_{\infty}.
\end{split}
\end{eqnarray*}
Moreover, notice $(\ref{eq16})$, it follows that
\begin{eqnarray}
\label{Qresidual-eq}
\begin{split}
\mid Z-\tilde{Z}\mid&\leq\int_{L^{\infty}(\Omega)}\mid e^{-J_B(q)}-e^{-\tilde{J}_{B(q)}}\mid dw_0(q)\\
&\leq\int_{L^{\infty}(\Omega)}\mid J_B(q)-\tilde{J_B}(q)\mid dw_0(q)\\
&\leq o(\parallel q-q_{\mathrm{MAP}}\parallel_{L^{\infty}(\Omega)}^2)\int_{L^{\infty}(\Omega)} 1dw_0(q)\\
&=o(\parallel q-q_{\mathrm{MAP}}\parallel_{L^{\infty}(\Omega)}^2)\mid Q_r\mid_\infty.
\end{split}
\end{eqnarray}
From the definition of Hellinger distance, we get
\begin{eqnarray*}
\label{Qresidual-eq}
\begin{split}
2d_{Hell}(w,\tilde{w})^2&=\int_{L^{\infty}(\Omega)}( Z^{-\frac{1}{2}}e^{-\frac{1}{2}{J}_B(q)}-\tilde{Z}^{-\frac{1}{2}}e^{-\frac{1}{2}\tilde{J}_{B}(q)})^2 dw_0(q)\\
&=\int_{L^{\infty}(\Omega)}[( Z^{-\frac{1}{2}}e^{-\frac{1}{2}{J}_{B}(q)}-Z^{-\frac{1}{2}}e^{-{\frac{1}{2}\tilde{J}_{B}(q)}})+({Z}^{-\frac{1}{2}}e^{-\frac{1}{2}\tilde{J}_{B}(q)}-\tilde{Z}^{-\frac{1}{2}}e^{-\frac{1}{2}\tilde{J}_{B}(q)})]^2dw_0(q)\\
&\leq \int_{L^{\infty}(\Omega)}\frac{2}{Z}(e^{-\frac{1}{2}J_B(q)}-e^{-\frac{1}{2}\tilde{J}_B(q)})^2dw_0(q)+\int_{L^{\infty}(\Omega)}2(Z^{-\frac{1}{2}}-\tilde{Z}^{-\frac{1}{2}})^2e^{-\tilde{J}_B(q)}dw_0(q)\\
&:=I_1+I_2.
\end{split}
\end{eqnarray*}
Now, using $(\ref{eq16})$ and $(\ref{eq17})$, it deduces
\begin{eqnarray}
\label{Qresidual-eq}
\begin{split}
I_{1}&\leq\frac{1}{2Z}\int_{L^{\infty}(\Omega)}\mid J_B(q)-\tilde{J_B}(q)\mid^2 dw_0(q)\\
&\leq\frac{1}{2Z}\cdot o(\parallel q-q_{\mathrm{MAP}}\parallel_{L^\infty(\Omega)}^4)\mid Q_r\mid_\infty\\
&\leq \frac{1}{2}e^{L(r)+\frac{\mu}{2\delta^2}q_{\mathrm{max}}^2\mid \Omega\mid}\mid Q_r\mid_\infty\cdot o(\parallel q-q_{\mathrm{MAP}}\parallel_{L^\infty(\Omega)}^4),
\end{split}
\end{eqnarray}

and then
\begin{eqnarray}
\label{Qresidual-eq}
\begin{split}
I_2&\leq 2Ce^{-M-\frac{\mu}{2\delta^2}q_{\mathrm{min}}^2\mid \Omega\mid}\mid Z^{-\frac{1}{2}}-\tilde{Z}^{-\frac{1}{2}}\mid^2\\
&\leq Ce^{-M-\frac{\mu}{2\delta^2}q_{\mathrm{min}}^2\mid \Omega\mid}\mid \xi^{-3}\mid\mid Z-\tilde{Z}\mid^2\\
&\leq C\max\{Z^{-3},\tilde{Z}^{-3}\}e^{-M-\frac{\mu}{2\delta^2}q_{\mathrm{min}}^2\mid \Omega\mid}\cdot o(\parallel q-q_{\mathrm{MAP}}\parallel_{L^{\infty}(\Omega)}^4)\mid Q_r\mid_\infty^2.
\end{split}
\end{eqnarray}
Hence, we can find that
\begin{eqnarray}
\label{Qresidual-eq}
\begin{split}
d_{Hell}(w,\tilde{w})=C\sqrt{e^{-M-\frac{\mu}{2\delta^2}q_{\mathrm{min}}^2\mid \Omega\mid}+e^{L(r)+\frac{\mu}{2\delta^2}q_{\mathrm{max}}^2\mid \Omega\mid}}\cdot o(\parallel q-q_{\mathrm{MAP}}\parallel_{L^{\infty}(\Omega)}^2).
\end{split}
\end{eqnarray}
\\
\textbf{Remark 5.1.4.} From the proof of Theorem 5.1.3, we can see that the approximation error is derived from the Taylor expansion (\ref{eq16}), and second-order convergence rate is
given. Moreover, from the convergence estimation (\ref{5.53}), we also find that when the noise level $\delta$ is decreased or the regularization parameter $\mu$ is increased, the Hellinger distance $d_{Hell}(w,\tilde{w})$ become small, and it provides some inspiration for selecting the regularization parameter.\\

\subsection{Numerical algorithm of Laplace approximation}
For effective numerical simulation, we take a finite dimensional approximation of the previous minimization problem as follows
\begin{eqnarray}
\label{}
\begin{split}
J_B^M(q^M)=\frac{1}{2}\parallel \mathcal{\tilde{F}}(q)-\Phi^\delta\parallel_B^2+\frac{1}{2}\parallel q^M\parallel_{B_\mu}^2,
\end{split}
\end{eqnarray}
where $q^M=(q_1,...,q_M),\ B=\delta^2I,\ B_\mu=\frac{\delta^2}{\mu}I$, and $\mathcal{\tilde{F}}$ is a finite-dimensional approximation of the continuous forward map $\mathcal{F}$. The Laplace approximation theory shows the posterior distribution $w\approx N(q_{\mathrm{MAP}},C_{\mathrm{MAP}})$. In the finite dimensional, the covariance matrix (see \cite{Reynolds2008}, Section 10.5)
\begin{eqnarray}
\label{5.65}
\begin{split}
C_{\mathrm{MAP}}=({J^M_B}''(q^M))^{-1}&=(B_\mu^{-1}+P^{T}B^{-1}P)^{-1}=(\frac{\mu}{\delta^2}I+\frac{1}{\delta^2}P^{T}P)^{-1},
\end{split}
\end{eqnarray}
where $P$ is Jacobian matrix of forward operator $\tilde{\mathcal{F}}$ at $q^M$ point. Notice that the covariance formula (\ref{5.65}) only uses the first order derivatives of $\tilde{\mathcal{F}}$. A standard implementation of Laplace approximation is presented in the following algorithm \cite{Reynolds2008}:
\begin{algorithm}
\caption{Laplace approximation (LA) for sampling.}
      1:~~Compute $q_{\mathrm{MAP}}$ from $(\ref{eq20})$ by using Algorithm 1 (CGM), and $C_{\mathrm{MAP}}$ from $(\ref{5.65})$, respectively;\\
      ~2:~~Compute the Cholesky factor L of $C_{\mathrm{MAP}}$, i.e.,
      \begin{eqnarray}
      C_{\mathrm{MAP}}=LL^{T};
      \end{eqnarray}
      ~3:~~For $j=\{1,...,N_e\}$, generate
      \begin{eqnarray}\label{5.67}
      q^{j}=q_{\mathrm{MAP}}+L^{T}z^{j},
      \end{eqnarray}
      $~~~~~$where $z^{j}\sim N(0,I)$.\\
   \label{algorithm-pgddf}
\end{algorithm}

Samples generated by (\ref{5.67}) are drawn from $N(q_{\mathrm{MAP}},C_{\mathrm{MAP}})$, and so the ensemble of $N_e$ realizations $\{q^j\}_{j=1}^{N_e}$ provides
an approximation to $N(q_{\mathrm{MAP}},C_{\mathrm{MAP}})$ and hence the posterior. Finally, we use the mean of the sample $\bar{q}_n=\frac{1}{n}\sum\limits_{j=1}^n q^j$ as an approximation of $q_{\mathrm{MAP}}$, where the convergence of $\bar{q}_n$ follows from the strong law of large numbers. From the classical Gaussian statistic theory, we find that $\bar{q}_n$ are consistent and best unbiased estimate of $q_{\mathrm{MAP}}$.

\subsection{Confidence region and skewness}
We can calculate the confidence region for the inferred parameters (such as degradation coefficient in IDCP) by the LA posterior samples. The confidence region is a set of points in an $d$-dimensional space, often represented as an ellipsoid around a point which is an estimated parameter. The confidence region shows that the real value of the identified parameter has a certain probability of falling around the numerical construction result and quantifies the level of confidence that the parameter lies in the region. Therefore, the confidence region gives the reliability of the construction result of the inferred parameters.\\
\\
\textbf{Definiton 4.3.1. (Confidence region)} \cite{draper1981} The confidence region of $(1-\alpha)\times100\%$ is defined as follows
\begin{eqnarray}
\label{eq19}
\begin{split}
\{\bar{q}_n:n(\bar{q}_n-q_{\mathrm{MAP}})^{T}C_{\mathrm{MAP}}^{-1}(\bar{q}_n-q_{\mathrm{MAP}})\leq\mathrm{\chi^2_\alpha}(M)\}.
\end{split}
\end{eqnarray}
It is an ellipsoid centered on $q_{\mathrm{MAP}}$ point, and assume the eigenvalue of $C_{\mathrm{MAP}}$ are
\begin{eqnarray}
\label{eq19}
\begin{split}
\lambda_1\geq\lambda_2\geq ...\geq\lambda_{M}>0.
\end{split}
\end{eqnarray}
We find that the axes of the ellipsoid are $\sqrt{\frac{\lambda_j\mathrm{\chi^2_\alpha}(M)}{n}}$. Since $C_{\mathrm{MAP}}=\delta^2(\mu I+P^{T}P)^{-1}$, it deduces that $\lambda_j(C_{\mathrm{MAP}})=\frac{\delta^2}{\mu+\lambda_j(P^{T}P)}$, and the length of each axis is $\delta\sqrt{\frac{\chi_\alpha^2(M)}{n(\mu+\lambda_j(P^{T}P))}}$.
Let $\bar{q}_{n,j}$ be the component of $\bar{q}_n$, then it satisfies $\mid \bar{q}_{n,j}-\bar{q}_{\mathrm{MAP},j}\mid\leq\delta\sqrt{\frac{\chi_\alpha^2(M)}{n(\mu+\lambda_j(P^{T}P))}}$.
Thus, we can see that the confidence region is a multi-dimensional generalization of a confidence interval. It also can be seen that when the number of samples $n$, the regularization parameter $\mu$ and the eigenvalue of $P^TP$ increase or the measurement noise level $\delta$ decrease, the size of confidence region will decrease and give higher reliability.
\\\indent

Skewness is a statistic that studies the symmetry of data distribution. By measuring the skewness, we can determine the degree and direction of the asymmetry of the data distribution. The definition of skewness is given below\\
\\
\textbf{Definition 4.3.2. (Skewness)} \cite{Everitt2006} Assume the third central moment of the random variable $X$ is exists, the following ratios:
\begin{eqnarray}\label{ske}
\begin{split}
\beta=E\left[ \left(\frac{X-\varsigma}{\sigma}\right)^3\right]=\frac{\gamma_3}{\sigma^3}
\end{split}
\end{eqnarray}
 is called  the skewness  coefficient for X, i.e., skewness. Here $\varsigma$ is the mean, $\sigma$ is the standard deviation, $E$ is the expectation operator, $\gamma_3$ is the third central moment.

By normalization we can transform the identified parameter into a well defined probability density function for a hypothetical random variable, and calculate the skewness of
random variable to estimate the symmetry of the inferred parameter (see Section 5).

\section{Numerical experiments and discussions}
\label{Numerical-result-FD}

In this section, we present some numerical examples to illustrate the feasibility of the Laplace approximation method for IDCP.

The noisy data are generated by adding random perturbations as follows \cite{jin2007}
\begin{eqnarray}
\label{eqq1}
\begin{split}
\varphi_{ij}^\delta=\varphi_{ij}+\max\limits_{i,j}\{\mid\varphi_{ij}\mid\}\varepsilon\zeta,\ \ i=1,2;\ j=1,2,\cdots,N,\\
\end{split}
\end{eqnarray}
here $\zeta$ is a Gaussian random variable with zero mean and unit standard deviation. $\varepsilon$ indicates the noise level.

To show the accuracy of numerical solutions, we compute the relative $L^2$ error denoted by
\begin{eqnarray}
\label{eqq2}
\begin{split}
r_{e}=\frac{\parallel \bar{q}_n(x)-q(x)\parallel_{L^2(\Omega)}}{\parallel q(x) \parallel_{L^2(\Omega)}},
\end{split}
\end{eqnarray}
where $\bar{q}_n(x)$ is the approximate degradation coefficient reconstructed by the LA algorithm, and $q(x)$ is the exact solution.\\\indent
In an iteration algorithm, aim to make $q_k$ to approximate $q_{\mathrm{MAP}}$ point, the important work is to find a suitable stopping rule. Using $E_k$ to represent the residuals of iterations of $k-1$ and $k$ steps, i.e.,
\begin{eqnarray}
\label{eqq2}
\begin{split}
E_k=\parallel q_k-q_{k-1}\parallel_{L^{\infty}(\Omega)},
\end{split}
\end{eqnarray}
when $E_k\leq eps$, we stop iterating as a stopping criterion.

\subsection{Inversion for one-dimensional degradation coefficient}

The domain $\Omega$ under consideration is a unite line $[0,1]$, and the boundary is indicated as $\Lambda_0=\{x=0\}$, $\Lambda_1=\{x=1\}$. The direct problem $(\ref{mix-eq1})$ is discretized by using $300$ uniform rectangular finite element. Set $T=1$, and the grid point on $[0,T]$ is $101$. The boundary data are obtained by solving the direct problem $(\ref{mix-eq1})$ with $a_{kl}(x)=1$, and setting positive function $h(x,t)=(1-t)$.  We solve the direct problem, sensitive problem and adjoint problem by using the finite element method and construction of data $\varphi_{ij}^\delta$ using difference method. \\\indent
The $q_{k}$ obtained by conjugate gradient method is an approximation to $q_{\mathrm{MAP}}$. Samples generated by algorithm 2 are drawn from $N(q_{\mathrm{MAP}},\ C_{\mathrm{MAP}})$, and so the ensemble $\{q^{(j)}\}_{j=1}^{Ne}$, the length $Ne$ is taken to be $10000$, provides an approximation to $N(q_{\mathrm{MAP}},\ C_{\mathrm{MAP}})$.

\subsubsection{Smooth solution}
$\mathbf{Example~~ 1.}$ In $(\ref{mix-eq1})$, the accessible boundary $\Lambda$ is taken to be $\Lambda_1$, and the degradation coefficient is given by
 \begin{eqnarray*}
\label{eqq3}
\begin{split}
q(x)=x^2(1-x^2),\ \  x\in\Omega.
\end{split}
\end{eqnarray*}\indent

\begin{figure}[H]
\centering
\subfigure[$\alpha=0.3$ ]{
    \label{fig:subfig:a}
  \includegraphics[width=3in, height=2.5in]{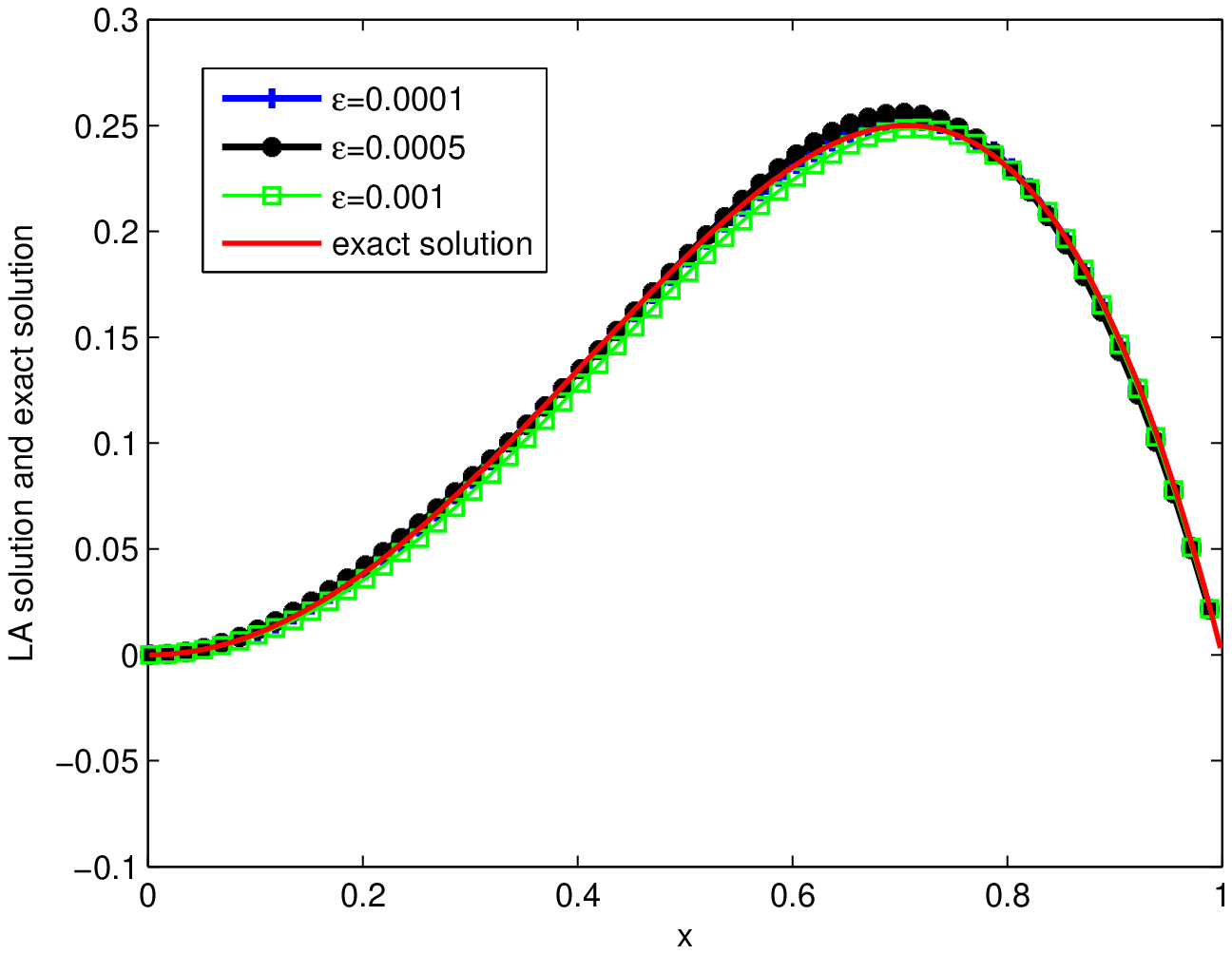}}
  \subfigure[$\alpha=0.7$]{
    \label{fig:subfig:b}
  \includegraphics[width=3in, height=2.5in]{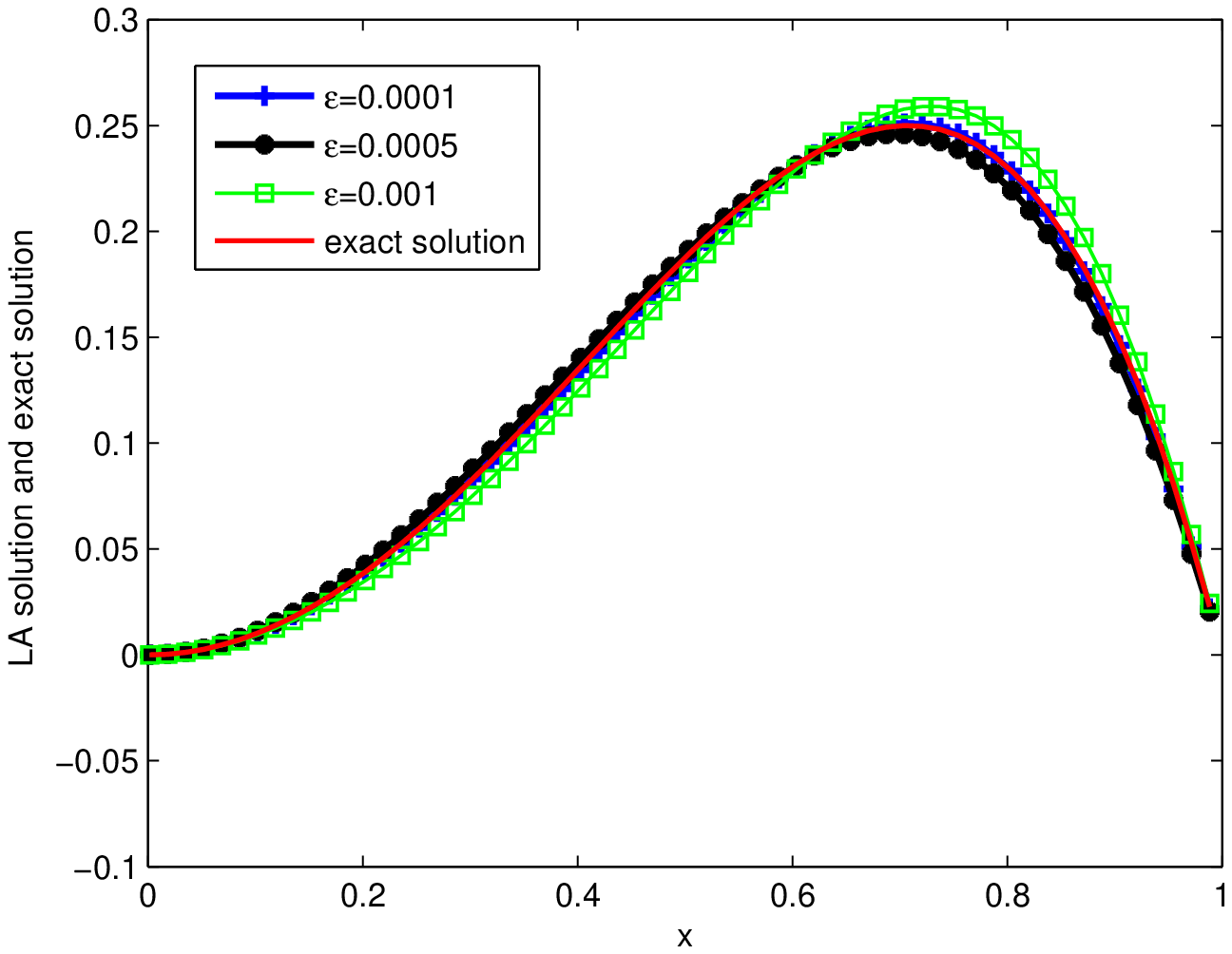}}\\
  \caption{The numerical result for Example 1 for various noise levels with
  $\mu=\delta^{\frac{3}{2}}$}
  \label{Example.1 for various noise le}
\end{figure}

First we investigate the effect of the amounts of basis functions. Trigonometric basis functions are used, i.e., $\phi_j\in span\{1, cos(2\pi x), sin(2\pi x),\cdots, cos(2N\pi x), sin(2N\pi x)\}$ in Example 1. In Table $\ref{tab1-N}$, for fractional order $\alpha=0.3$, we list the $L^2$ error $r_{e}$ of LA solution and exact solution under different amounts of basis functions and different noise levels, $N$ is the number of basis functions. From the Table $\ref{tab1-N}$, we can see that the error becomes smaller at the same noise level as the number of basis function increases. When the number of basis functions is the same, the result become worse with the increase of noise levels. In the following calculations, we chose $N=5$.\\
\begin{table}[H]
\small
\centering
\caption{ Numerical results for Example 1 with various N and $\varepsilon$ ($\alpha=0.3$).}
\begin{tabular}{c|c|c|c}
\Xhline{1pt}
 $N\setminus\varepsilon$   & $0.0001$ & $0.0005$ & $0.001$ \\
  \hline
 $1$ & $0.0525$ & $0.0528$ &  $ 0.0590$ \\
 \hline
 $2 $& $ 0.0444$&  $ 0.0477$ &  $ 0.0506$ \\
 \hline
 $3$ & $0.0111$& $0.0287$ &  $ 0.0409$ \\
 \hline
 $4$ & $ 0.0095$&  $ 0.0257$ &  $ 0.0344$ \\
 \hline
 $5$ & $0.0090$&  $0.0197$ &  $ 0.0295$ \\
\Xhline{1pt}
\end{tabular}
\label{tab1-N}
\end{table}
\begin{table}[H]
\small
\centering
\caption{ Numerical results for Example 1 with various $\mu$ and $\varepsilon$ ($\alpha=0.3$).}
\begin{tabular}{c|c|c|c}
\Xhline{1pt}
 $\mu\setminus\varepsilon$   & $0.0001$ & $0.0005$ & $0.001$ \\
  \hline
 $\delta$ & $0.1086$ & $0.2090$ &  $ 0.2970$   \\
 \hline
 $\delta^{\frac{3}{2}} $& $ 0.0090$&  $ 0.0197$ &  $ 0.0295$ \\
 \hline
 $\delta^{2}$ & $0.0386$& $0.0608$ &  $ 0.1437$ \\
\Xhline{1pt}
\end{tabular}
\label{tab2-re parameter}
\end{table}
\begin{table}[H]
\small
\centering
\caption{ Numerical results for Example 1 with various $\alpha$ and $\varepsilon$.}
\begin{tabular}{c|c|c|c}
\Xhline{1pt}
 $\alpha\setminus\varepsilon$   & $0.0001$ & $0.0005$ & $0.001$ \\
  \hline
 $0.1$ & $0.0086$ & $0.0117$ &  $ 0.0252$ \\
 \hline
 $0.3 $& $ 0.0090$&  $ 0.0197$ &  $ 0.0296$ \\
 \hline
 $0.5$ & $0.0079$& $0.0204$ &  $ 0.0398$   \\
 \hline
 $0.7$ & $ 0.0074$&  $ 0.0351$ &  $ 0.0517$ \\
 \hline
 $0.9$ & $0.0071$&  $0.0447$ &  $ 0.0642$ \\
\Xhline{1pt}
\end{tabular}
\label{tab1- fractional order}
\end{table}
In Theorem 3.3.4, we introduce a slightly crude rule of regularization parameter selection, i.e.,
\begin{eqnarray}
\label{eqq8}
\begin{split}
\frac{\delta^2}{\mu(\delta)}\rightarrow 0,\ as\ \delta\rightarrow 0.
\end{split}
\end{eqnarray}
The numerical results for Example 1 with various $\mu$ and $\varepsilon$ are shown in the Table $\ref{tab2-re parameter}$. We find that at the same noise level, $\mu=\delta$ was selected to meet the rule $(\ref{eqq8})$, but the results are unsatisfactory. By choosing $\mu=\delta^{\frac{3}{2}}$ which satisfies $(\ref{eqq8})$, the error between the LA solution and the exact solution is small, and the result was desirable. However, the crude rule $(\ref{eqq8})$ for regularization parameter selection is not satisfied with $\mu=\delta^2$, and the error between LA solution and exact solution will gradually increase and result will be shock. Thus, the regularization parameter can be chosen neither too large nor too small, even if it satisfies $(\ref{eqq8})$.
\\\indent
The numerical result for Example 1 for various noise levels $\varepsilon=0.0001,0.0005,0.001$, and in the case of $\alpha=0.3, 0.7$ are shown in Figure $\ref{Example.1 for various noise le}$. The basis function is the trigonometric basis function in $L^2(\Omega)$. Compare the Figure $\ref{Example.1 for various noise le}(a)$ with Figure $\ref{Example.1 for various noise le}(b)$, the results in Figure $\ref{Example.1 for various noise le}(b)$ is worse than  Figure $ \ref{Example.1 for various noise le}(a)$. In Table $\ref{tab1- fractional order}$, we further show the numerical errors $r_{e}$ of the Example 1 for different $\alpha$ and $\varepsilon$. It can be seen that the numerical results become worse as the noise levels increase.  The error become larger as fractional order $\alpha$ increases. However, when the noise level is 0.0001, the numerical result is insensitive to fractional order $\alpha$.\\\indent
Next, we change trigonometric basis function into polynomial basis function in $L^2(\Omega)$, i.e., $\phi_j\in span\{1, x, x^2 ,\cdots, x^N\}$. We show the reconstruction results for Example 1 under various error levels and different types of basis function with $\alpha=0.3$ in the Table $\ref{tab1-basis}$. The data indicate that the difference in types of basis function will affect the numerical results. The numerical results  obtained by trigonometric basis functions are better than those obtained by polynomial basis functions. In the following example, we only consider the triangle basis functions.\\

\begin{table}[H]
\small
\centering
\caption{ Numerical results for Example 1 for different basis functions and $\varepsilon$ with $\alpha=0.3$.}
\begin{tabular}{c|c|c|c}
\Xhline{1pt}
 Type of basis function$\setminus\varepsilon $   & $0.0001$ & $0.0005$ & $0.001$ \\
  \hline
 polynomial basis & $0.0129$ & $0.0279$ &  $ 0.0436$   \\
 \hline
 Trigonometric basis & $ 0.0090$&  $ 0.0197$ &  $ 0.0295$ \\
\Xhline{1pt}
\end{tabular}
\label{tab1-basis}
\end{table}

In the Figure 5.3
\ show that the results by using the direct flux data is not sensitive to the fractional order $\alpha$. From Table $\ref{tab1-data}$, we find that the reconstruction results by using the direct flux data are better than the ones by using the average flux data. However, we also get satisfactory numerical results by using the average flux data when the noise levels are $\varepsilon=0.0001,\ 0.0005,\ 0.001$ respectively. The reconstruction error caused by using average flux data increases sharply when the noise level exceeds $0.001$, but using the direct flux data still give good results when the noise level exceeds $0.005$ and even reaches $\varepsilon=0.05$. This is because, compared to direct flux data, the amount of the average flux data is less and provides limited information. Moreover, the limited measurement data lead to higher sensitivity to noise and severally ill-posedness of IDCP. But, the average flux data is rather easier to measure as a practical matter, and it has been widely used recently. Hence, recovering the degradation coefficient accurately by using limited measurement data can be a real challenge. \\\indent
\begin{figure}[H]
\centering
\subfigure[$\alpha=0.3$ ]{
    \label{fig:subfig:a}
  \includegraphics[width=3in, height=2.5in]{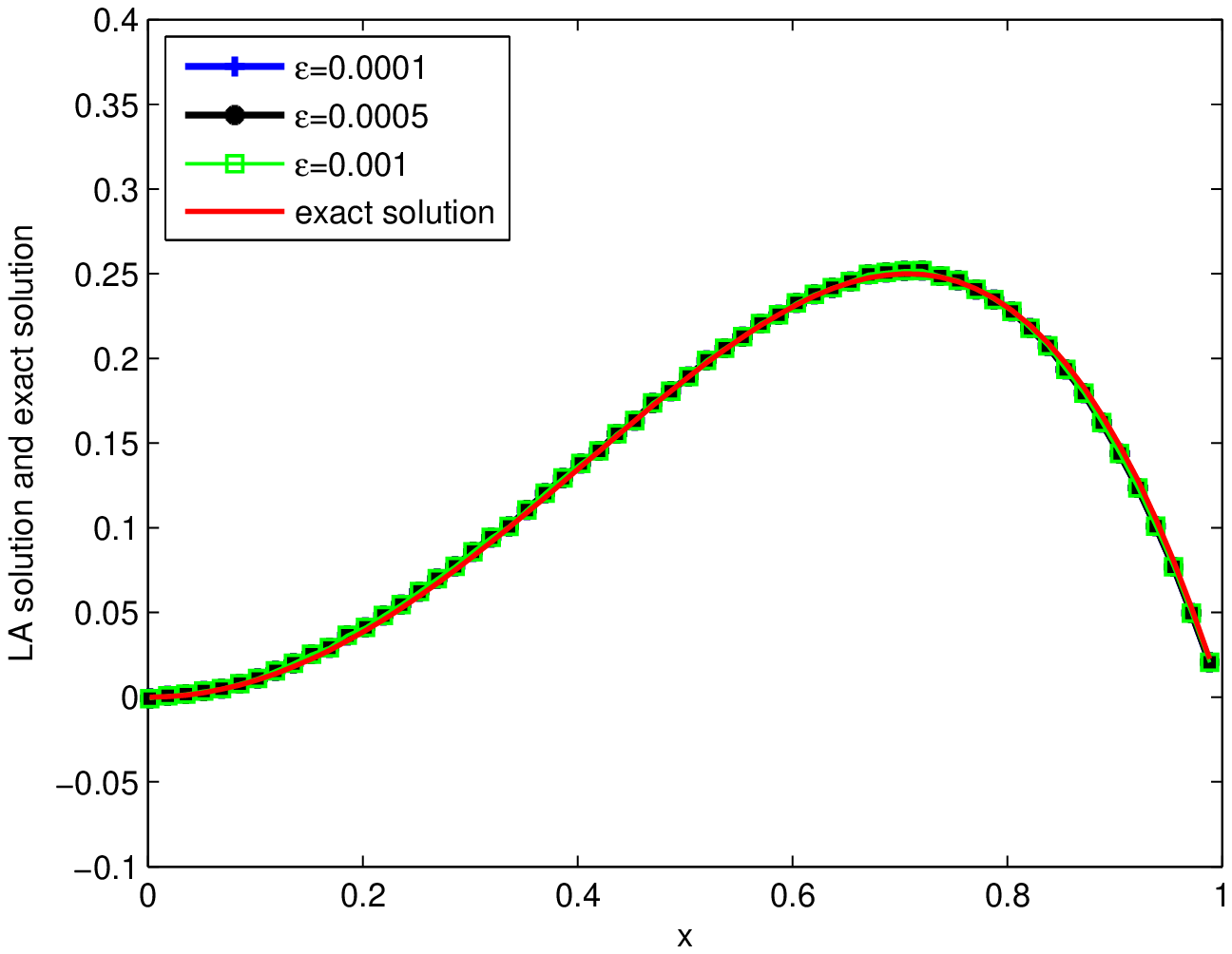}}
  \subfigure[$\alpha=0.7$]{
    \label{fig:subfig:b}
  \includegraphics[width=3in, height=2.5in]{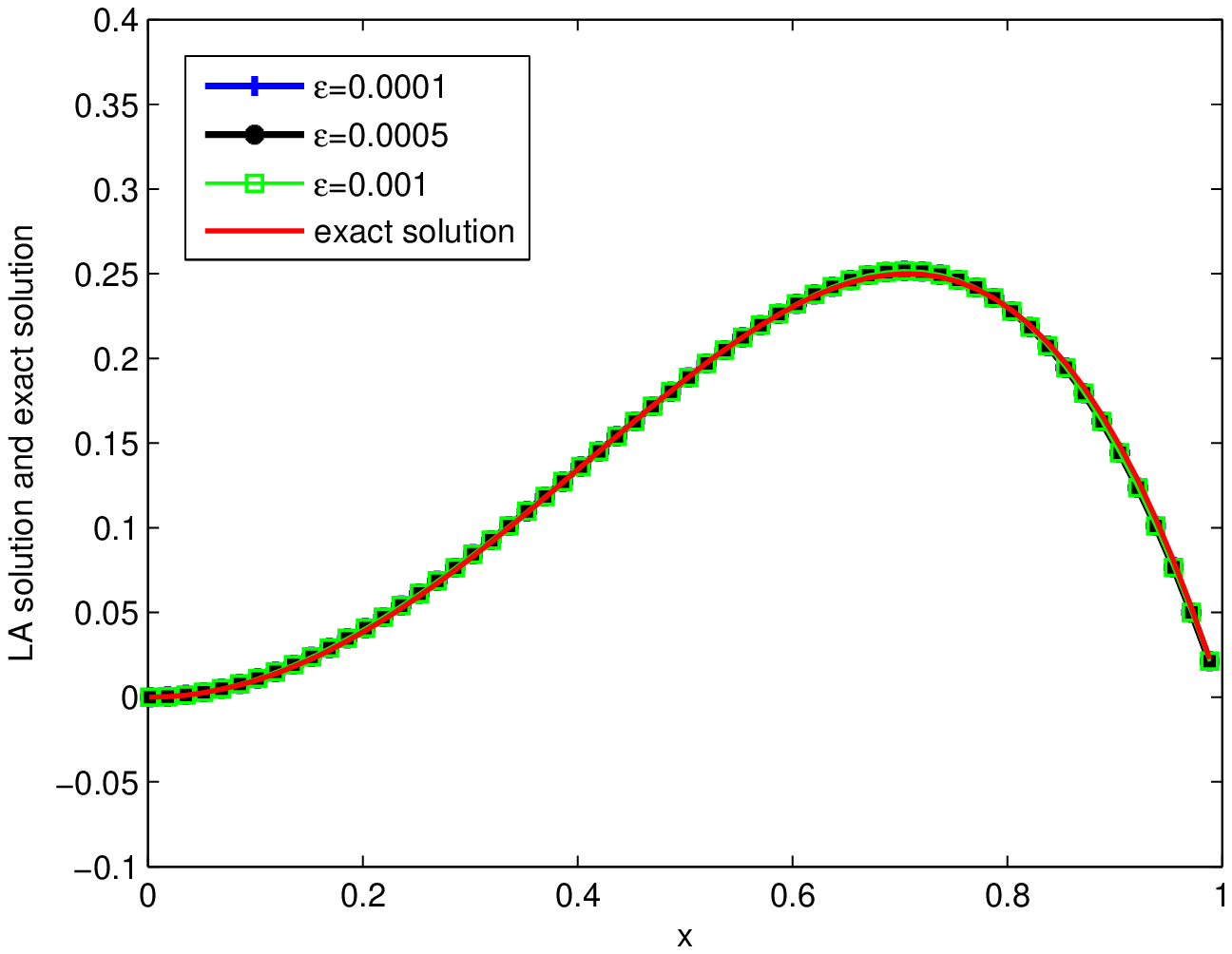}}\\
  \caption{The numerical result for Example 1 for different noise levels, and the type of data is Direct flux data}
  \label{ direct data }
\end{figure}

\begin{table}[H]
\small
\centering
\caption{ Numerical results for Example 1 for various type of data and $\varepsilon$ with $\alpha=0.3$.}
\begin{tabular}{c|c|c|c|c|c|c}
\Xhline{1pt}
 Type of data$\setminus\varepsilon $   & $0.0001$ & $0.0005$ & $0.001$ & $0.005$& $0.01$& $0.05$\\
  \hline
 Average flux data & $0.0090$ & $0.0197$ &  $ 0.0295$ &  $ 0.0662$&  $ 0.1371$& $0.2307$\\
 \hline
 Direct flux data & $ 0.0116$&  $ 0.0125$ &  $ 0.0203$ &  $ 0.0445$&  $ 0.0497$& $0.0539$\\
\Xhline{1pt}
\end{tabular}
\label{tab1-data}
\end{table}

\begin{figure}[H]
\centering
\subfigure[$\varepsilon=0.0001$]{
    \label{fig:subfig:a}
  \includegraphics[width=3in, height=2.4in]{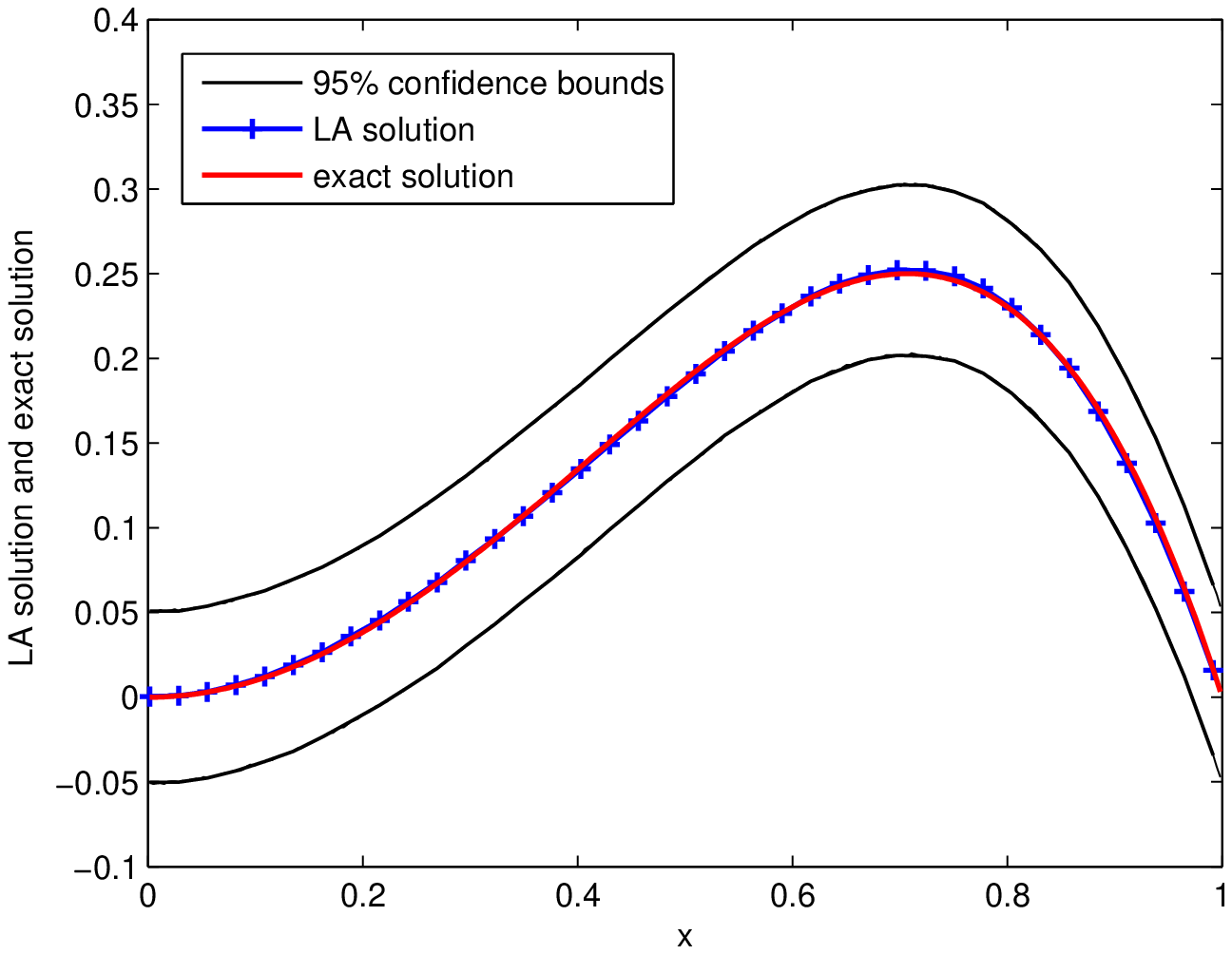}}
  \subfigure[$\varepsilon=0.0005$]{
    \label{fig:subfig:b}
  \includegraphics[width=3in, height=2.4in]{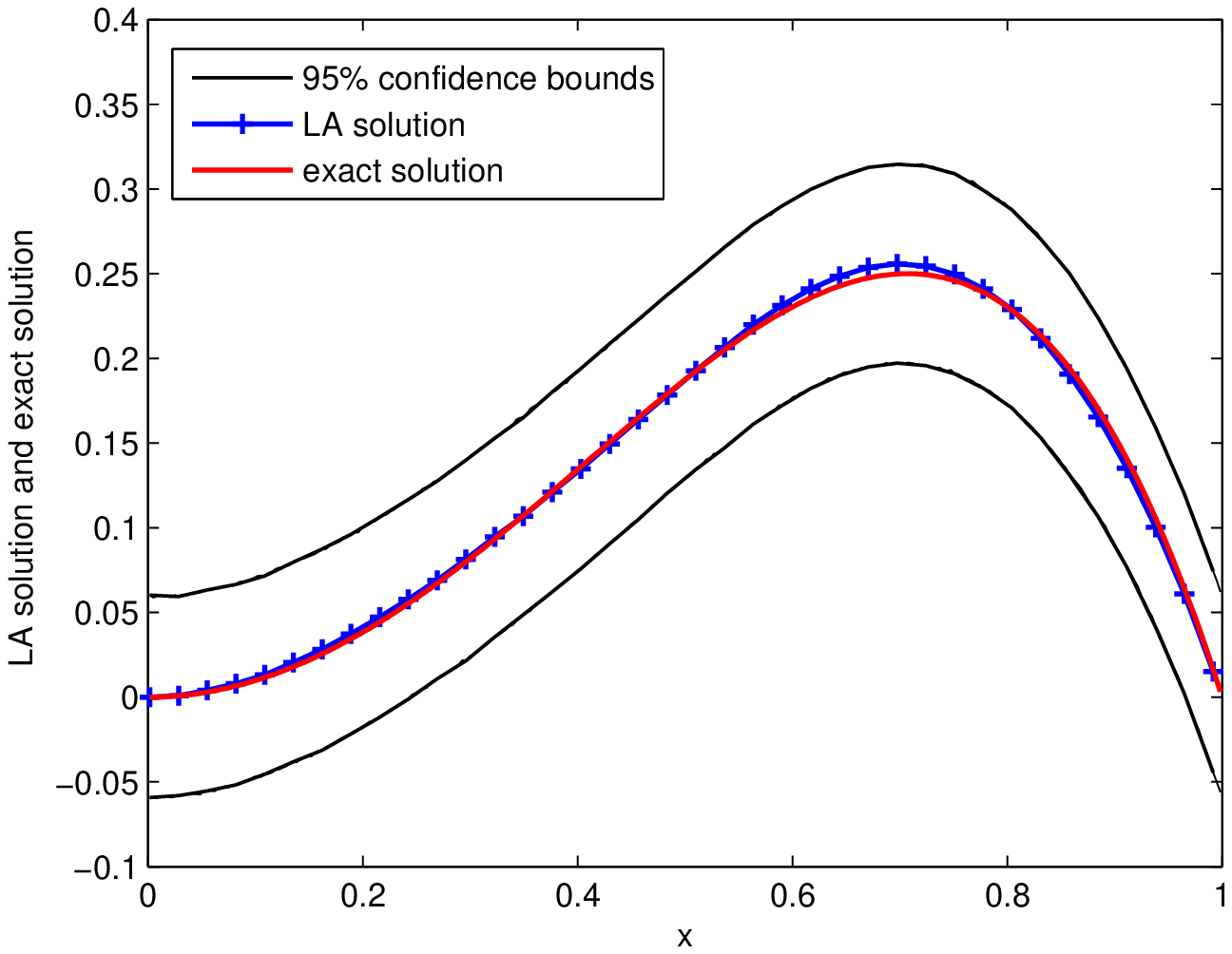}}\\
  \caption{The numerical result for Example 1 with $\mu=\delta^{\frac{3}{2}}$ and $\alpha=0.3$}
  \label{tu-confidence interval }
\end{figure}
\begin{figure}[H]
\centering
\subfigure[$\alpha=0.3$ ]{
    \label{fig:subfig:a}
  \includegraphics[width=3in, height=2.4in]{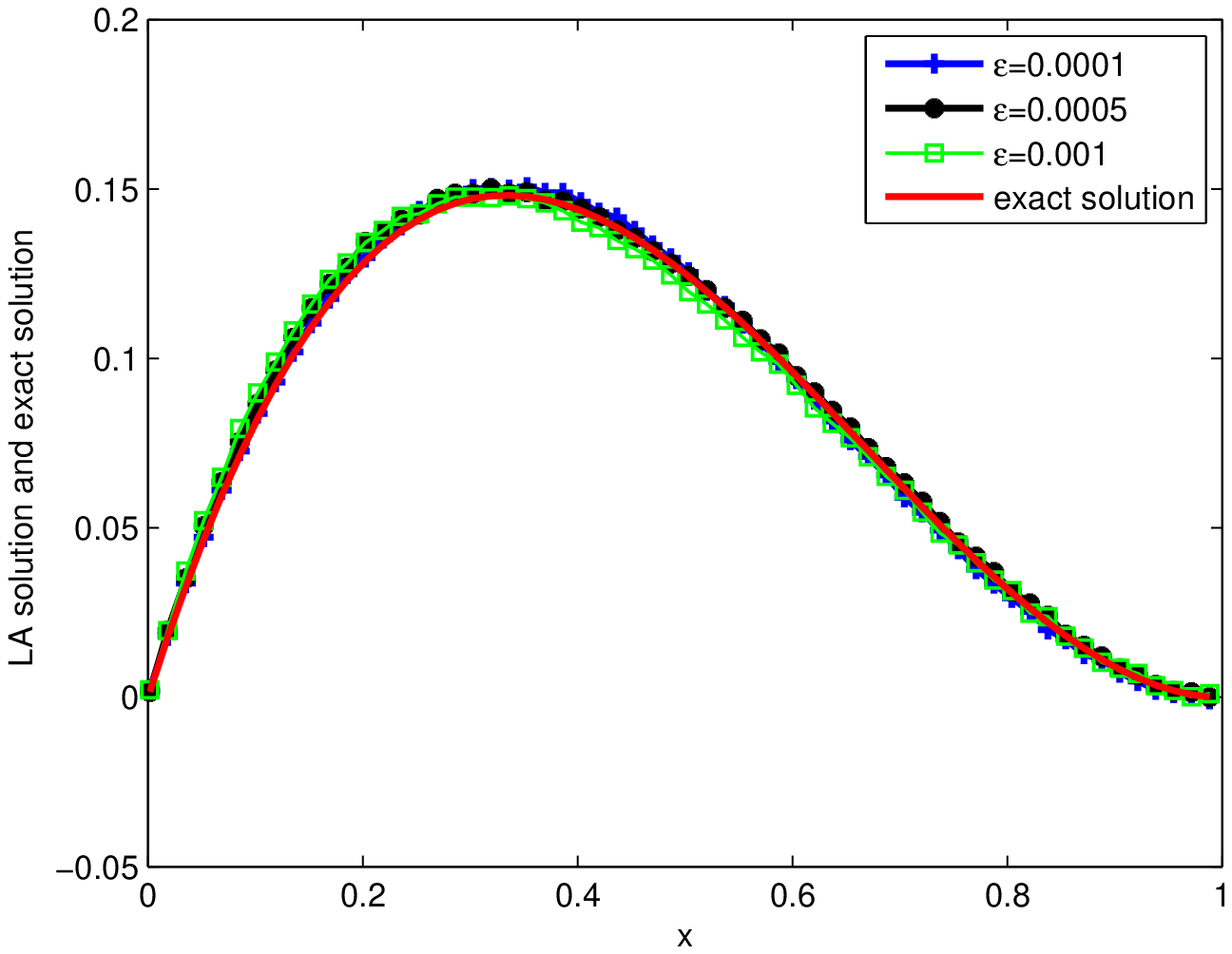}}
  \subfigure[$\alpha=0.7$]{
    \label{fig:subfig:b}
  \includegraphics[width=3in, height=2.4in]{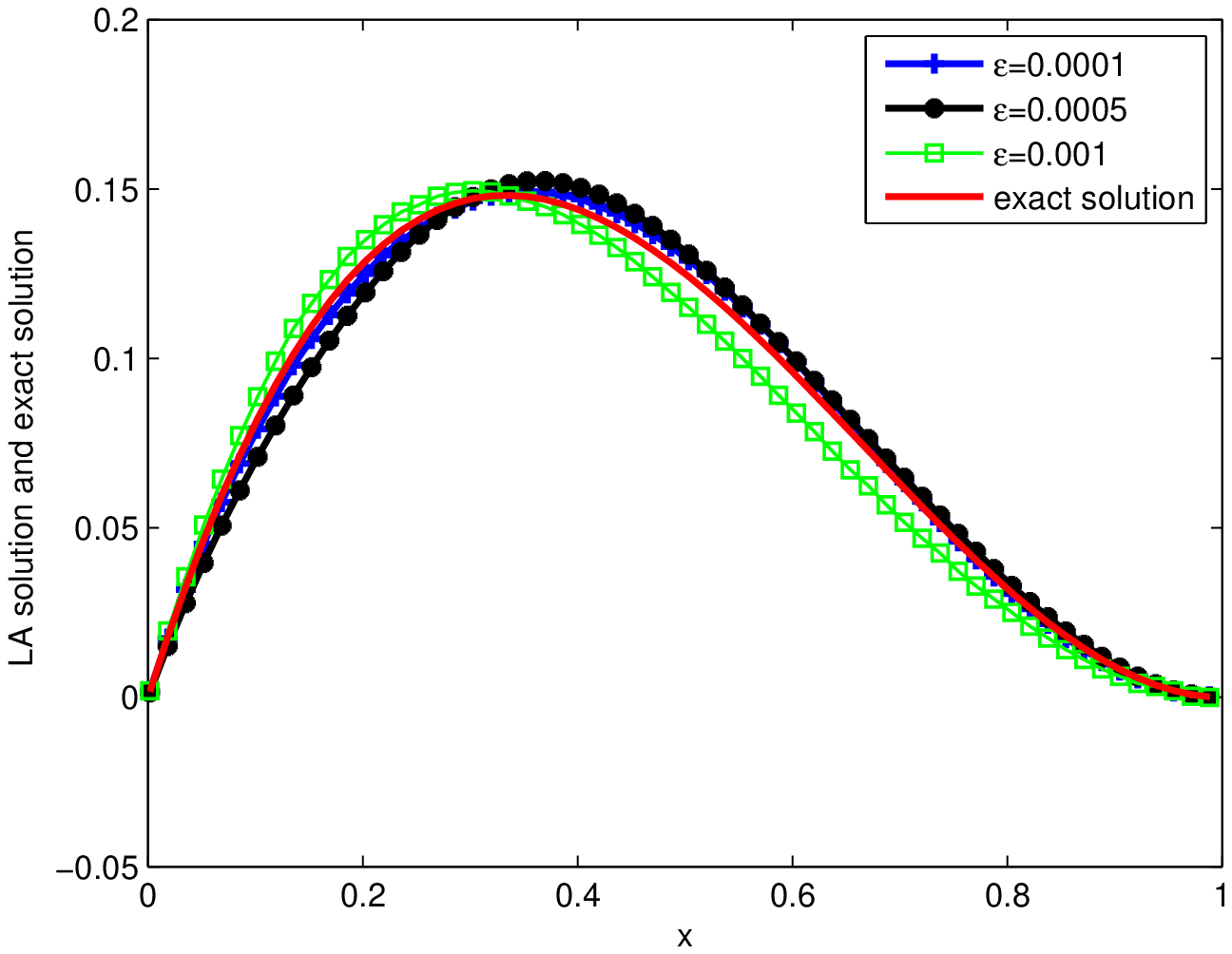}}\\
  \caption{The numerical result for Example 2 for various noise levels with $\mu=\delta^{\frac{3}{2}}$}
  \label{tu-left}
\end{figure}
\begin{figure}[H]
\centering
\subfigure[$\varepsilon=0.0001$]{
    \label{fig:subfig:a}
  \includegraphics[width=3in, height=2.4in]{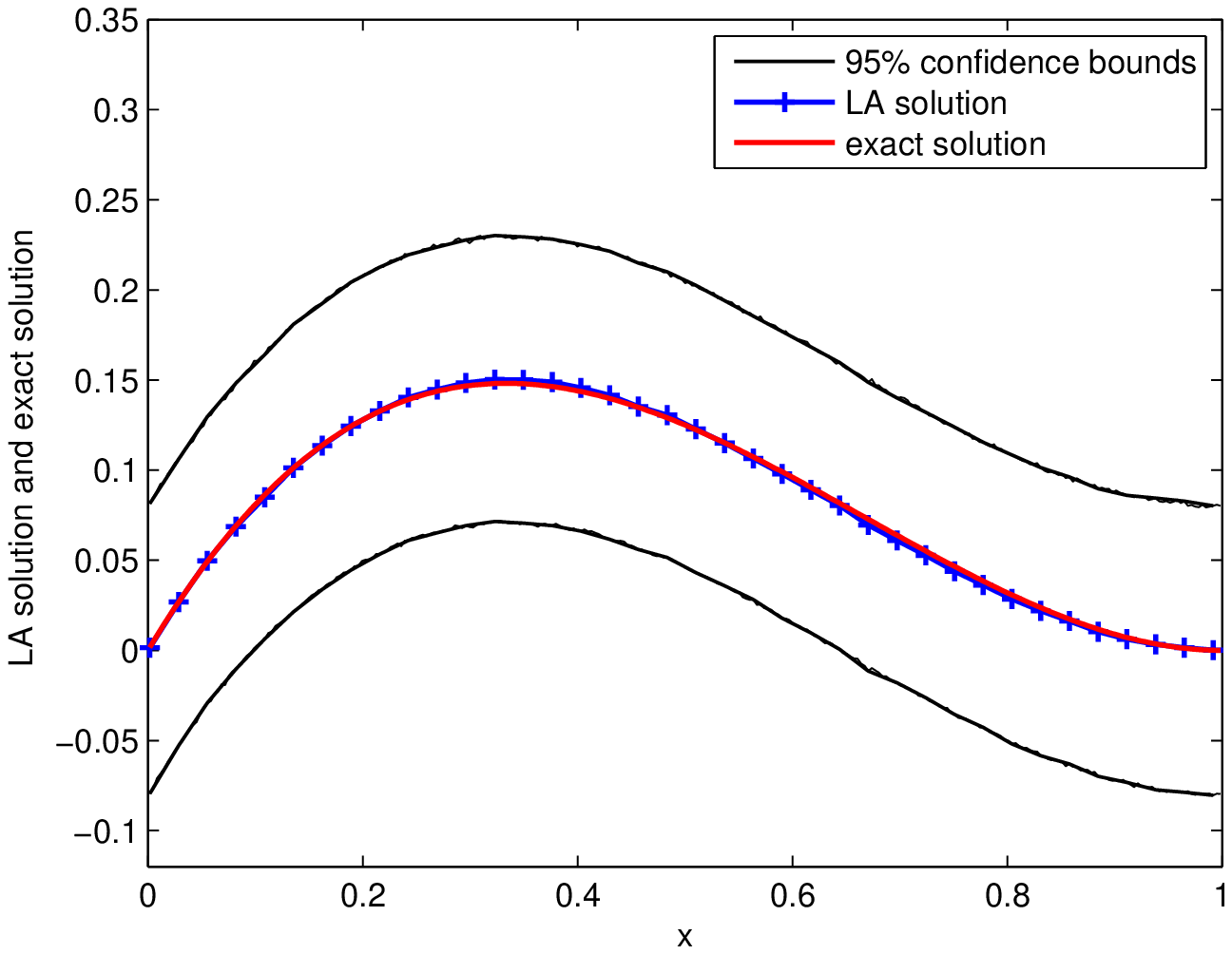}}
  \subfigure[$\varepsilon=0.0005$]{
    \label{fig:subfig:b}
  \includegraphics[width=3in, height=2.4in]{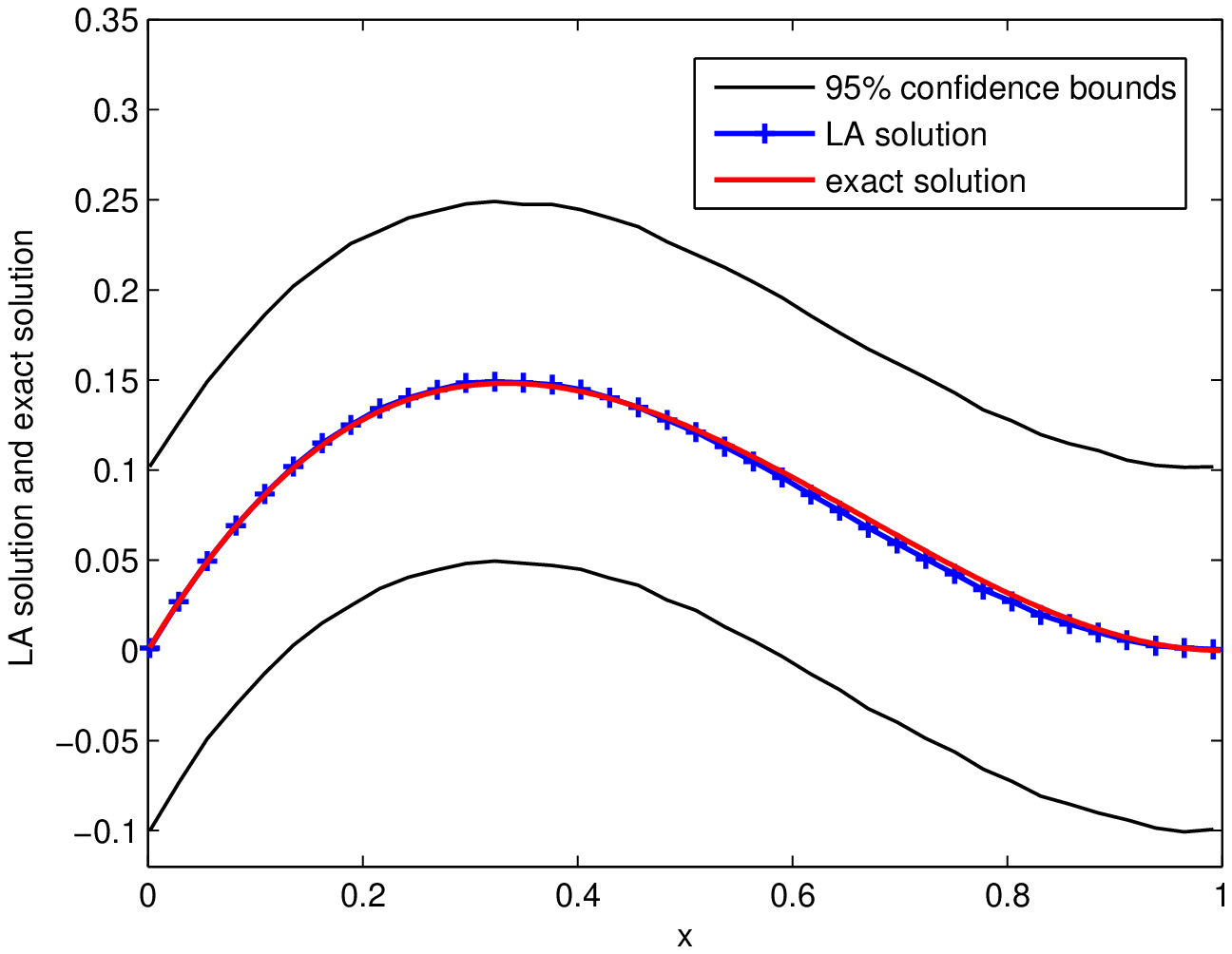}}\\
  \caption{The numerical result for Example 2 with $\mu=\delta^{\frac{3}{2}}$ and $\alpha=0.3$}
  \label{tu-conf left}
\end{figure}
\begin{figure}[H]
\label{Example.3.1}
\centering
\subfigure[$\alpha=0.3$ ]{
    \label{fig:subfig:a}
  \includegraphics[width=3in, height=2.4in]{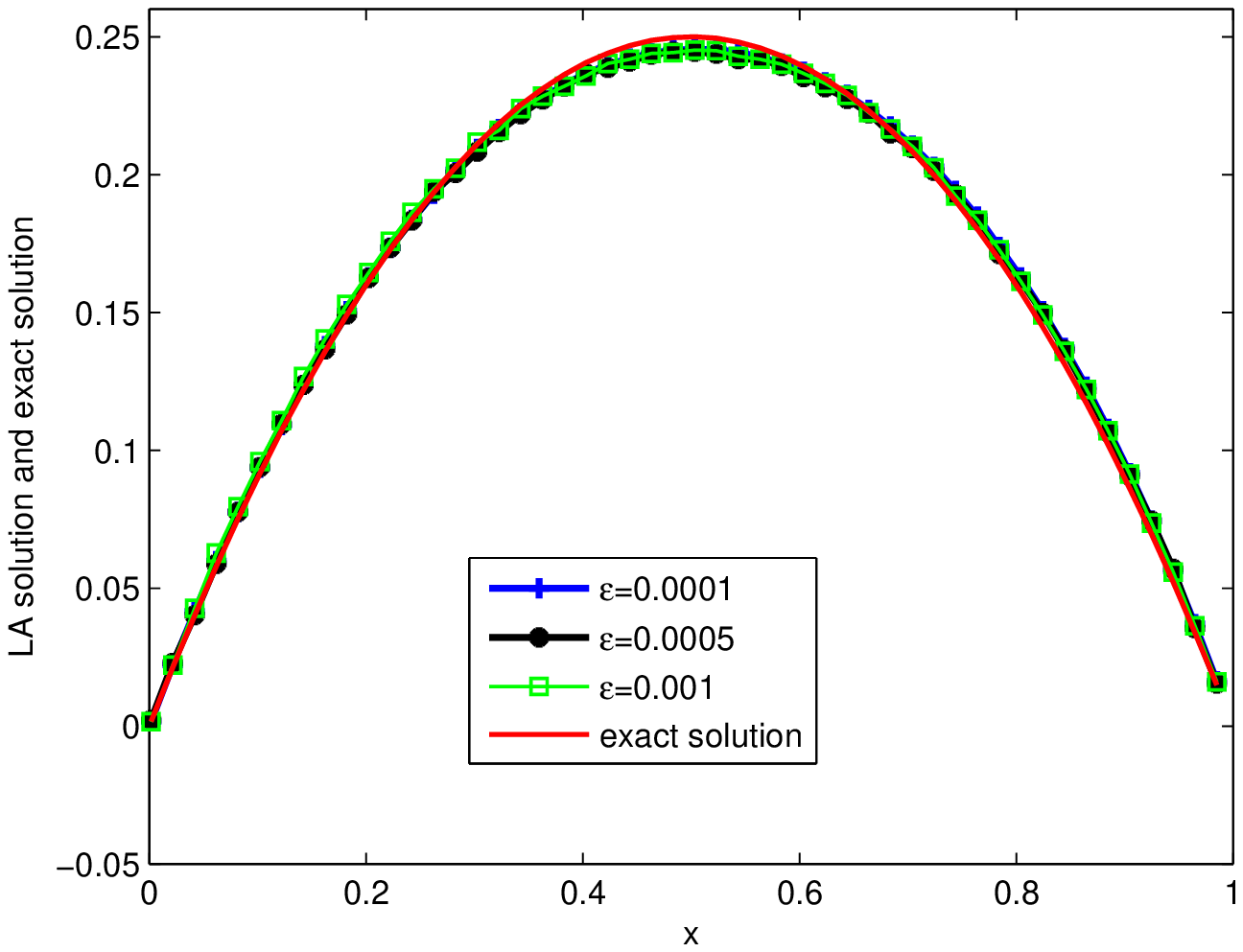}}
  \subfigure[$\alpha=0.7$]{
    \label{fig:subfig:b}
  \includegraphics[width=3in, height=2.4in]{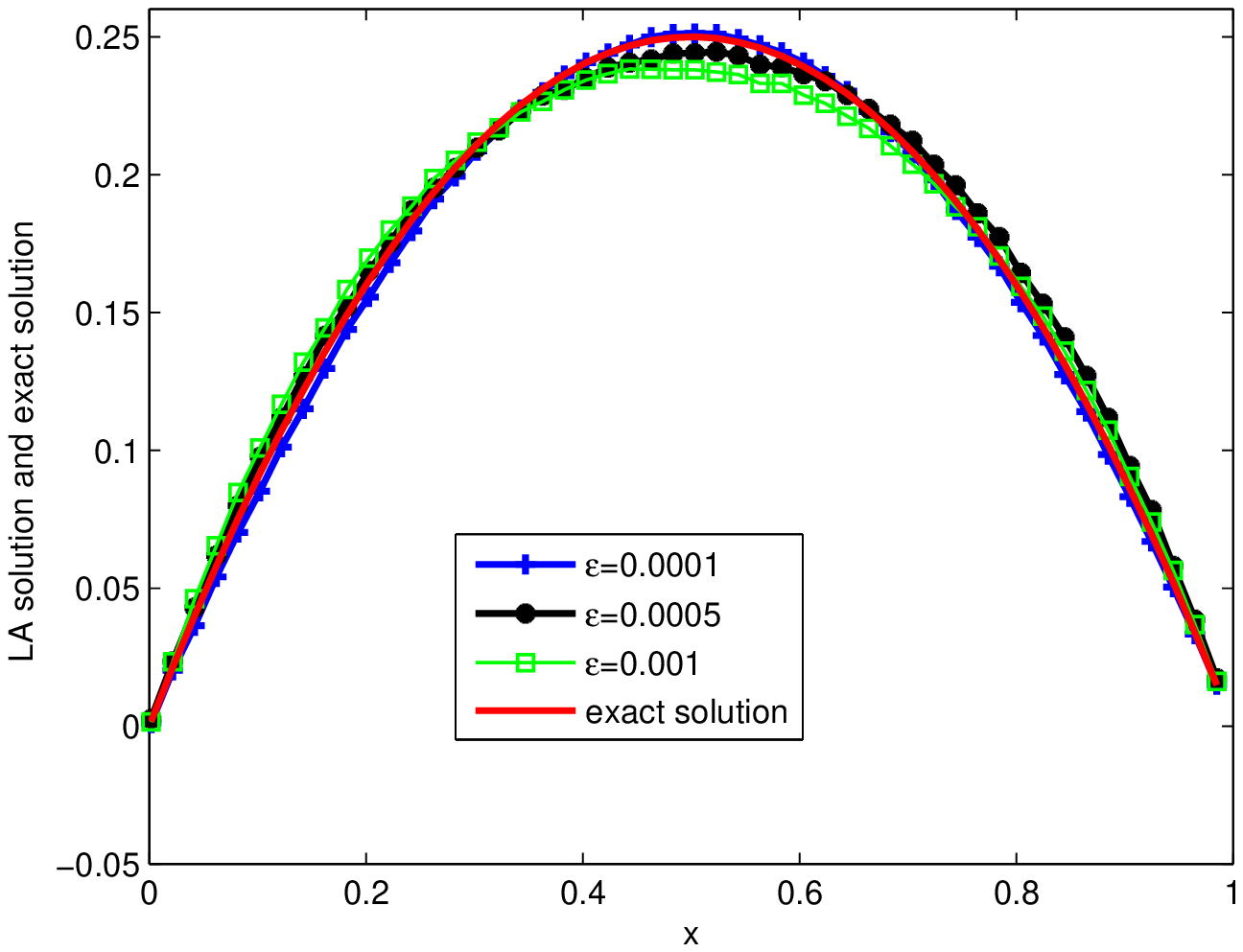}}\\
  \caption{The numerical result for Example 3 for various noise levels with $\mu=\delta^{\frac{3}{2}}$}
  \label{tu-symmetry}
\end{figure}
\begin{figure}[H]
\centering
\subfigure[$\varepsilon=0.0001$ ]{
    \label{fig:subfig:a}
  \includegraphics[width=3in, height=2.4in]{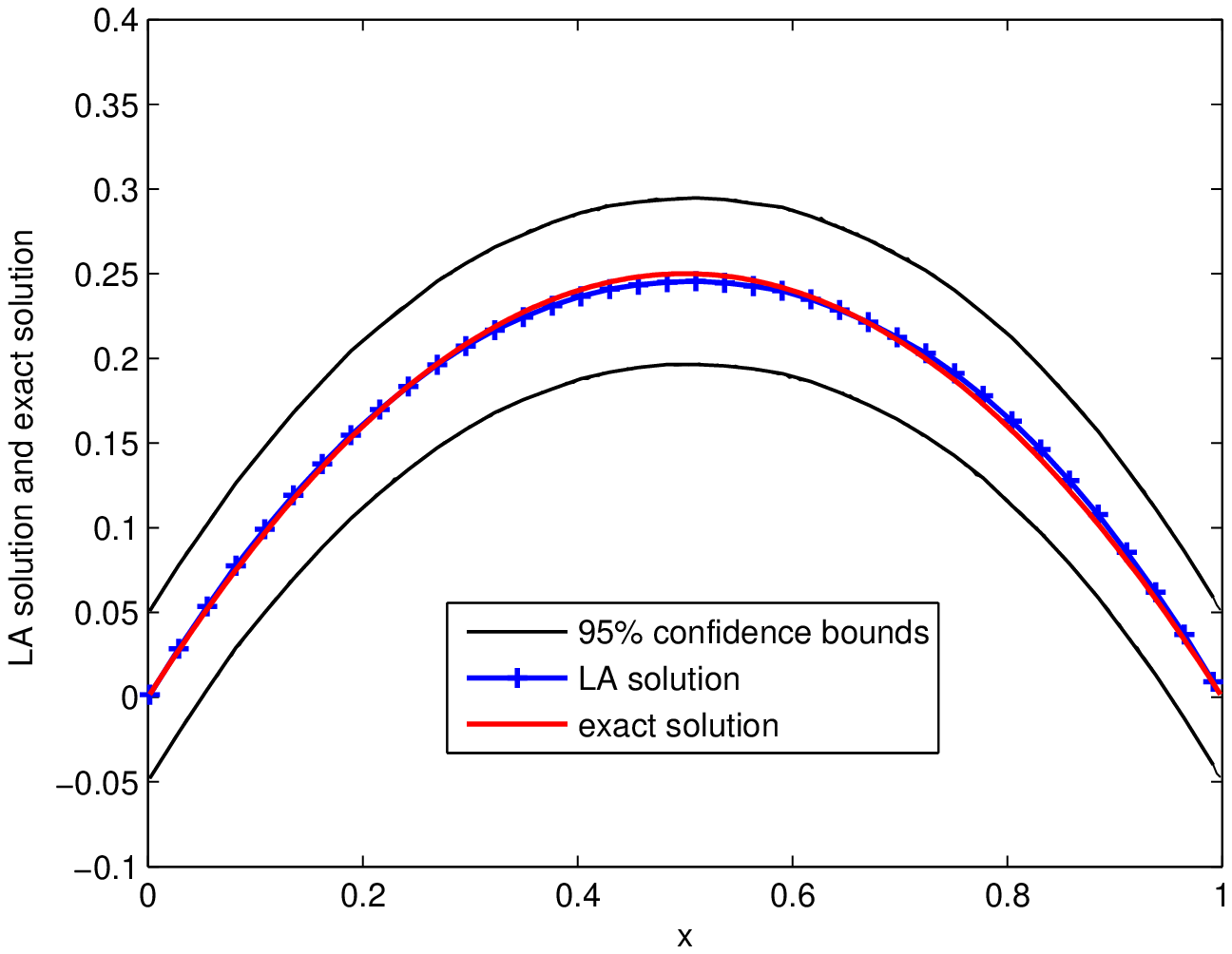}}
  \subfigure[$\varepsilon=0.0005$]{
    \label{fig:subfig:b}
  \includegraphics[width=3in, height=2.4in]{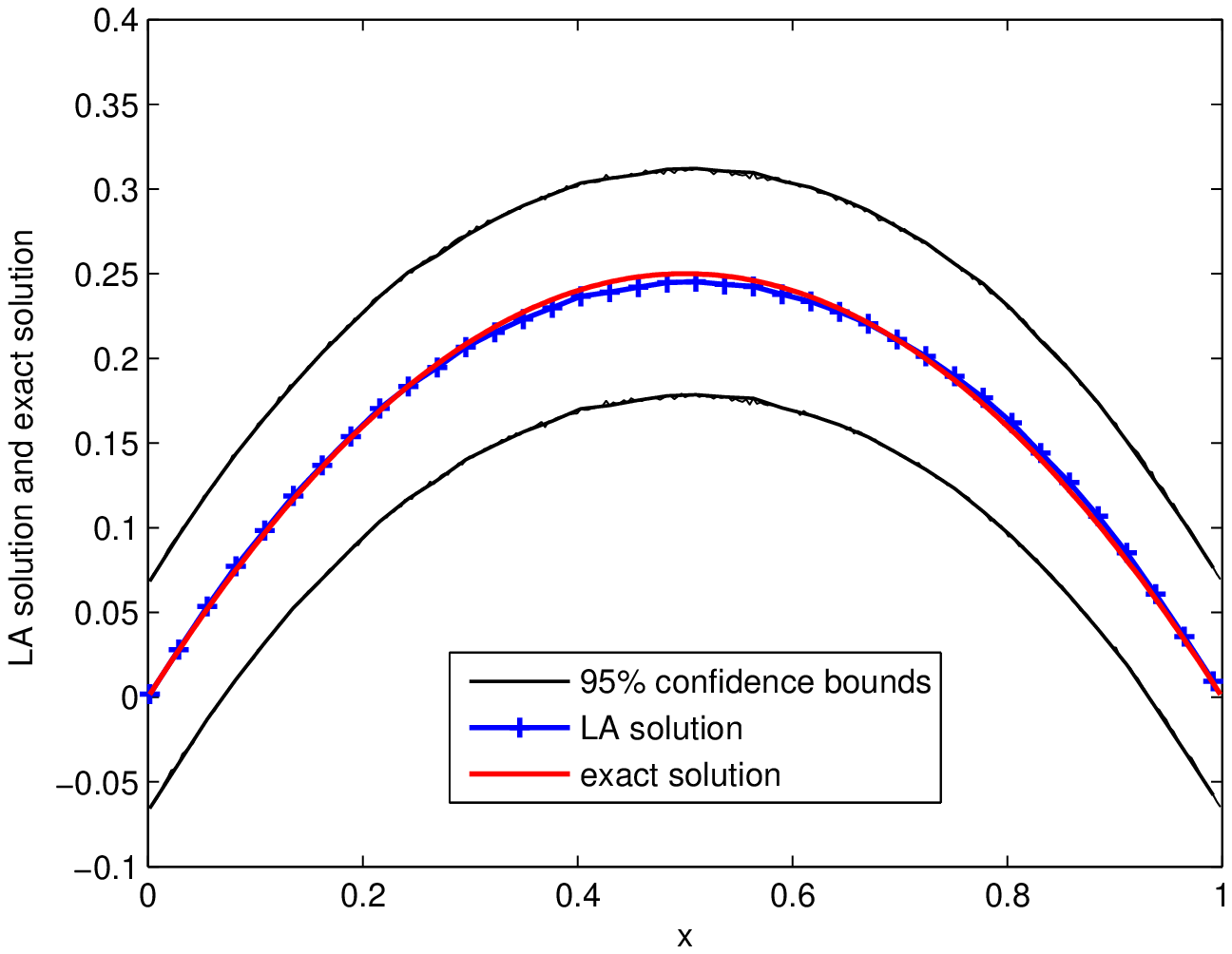}}\\
  \caption{The numerical result for Example 3 with $\mu=\delta^{\frac{3}{2}}$ and $\alpha=0.3$}
  \label{tu-symmetry2}
\end{figure}

  For Example 1, we draw the 95\% confidence interval in Figure $\ref{tu-confidence interval }$ for the noise level $\varepsilon=0.0001,\ 0.0005$ with $\alpha=0.3$. The posterior mean $\bar{q}_n(x)$ is in excellent agreement with the exact solution, and the confidence interval quantifies its associated uncertainty. The confidence interval shrinks as the noise level $\varepsilon$ decreases. We observe that the confidence interval on one side near the observation data is relatively narrow, and the corresponding confidence interval is also relatively accurate. The confidence interval far from the observation data is wide and the corresponding confidence interval is relatively inaccurate. It's a surprise that the position of the data affects the accuracy of our reconstruction results. In Example 2, Example 3, we will continue to verify this result.\\
  \\
$\mathbf{Example ~~2}$. In $(\ref{mix-eq1})$, the accessible boundary $\Lambda$ is taken to be $\Lambda_0$, and the degradation coefficient is given by
\begin{eqnarray}
\label{eqq4}
\begin{split}
q(x)=x(1-x)^2,\ \  x\in\Omega.
\end{split}
\end{eqnarray}
The numerical results for Example 2 for various levels of noise in the data are shown in Figure $\ref{tu-left}$ with $\alpha=0.3,\ 0.7$. The results in Figure $\ref{tu-left}(b)$  is worse than  Figure $\ref{tu-left}(a)$ as the fractional order $\alpha$ increase. Moreover, we also find that the position of the measurement data have a great influence on the reconstruction results. In the Figure $\ref{Example.1 for various noise le}$, we let $\Lambda=\Lambda_1$, the average flux data are measured on the $\Lambda_1$ for Example 1, when the exact solution $q$ is biased to the right of the region, we obtain good numerical results. In the Figure $\ref{tu-left}$, setting $\Lambda=\Lambda_0$, when the average flux data are measured on the $\Lambda_0$ for Example 2, the left biased $q$ give desired approximate result (see Example 3 for further discussion). We used LA algorithm to sample, calculated the corresponding posterior mean value $\bar{q}_n$, and drew 95\% confidence interval in the Figure $\ref{tu-conf left}$, we can get similar results of Example 1.  The confidence interval shrinks as the noise level $\varepsilon$ decreases and  near observation data is more accurate, far from the observation data is wide and inaccurate, see the Figure $\ref{tu-conf left}$.
$\mathbf{Example ~~3}$. In $(\ref{mix-eq1})$, the accessible boundary $\Lambda$ is taken to be $\Lambda=\Lambda_0\cup \Lambda_1$, and the degradation coefficient is given by
\begin{eqnarray}
\label{eqq4}
\begin{split}
q(x)=x(1-x), x\in\Omega.
\end{split}
\end{eqnarray}
The numerical results for Example 3 with various fractional order and noise levels are presented in Figure $\ref{tu-symmetry}$, it shows that the fractional order $\alpha$ increases but the result is not good.  When we measure average flux data from both sides, and use LA method to recover the degradation coefficient, the error is small with noise levels are $\varepsilon=0.0001,\ 0.0005,\ 0.001$ respectively. Furthermore, the confidence interval has higher precision on both sides and lower precision in the middle of the region. The confidence interval are shown in Figure $\ref{tu-symmetry2}$. This is consistent with our previous conclusions.

\begin{figure}[ht]
\centering
\includegraphics[scale=0.6]{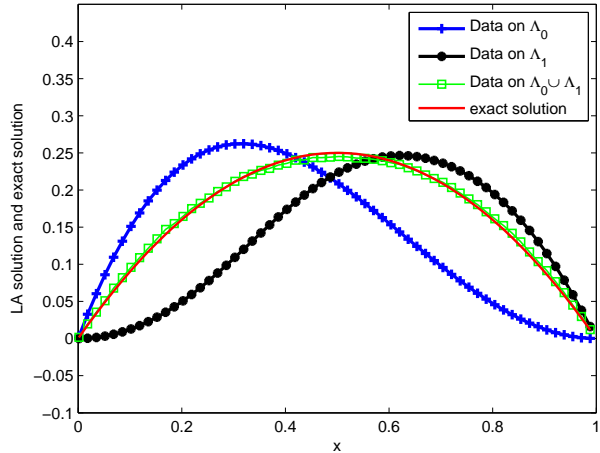}
\caption{Numerical results for Example 3 for different  position  of average flux data with noise level $\varepsilon=0.001$}
\label{The figure of skewness}
\end{figure}

 Now, we introduce the concept of skewness to describe the degree of LA solution's deviation from the mean value, and further demonstrate the relation between the symmetry of solution and the symmetry of measurement data. The degradation coefficient is taken as Example 3, and the boundary average flux data are measured on different positions. Then the corresponding numerical results for Example 3 are displayed in Figure $\ref{The figure of skewness}$. As the  error level increase to 0.001, when measuring the data on the boundary of $\Lambda_1$, the degradation coefficient recovered by LA method is left skewed, i.e., its skewness is positive. When the average flux data are measured on the boundary of $\Lambda_0$, the degradation coefficient is right skewed, i.e., its skewness is negative. When we measure average flux data on both sides of the boundary (symmetric data), that is to say $\Lambda=\Lambda_1\cup\Lambda_0$, the reconstructed degradation coefficient is also symmetric and give a accurate approximation, and its skewness is zero. Thus, we can use the symmetry of measurement data to capture the symmetry feature of the identified object. By using (\ref{ske}), the calculation of the corresponding skewness can be referred to Table $\ref{tab1-skewness}$.

\begin{table}[H]
\small
\centering
\caption{ Numerical results for Example 3 for the skewness of  exact solution is 0}
\begin{tabular}{c|c}
\Xhline{1pt}
  Measurement data   &  skewness \\
  \hline
 Data on $\Lambda_0$ & $~~0.3285$ \\
 \hline
 Data on $\Lambda_1$  & $-0.2088$ \\
 \hline
 Data on $\Lambda_0\cup\Lambda_1$  & $ -0.0033$\\
\Xhline{1pt}
\end{tabular}
\label{tab1-skewness}
\end{table}

\subsubsection{Nonsmooth function}
\label{Numerical-result-elliptic}
$\mathbf{Example ~~4}. $ In this example, we consider the more challenging case of reconstructing a nonsmooth example with a cusp, and the degradation coefficient is prescribed as follows:
\begin{eqnarray*}
q(x)&=
\begin{cases}
\begin{split}
&x,& ~0\leq x\leq \frac{2}{3}, \\
&-2x+2,  & ~\frac{2}{3}<x\leq 1,
\end{split}
\end{cases}
\end{eqnarray*}
$\mathbf{Example ~~5}. $ We also consider an discontinuous example, and the degradation coefficient is:
\begin{eqnarray*}
q(x)&=
\begin{cases}
\begin{split}
&0,& ~0\leq x\leq \frac{1}{2}, \\
&0.4,  & ~\frac{1}{2}<x\leq \frac{4}{5},\\
&0 , & ~\frac{4}{5}<x\leq 1,\\
\end{split}
\end{cases}
\end{eqnarray*}
\begin{figure}[H]
\centering
\subfigure[]{
    \label{fig:subfig:a}
  \includegraphics[width=3in, height=2.4in]{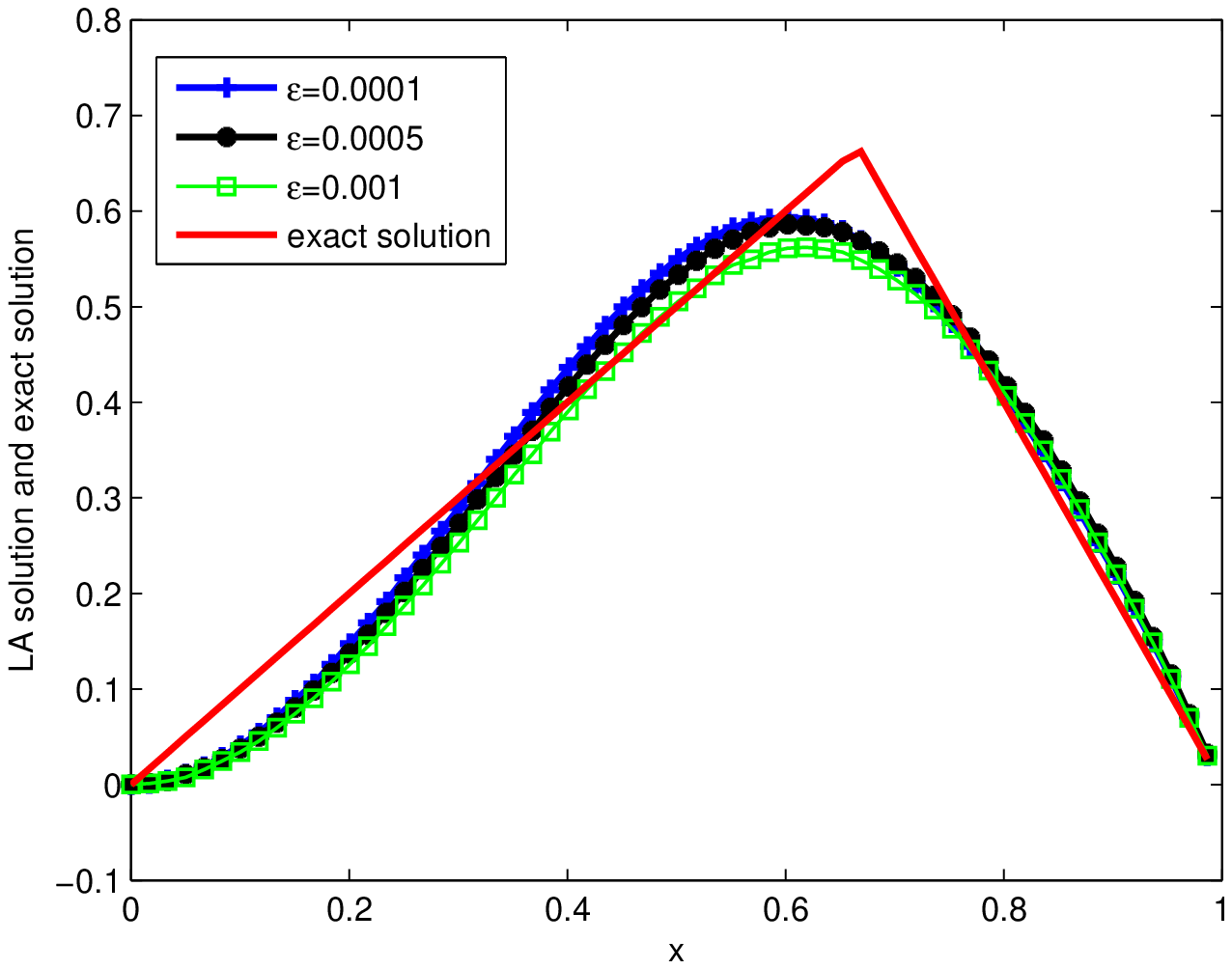}}
  \subfigure[]{
    \label{fig:subfig:b}
  \includegraphics[width=3in, height=2.4in]{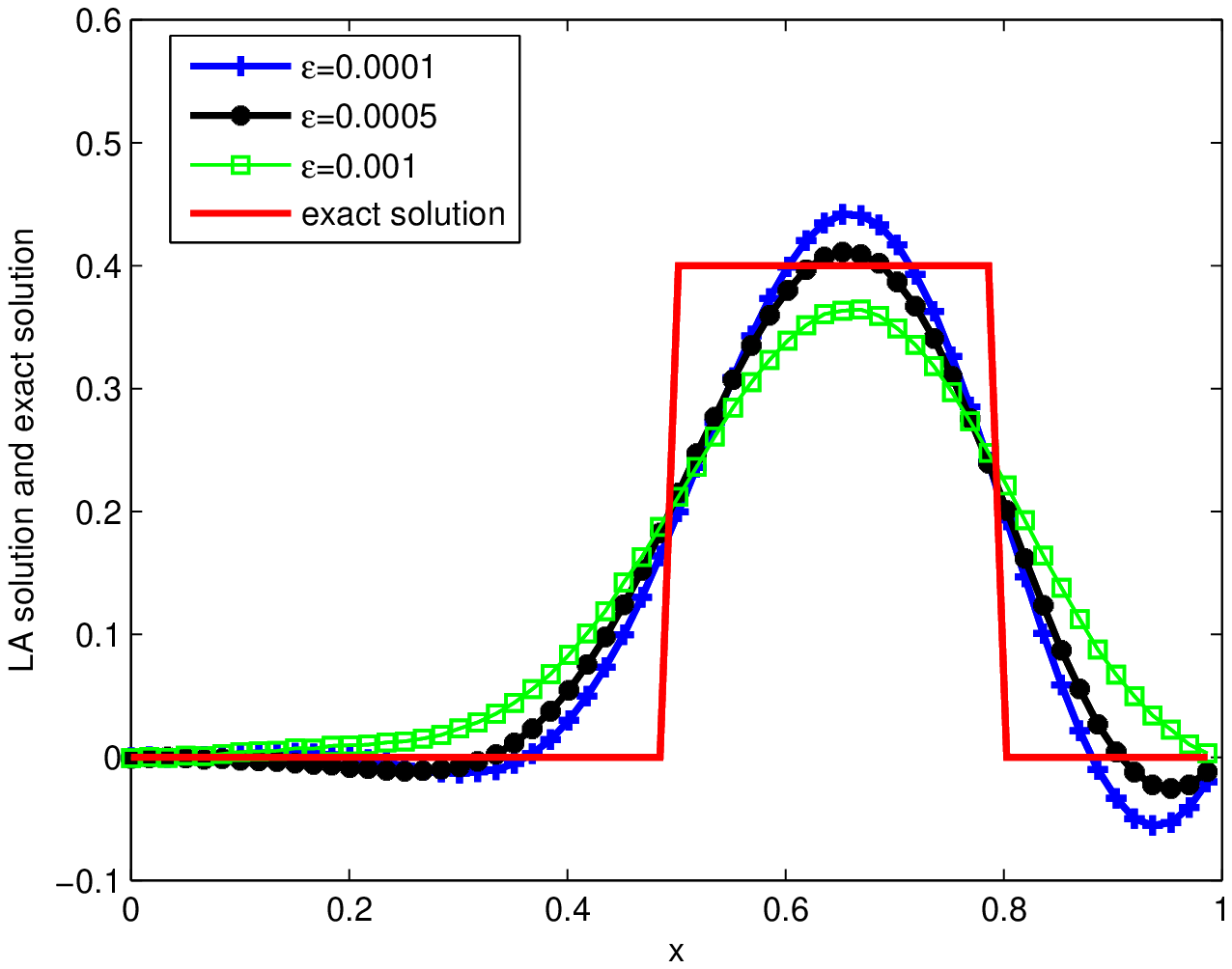}}\\
  \caption{The numerical result for (a) Example 4 and (b) Example 5 for various noise levels with regularization parameter
  $\mu=\delta^{\frac{3}{2}}$}
  \label{Nonsmooth function}
\end{figure}
The numerical results for Example 4 and Example 5 for various noise levels in the case of $\alpha=0.3$ are shown in the Figure $\ref{Nonsmooth function}$ with the regularization parameter $\mu=\delta^{\frac{3}{2}}$. It can be seen that the smaller the noise level, the better the numerical results. The error for  nonsmooth Example 4 in the neighborhood of the cusp is large because of the smoothing nature of the prior. And the same result can be obtained from Example 5, the numerical results obtained in the neighborhood of segment point are not very desired. For the nonsmooth case, different regularization methods and different prior information are needed, such as TV prior \cite{vogel2002}, and we do not discuss the details here.

\subsection{Inversion for two-dimensional degradation coefficient}
\label{Numerical-result-elliptic}
\begin{figure}[H]
\centering
\subfigure[]{
    \label{fig:subfig:a}
  \includegraphics[width=3in, height=2.4in]{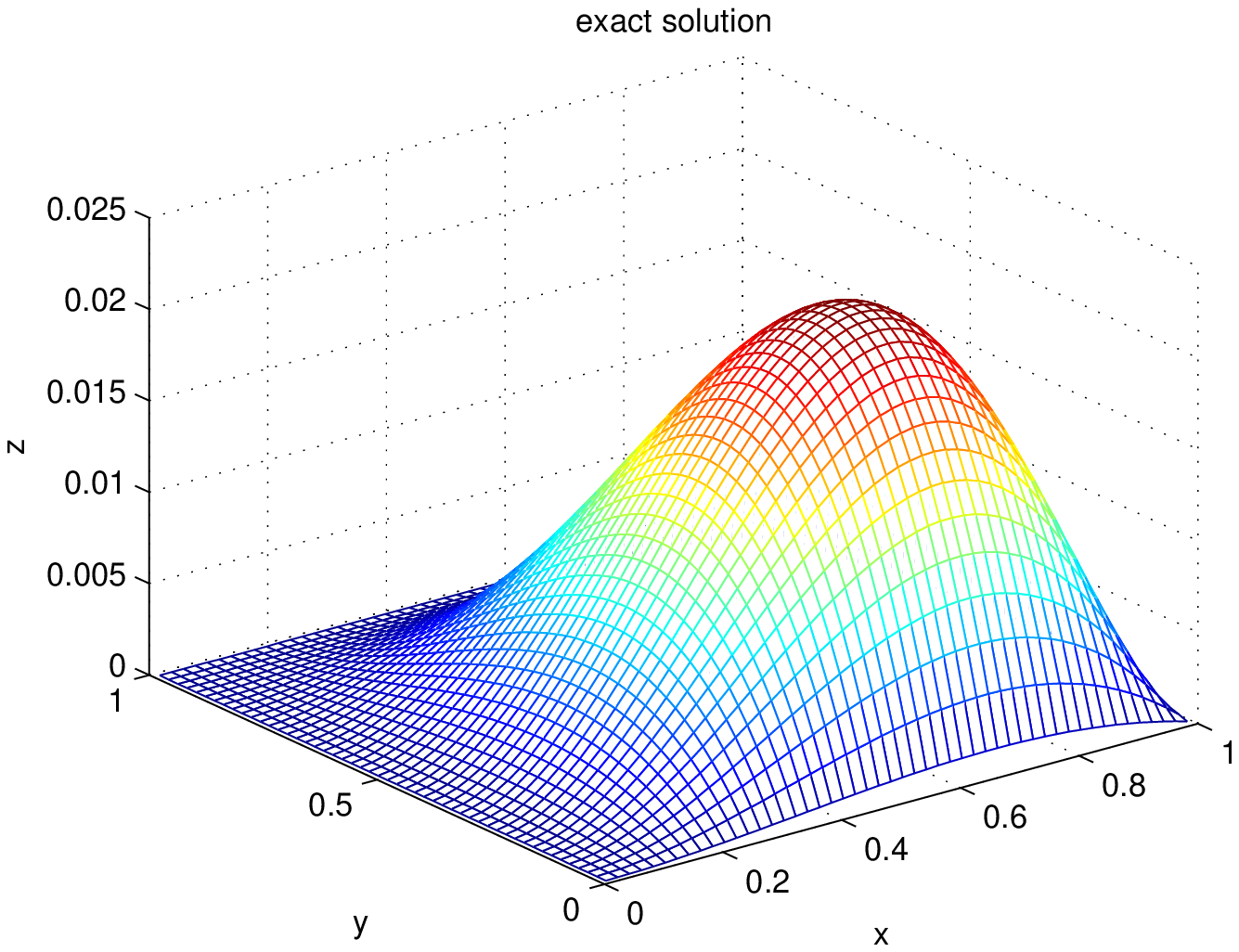}}
  \subfigure[]{
    \label{fig:subfig:b}
  \includegraphics[width=3in, height=2.4in]{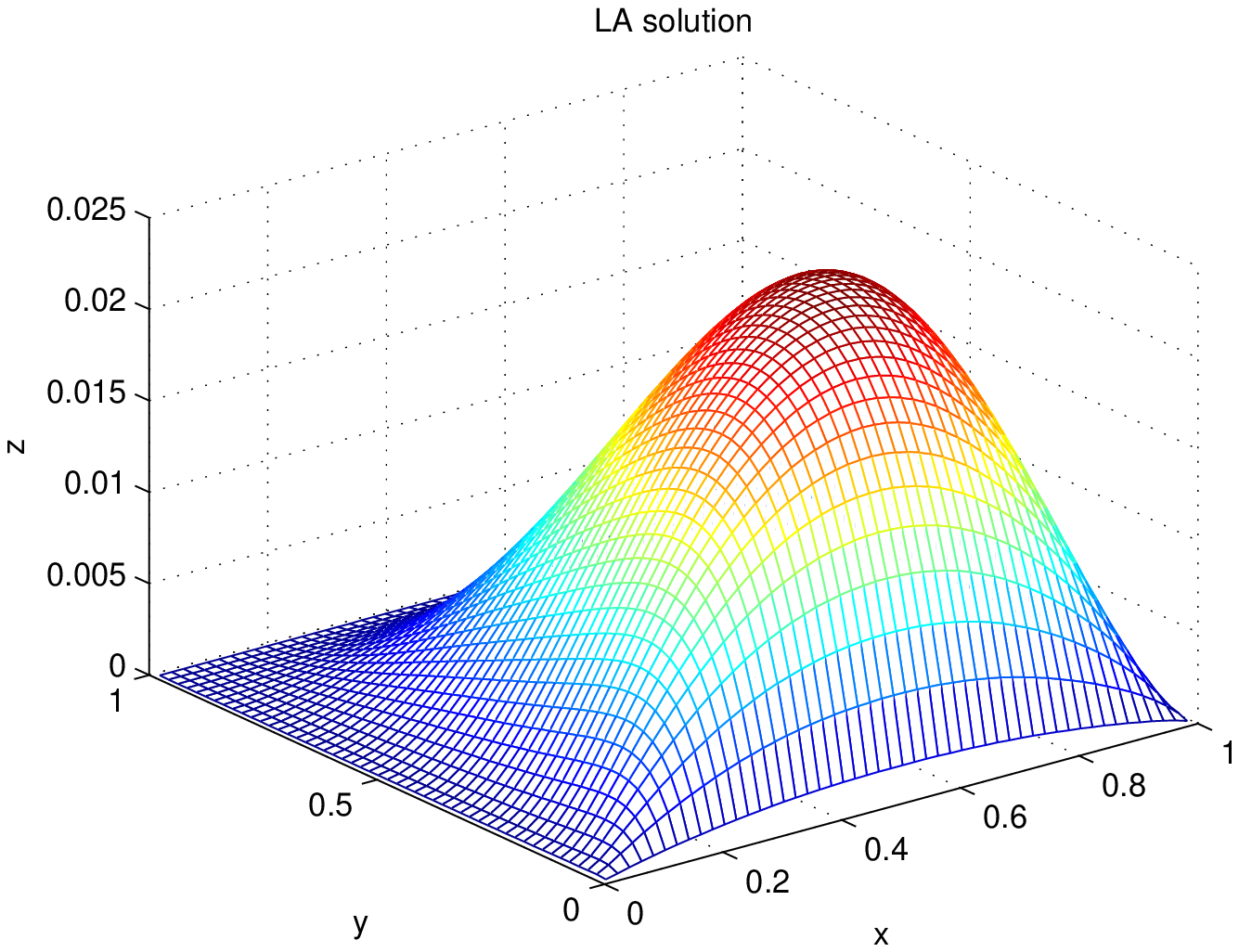}}\\
  \caption{The numerical result for Example 6 for noise level $\varepsilon=0.0001$ and $\alpha=0.3$, with $\mu=\delta^{\frac{1}{2}}$}
  \label{two-1}
\end{figure}

The domain $\Omega$ under consideration is a unit square $[0,1]\times[0,1]$. Set $T=1$, and the grid point on $[0,T]$ is $101$. In the equation $(\ref{mix-eq1})$, we take $d=2$, the diffusion coefficient matrix is unitary. The forward problem is discretized using $2500$ uniform rectangular finite element. The number of basis functions is take N=5 and sample size $Ne=10000$. Set positive function $h(x,t)=(t-1)$. The $\partial\Omega=\Lambda_1\cup\Lambda_2\cup\Lambda_3\cup\Lambda_4$, where $\Lambda_1=(0,1]\times\{0\}$, $\Lambda_2=\{1\}\times(0,1]$, $\Lambda_3=[0,1)\times\{1\}$, $\Lambda_4=\{0\}\times (0,1)$.
 Next, we presents the numerical results for 2D cases. \\
\\
$\mathbf{Example ~~6}.$  We take the $\Lambda=\Lambda_1\cup\Lambda_2$, the degradation coefficient is:
\begin{eqnarray*}
\label{Qresidual-eq}
\begin{split}
q(x,y)=x(x-x^2)y(1-y)^2,\ \ x,y\in\Omega.
\end{split}
\end{eqnarray*}
 The LA solution of Example 6 tends to the front right of the region, we can see the Figure $\ref{two-1}$. In this case, we use the asymmetric boundary data to reconstruct the degradation coefficients, here the data is on the $\Lambda_1\cup\Lambda_2$, and it verifies the LA solution is close to the data location. The solution dependence on data location is more obvious in two dimensions.

\begin{figure}[H]
\centering
\subfigure[Exact solution]{
    \label{fig:subfig:a}
  \includegraphics[width=3in, height=2.4in]{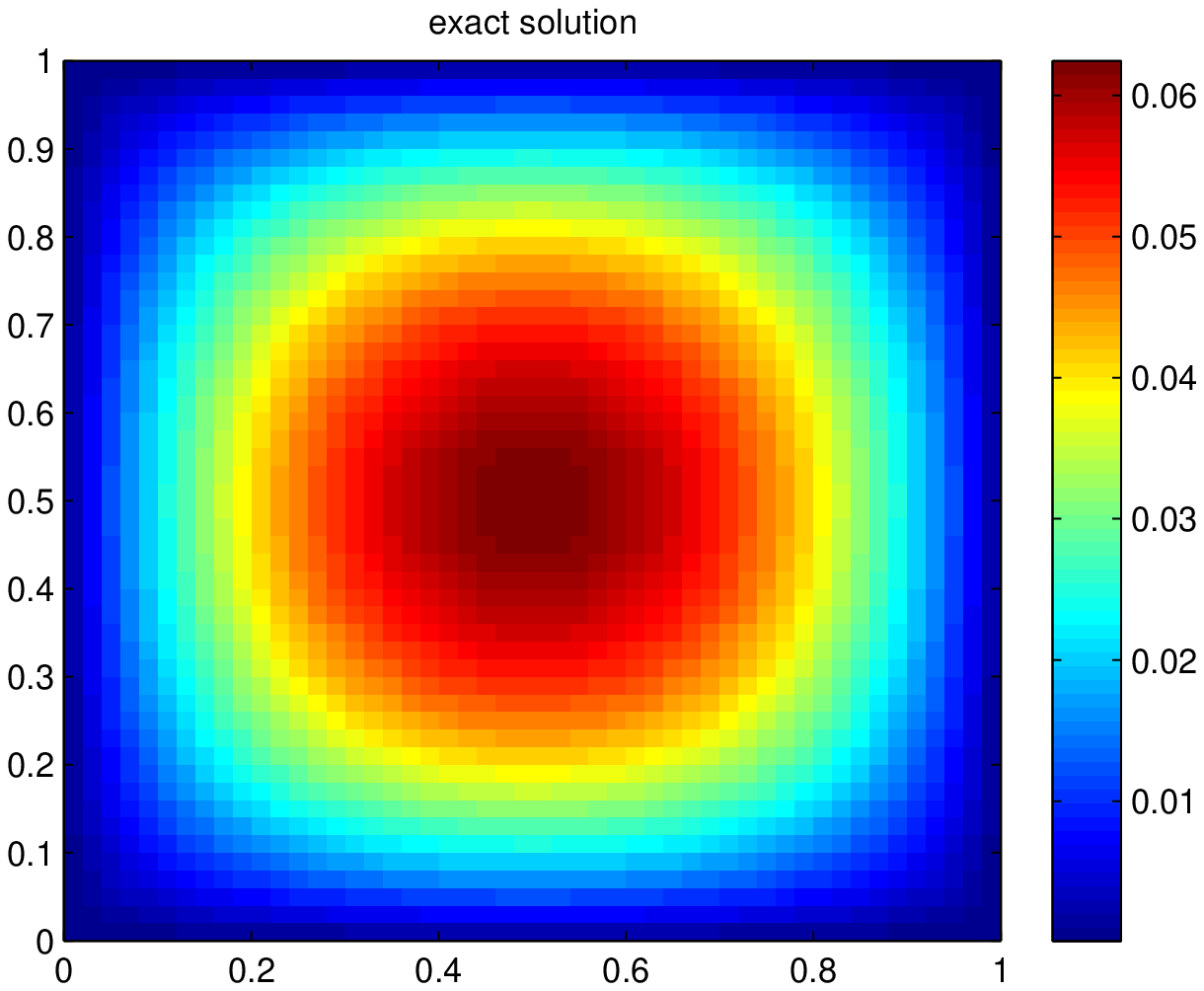}}
  \subfigure[Data on $\Lambda_1\cup\Lambda_2$]{
    \label{fig:subfig:b}
  \includegraphics[width=3in, height=2.4in]{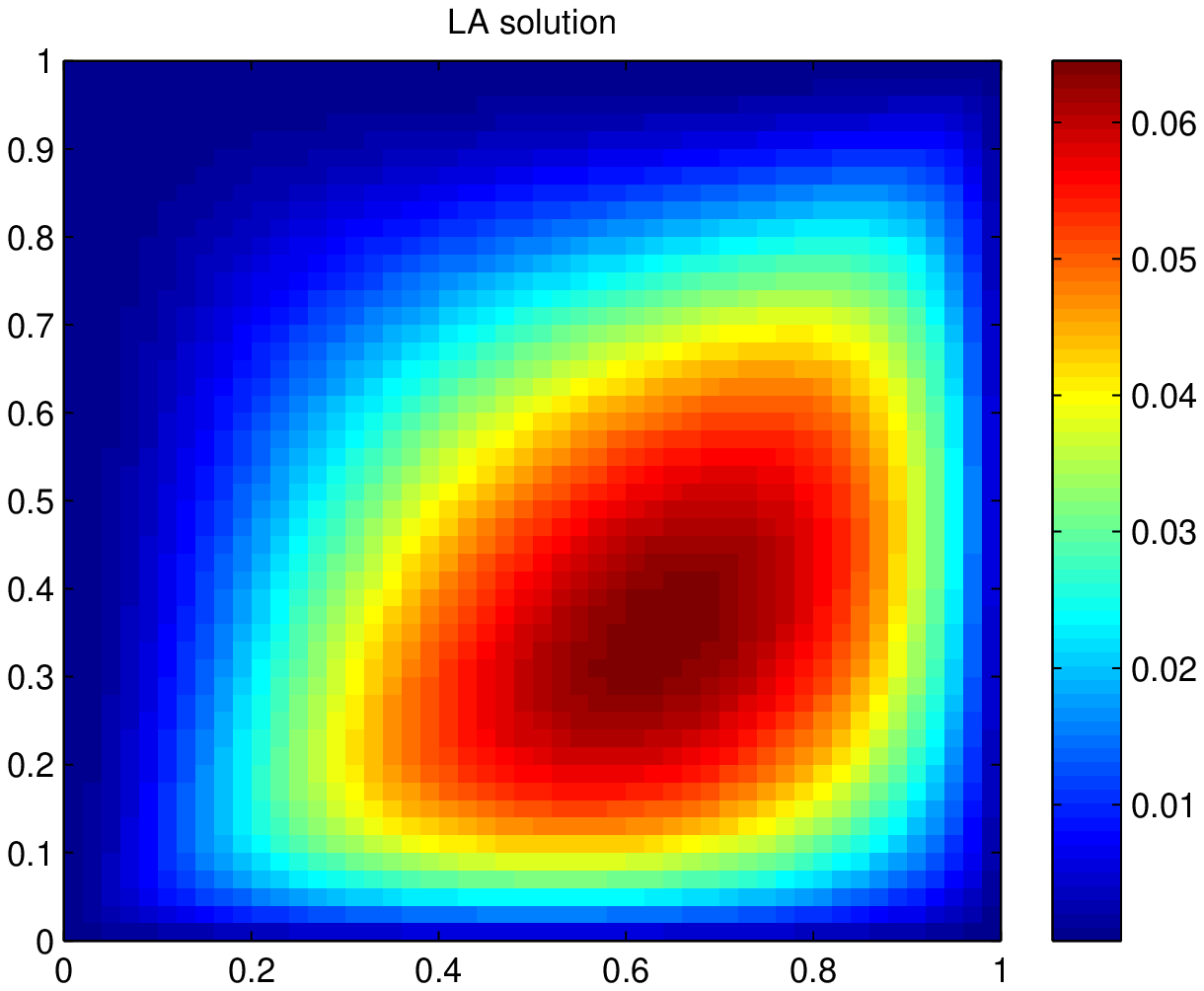}}\\
  \subfigure[Data on $\Lambda_1\cup\Lambda_2\cup\Lambda_4$]{
    \label{fig:subfig:c}
  \includegraphics[width=3in, height=2.4in]{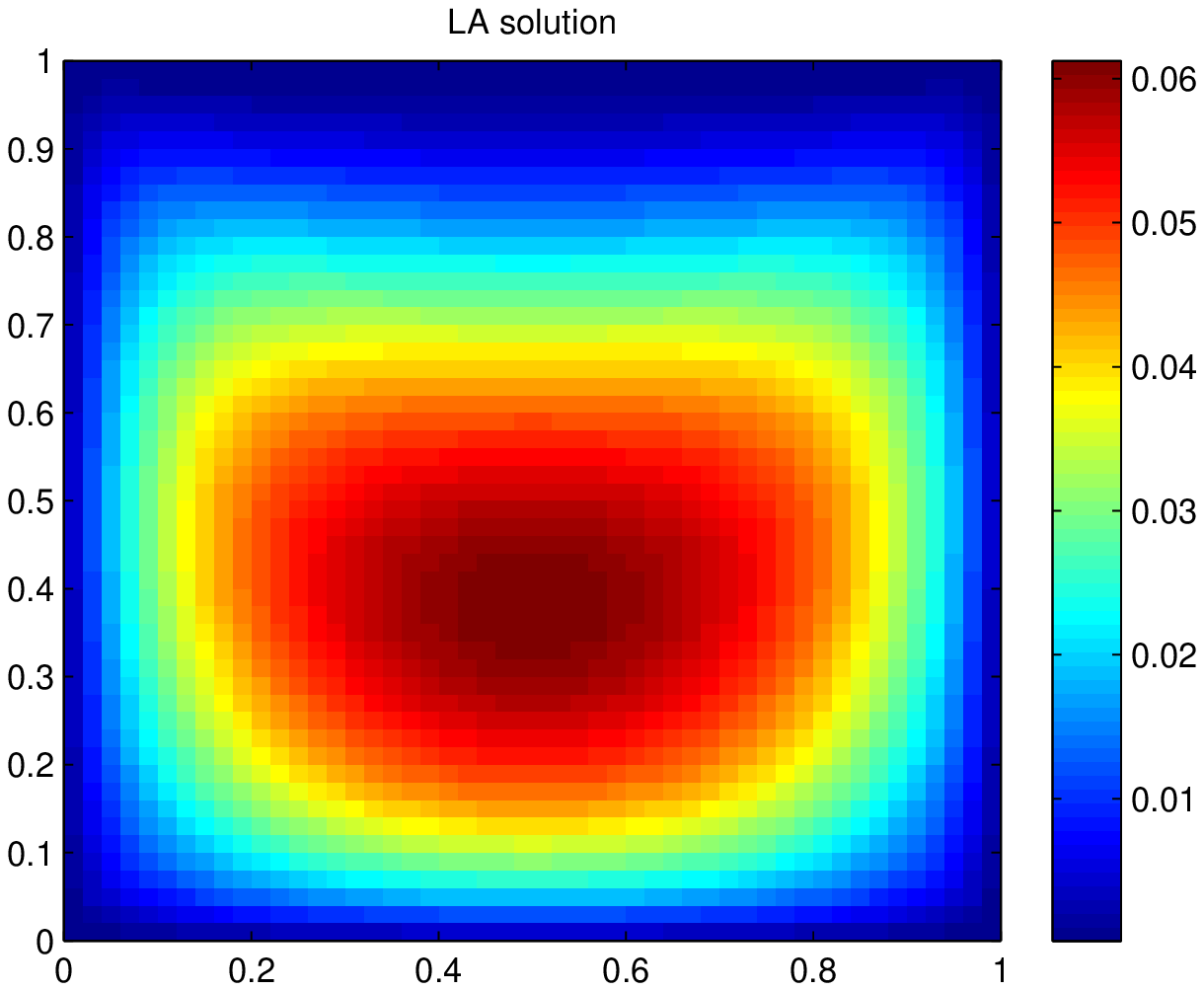}}
  \subfigure[Data on $\partial\Omega$]{
    \label{fig:subfig:d}
  \includegraphics[width=3in, height=2.4in]{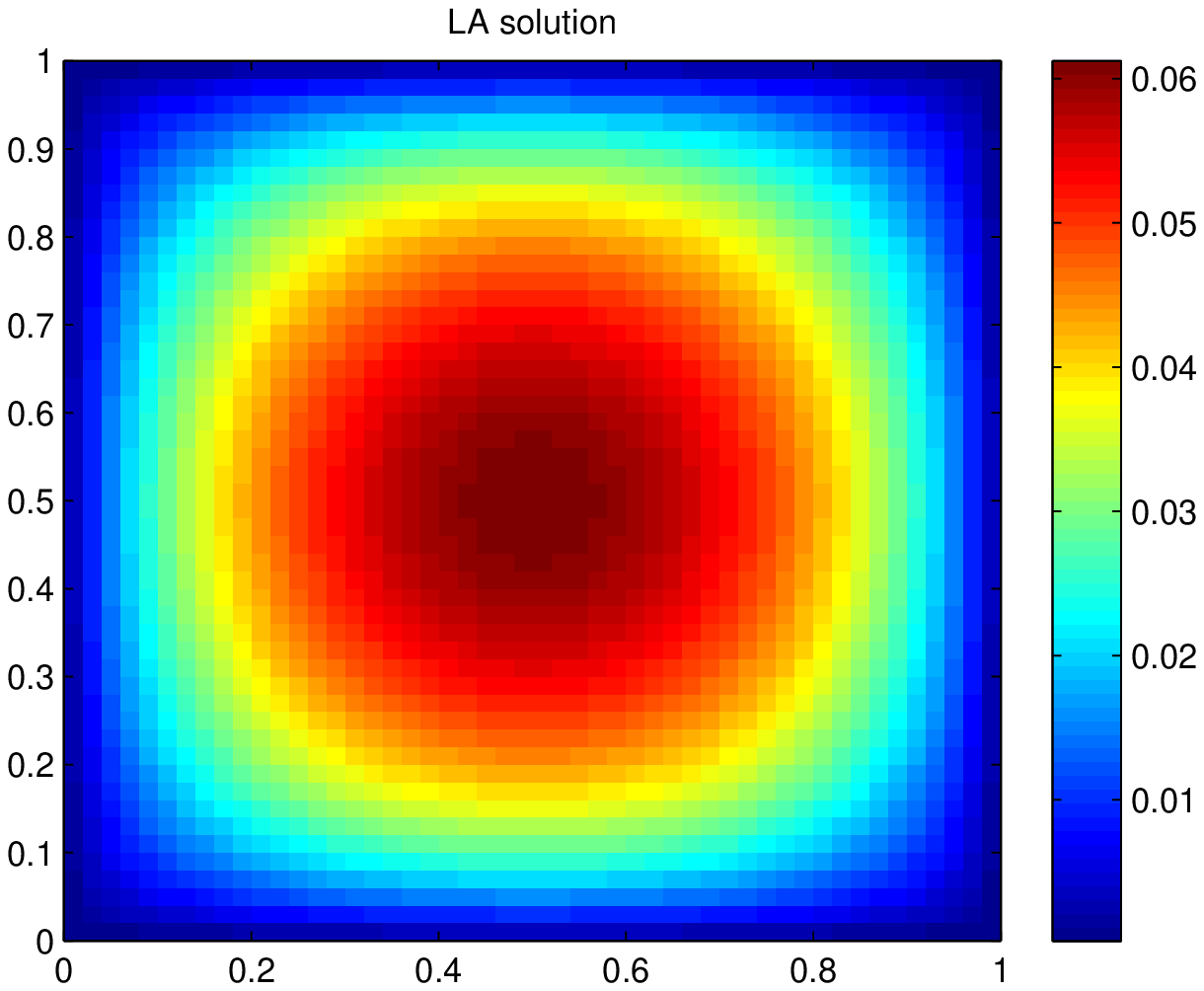}}\\
  \caption{The numerical result for Example 7 for noise levels $\varepsilon=0.0001$ with $\mu=\delta^{\frac{1}{2}}$}
  \label{two-example 2}
\end{figure}

$\mathbf{Example ~~7}. $ We take the $\Lambda=\partial\Omega$, the degradation coefficient is:
\begin{eqnarray*}
\label{Qresidual-eq}
\begin{split}
q(x,y)=x(1-x)y(1-y),\ \ x,y\in\Omega.
\end{split}
\end{eqnarray*}
It is easy to see $q(x,y)$ is central symmetry, and the reconstruction results are displayed in Figure $\ref{two-example 2}$. Here, Figure $\ref{two-example 2}(a)$ is the exact solution and the Figure $\ref{two-example 2}(b)$, $\ref{two-example 2}(c)$, $\ref{two-example 2}(d)$ are correspond to average flux data measured on the different part of boundary. The skewness is calculated from the marginal density function.  The data measured on boundary $\Lambda_1\cup\Lambda_2$, the graph is skewed to the front right of the region, that is, the skewness is negative in the x direction and positive in the y direction. When the average flux data measured on $\Lambda_1\cup\Lambda_2\cup\Lambda_4$, the graph slant toward the front of the area i.e., the skewness is zero in the x direction and positive in the y direction. We use the average flux data on the whole boundary $\partial\Omega$, the graph is accordance with the exact solution. The skewness correspond to the average flux data collected on different position of boundary can be seen in Table $\ref{tab1-two ske}$. Notice that here the corresponding marginal distributions of 2-D random variable are used to calculate the skewness by applying formula (\ref{ske}). In addition, the  95\% confidence intervals of LA solutions at different $(x,y)\in\Omega$ are shown in the Table 8. It can be seen that, owing to the appropriate regularization parameter selection and average flux measurement on whole boundary, the size of confidence region is very small, i.e., the proposed LA algorithm gives a higher reliability.

\begin{table}[H]
\small
\centering
\caption{ Numerical results for Example 7 for the skewness of  exact solution is 0}
\begin{tabular}{c|c}
\Xhline{1pt}
  Measurement data   &  skewness \\
  \hline
 Data on $\Lambda_1\cup\Lambda_2$ & $(-0.2905,0.2551)$ \\
 \hline
 Data on $\Lambda_1\cup\Lambda_2\cup\Lambda_4$  & $(-0.0170,0.1308)$ \\
 \hline
 Data on $\partial\Omega$  & $ (-0.0150,-0.0150)$\\
\Xhline{1pt}
\end{tabular}
\label{tab1-two ske}
\end{table}
\begin{table}[H]
\small
\centering
\caption{ Results for Example 7 for the 95\% confidence interval at different $(x,y)\in\Omega$}
\begin{tabular}{c|c|c|c}
\Xhline{1pt}
  y$\setminus x$ & $0.2$ & $0.5$& $0.8$\\
  \hline
 0.2 & $[3.1721\times10^{-2},3.1814\times10^{-2}]$  & $[4.5196\times10^{-2},4.5230\times10^{-2}]$ & $[3.1716\times10^{-2},3.1810\times10^{-2}]$\\
 \hline
 0.5 & $[4.4759\times10^{-2},4.4853\times10^{-2}]$  & $[6.1178\times10^{-2},6.1279\times10^{-2}]$ & $[4.4749\times10^{-2},4.4851\times10^{-2}]$ \\
 \hline
 0.8 & $[3.1713\times10^{-2},3.1812\times10^{-2}]$  & $[4.5174\times10^{-2},4.5268\times10^{-2}]$ & $[3.1711\times10^{-2},3.1801\times10^{-2}]$\\
\Xhline{1pt}
\end{tabular}
\label{confidence95 }
\end{table}

\section{Conclusions}
In this paper, we study the identification of degeneracy coefficients in time-fractional diffusion equations (TFDE) by using the average flux data at the accessible part of boundary. We mainly prove that the average flux measurement data can uniquely determine the degradation coefficient. The Lipchitz continuity of the corresponding forward operator is obtained. Due to the average flux measurement data only provide very limited information, and lead to serious ill-posedness of IDCP. This paper combines Tikhonov regularization with Laplace method to overcome the ill-posedness.  The existence, stability and convergence of solutions of variational problems are given. The paper introduce the sensitivity problem and the adjoint problem to find the minimizer of the variational problem by using the conjugate gradient method, and derive the mean and variance of the approximate posterior distribution by applying Bayesian theory and Laplace approximation. The Hellinger distance between the exact posterior measure and Laplace approximation is analyzed, and the second order convergence rate at MAP point is proved. The symmetry of the LA solution described by skewness is proposed, and find that the symmetry of solution is closely related to the symmetry of data. Finally, some numerical examples show that the method is not only accurate and flexible, but also can capture statistical information and quantify the uncertainty of the solution.

\section*{Acknowledgments}
The work described in this paper was supported by the NSF of China (11301168).

\label{ssec:Conclusions}


\begin{thebibliography}{99}
\bibitem{R. Metzler2000}
{\sc R. Metzler, J. Klafter}, {\em The random walk¡¯s guide to anomalous diffusion: a fractional dynamics
approach}, Phys. Rep. 339 (2000) 1-77.

\bibitem{I. Podlubny1999}
{\sc I. Podlubny}, {\em Fractional differential equations. An introduction to fractional
derivatives, fractional differential equations, some methods of their solution
and some of their applications}, Academic Press, 1999.

\bibitem{A.A.Kilbas2006}
{\sc A. Kilbas, H. Srivastava, J. Trujillo}, {\em Theory and Applications of Fractional Differential Equations}, Elsevier, Amsterdam, 2006.

\bibitem{S.G. Samko1993}
{\sc S. Samko, A. Kilbas, O. Marichev}, {\em Fractional Integrals and Derivatives}, Gordon and Breach Science Publishers, Philadelphia, 1993.

\bibitem{V. V. Uchaikin2013}
{\sc V. Uchaikin}, {\em Fractional Derivatives for Physicists and Engineers}, Springer, 2013.


\bibitem{D. Brockmann2006}
{\sc D. Brockmann, L. Hufnagel, T. Geisel}, {\em The scaling laws of human travel}, Nature 439 (2006) 462-465 .

\bibitem{E. Scalas2000}
{\sc E. Scalas, R. Gorenflo, F. Mainardi}, {\em Fractional calculus and continuous-time finance}, Phys. A. 284 (2000) 376-384.

\bibitem{D. Benson2000}
{\sc D. Benson, S. Wheatcraft, M. Meerschaert}, {\em Application of a fractional advection-dispersion equation}, Water Resour. Res. 36 (6) (2000) 1403-1412.

\bibitem{M. Hall2008}
{\sc M. Hall, T. Barrick}, {\em From diffusion-weighted MRI to anomalous diffusion imaging}, Magn. Reson. Med. 59 (3) (2008) 447-455.

\bibitem{B. Henry2008}
{\sc B. Henry, T. Langlands, S. Wearne}, {\em Fractional cable models for spiny neuronal dendrites}, Phys. Rev. Lett. 100 (12) (2008) 128103.


\bibitem{S.Y. Lukashchuk2015}
{\sc S. Lukashchuk}, {\em Conservation laws for time-fractional subdiffusion and diffusion-wave equations}, Nonlinear Dyn. 80 (2015) 1-12.

\bibitem{G. Wang2015}
{\sc G. Wang, A. Kara, K. Fakhar}, {\em Symmetry analysis and conservation laws for the class of time-fractional nonlinear dispersive equation}, Nonlinear Dyn. 82 (2015) 281-287.

\bibitem{R. Schumer2003}
{\sc R. Schumer, D. Benson, M. Meerschaert, B. Baeumer}, {\em Fractal mobile/immobile solute transport}, Water Resour. Res. 39 (2003) 1269.

\bibitem{P. Chakraborty2009}
{\sc P. Chakraborty, M. Meerschaert, C. Lim}, {\em Parameter estimation for fractional transport: a particle tracking approach}, Water Resour. Res. 45 (2009) W10415.

\bibitem{K. Sakamoto2011}
{\sc K. Sakamoto, M. Yamamoto}, {\em Initial value/boundary value problems for fractional diffusion-wave equations and applications to some inverse problems}, J. Math. Anal. Appl. 382 (2011) 426-447.

\bibitem{S.D. Eidelman2004}
{\sc S. Eidelman, A. Kochubei}, {\em Cauchy problem for fractional diffusion equations}, J. Differential Equations 199 (2004) 211-255.

\bibitem{Y. Luchko2010}
{\sc Y. Luchko}, {\em Some uniqueness and existence results for the initial-boundary value problems for the generalized time-fractional diffusion equation}, Comput. Math. Appl. 59 (2010) 1766-1772.

\bibitem{R. Gorenflo2015}
{\sc R. Gorenflo, Y. Luchko, M. Yamamoto}, {\em Time-fractional diffusion equation in the fractional Sobolev spaces}, Fract. Calc. Appl. Anal. 18 (2015) 799-820.

\bibitem{Y. Lin and C. Xu2007}
{\sc Y. Lin, C. Xu }, {\em Finite difference/spectral approximations for the time-fractional diffusion equation}, J. Comput. Phys. 225 (2) (2007) 1533-1552.

\bibitem{Y. N. Zhang2014}
{\sc Y. Zhang, Z. Sun, H. Liao}, {\em Finite difference methods for the time fractional diffusion equation on non-uniform meshes}, J. Comput. Phys. 265 (2014) 195-210.

\bibitem{F. Zeng2013}
{\sc F. Zeng, C. Li, F. Liu, I. Turner }, {\em The use of finite difference/element approaches for solving the time-fractional subdiffusion equation}, SIAM J. Sci. Comput. 35 (6) (2013) A2976-A3000.

\bibitem{K. Mustapha2016}
{\sc K. Mustapha, M. Nour, B. Cockburn}, {\em Convergence and superconvergence analyses of HDG methods for time fractional diffusion problems,}, Adv. Comput. Math. 42 (2) (2016) 377-393.

\bibitem{Q. Xu2013}
{\sc Q. Xu, Z. Zheng}, {\em Discontinuous Galerkin method for time fractional diffusion equation}, J. Informat. Comput. Sci. 10 (2013) 3253-3264.

\bibitem{B. Jin2015}
{\sc B. Jin, R. Lazarov, J. Pasciak, Z. Zhou }, {\em Error analysis of semidiscrete finite element methods for inhomogeneous time-fractional diffusion}, IMA J. Numer. Anal. 35 (2) (2015) 561-582.

\bibitem{B. JinRundell2015}
{\sc B. Jin, W. Rundell}, {\em A tutorial on inverse problems for anomalous diffusion processes}, Inverse Probl. 31 (3) (2015) 035003.

\bibitem{J.J. Liu2010}
{\sc J. Liu, M. Yamamoto}, {\em A backward problem for the time-fractional diffusion equation }, Appl. Anal. 89(11) (2010) 1769-1788.

\bibitem{D.A. Murio2008}
{\sc D. Murio}, {\em Time fractional IHCP with Caputo fractional derivatives }, Comput. Math. Appl. 56 (2008) 2371-2381.

\bibitem{J. Liu 2015}
{\sc J. Liu, M. Yamamoto, L. Yan}, {\em On the uniqueness and reconstruction for an inverse problem of the fractional diffusion process}, Appl. Numer. Math. 87 (2015) 1-19.

\bibitem{G.H. Zheng2012}
{\sc G. Zheng, T. Wei}, {\em A new regularization method for a Cauchy problem of the time fractional diffusion equation }, Adv. Comput. Math. 36 (2) (2012) 377-398.

\bibitem{W. Rundell 2013}
{\sc W. Rundell, X. Xu, L. Zuo}, {\em The determination of an unknown boundary condition in a fractional diffusion equation}, Appl. Anal. 92 (7) (2013) 1511-1526.

\bibitem{T. Wei2016(1)}
{\sc T. Wei, X. Li, Y. Li}, {\em An inverse time-dependent source problem for a time-fractional diffusion equation}, Inverse Probl. 32 (8) (2016) 085003.

\bibitem{Y. Zhang2011}
{\sc Y. Zhang, X. Xu}, {\em Inverse source problem for a fractional diffusion equation}, Inverse Probl. 27(3) (2011) 035010.

\bibitem{Liu2016(1)}
{\sc Y. Liu, W. Rundell, M. Yamamoto}, {\em Strong maximum principle for fractional diffusion equations and an application to an inverse source problem}, Fract. Calc. Appl. Anal. 19 (2016) 888-906.

\bibitem{Wei2014(1)}
{\sc T. Wei, J. Wang}, {\em A modified quasi-boundary value method for an inverse source problem of the time-fractional diffusion equation}, Appl. Numer. Math. 78 (2014) 95-111.

\bibitem{N.H. Tuan2017}
{\sc N. Tuan, M. Kirane, L. Hoan, L. Long}, {\em Identification and regularization for unknown source for a time-fractional diffusion equation}, Comput. Math. Appl. 73 (2017) 931-950.

\bibitem{J. Cheng2009}
{\sc J. Cheng, J. Nakagawa, M. Yamamoto, T. Yamazaki}, {\em Uniqueness in an inverse problem for a one-dimensional fractional diffusion equation},
Inverse Probl. 25 (11) (2009) 115002, 16.

\bibitem{G.S. Li2013}
{\sc G. Li, D. Zhang, X. Jia, M. Yamamoto}, {\em Simultaneous inversion for the space-dependent diffusion coefficient and the fractional order in
the time-fractional diffusion equation}, Inverse Probl. 29 (6) (2013) 065014.

\bibitem{Z.D. Zhang2016}
{\sc Z. Zhang}, {\em An undetermined coefficient problem for a fractional diffusion equation}, Inverse Probl. 32 (1) (2016) 015011.

\bibitem{B.T. Jin2012}
{\sc B. Jin, W. Rundell}, {\em An inverse problem for a one-dimensional time-fractional diffusion problem}, Inverse Probl. 28 (7) (2012) 075010.

\bibitem{L. MillerYamamoto2013}
{\sc L. Miller, M. Yamamoto}, {\em Coefficient inverse problem for a fractional diffusion equation}, Inverse Probl. 29 (7) (2013) 075013.

\bibitem{V.K. Tuan2011}
{\sc V. Tuan}, {\em Inverse problem for fractional diffusion equation}, Fract. Calc. Appl. Anal. 14 (1) (2011) 31-55.

\bibitem{Z. Li2016}
{\sc Z. Li, O. Imanuvilov, M. Yamamoto}, {\em Uniqueness in inverse boundary value problems for fractional diffusion equations}, Inverse Probl. 32 (2016) 015004.

\bibitem{T. Wei2018}
{\sc T. Wei, Y. Li}, {\em Identifying a diffusion coefficient in a time-fractional diffusion equation}, Math. Comput. Simul. 151 (2018) 77-95

\bibitem{L.sun2017}
{\sc L. Sun, T. Wei}, {\em Identification of the zeroth-order coefficient in a time fractional diffusion equation}, Applied Numerical Mathematics 111 (2017) 160-180.

\bibitem{SunYanWei2018}
{\sc L. Sun, X. Yan, T. Wei}, {\em Identification of time-dependent convection coefficient in a time-fractional diffusion equation}, J. Appl. Math. Comput. (2018), https://doi.org/10.1016/j.cam.2018.07.029

\bibitem{T. Wei2016}
{\sc T. Wei, J. Wang}, {\em Determination of Robin coefficient in a fractional diffusion problem}, Appl. Math. Model. 40 (2016) 7948-7961.

\bibitem{Weizhang2016}
{\sc T. Wei, Z. Zhang}, {\em Robin coefficient identification for a time-fractional diffusion equation}, Inverse Probl. Sci. Eng. 24 (4) (2016) 647-666.

\bibitem{YamamotoZhang2012}
{\sc M. Yamamoto, Y. Zhang}, {Conditional stability in determining a zeroth-order coefficient in a half-order fractional diffusion equation by a Carleman estimate}, Inverse Probl. 28 (10) (2012) 105010.

\bibitem{Wang2006}
{\sc C. Wang, Y. Wang, P. Wang}, {\em Water quality modeling and oollution control for the eastern route of South to North Water Transfer Project in China}, Journal of Hydrodynamics, Ser. B 18 (3) (2006) 253-261.

\bibitem{Huang2017}
{\sc B. Huang, C. Hong, H. Du, J. Qiu, X. Liang, C. Tan, D. Liu}, {\em Quantitative study of degradation coefficient of pollutant against the flow
velocity}, Journal of Hydrodynamics 29 (1) (2017) 118-123.

\bibitem{Dostert2009}
{\sc P. Dostert, Y. Efendiev, B. Mohanty}, {\em Efficient uncertainty quantification techniques in inverse problems for Richards¡¯ equation using coarse-scale simulation models}, Adv. Water Resour. 32 (3) (2009) 329-339.

\bibitem{McEnroe2009}
{\sc N. McEnroe, N. Roulet, T. Moore, M. Garneau}, {\em Do pool surface areaand depth control CO2 and CH4 fluxes from an ombrotrophic raised bog, James Bay, Canada?}, J. Geophys. Res. Biogeosci. 114 (2009) http://dx.doi.org/10.1029/2007jg000639.

\bibitem{Bao2013}
{\sc G. Bao, P. Li, J. Lv}, {\em Numerical solution of an inverse diffraction grating problem from
phaseless data}, J. Opt. Soc. Amer. A 30 (2013) 293-299.

\bibitem{Ammari2016}
{\sc H. Ammari, Y. Chow, J. Zou}, {\em Phased and phaseless domain reconstructions in the inverse scattering problem via scattering coefficients}, SIAM J. Appl. Math. 76 (2016) 1000-1030.

\bibitem{Zhang2017}
{\sc B. Zhang, H. Zhang}, {\em Recovering scattering obstacles by multi-frequency phaseless far-field data}, J. Comput. Phys. 345 (2017) 58-73.

\bibitem{Al-Refai2015}
{\sc M. Al-Refai, Y. Luchko}, {\em Maximum principle for the multi-term time-fractional diffusion equations with the Riemann-Liouville fractional derivatives},
Appl. Math. Comput. 257 (2015) 40-51.

\bibitem{Badia1999}
{\sc A. Badia}, {\em Coefficient identification in some partial differential equations from partial boundary measurements}
Inverse Probl. 15 (1999) 11-18.

\bibitem{ford2011}
{\sc N.Ford, J. Xiao, Y. Yan}, {\em A finite element method for time fractional partial differential equations}, Fract. Calc. Appl. Anal. 14 (2011) 454-474.

\bibitem{hofmann2007}
{\sc B. Hofmann, B. Kaltenbacher, C. P\"{o}schl, O. Scherzer}, {\em A convergence rates for Tikhonov regularization in Banach spaces with non-smooth operators}, Inverse Probl. 23 (3) (2007) 987-1010.

\bibitem{brezis2011}
{\sc H. Brezis}, {\em Functional Analysis, Sobolev Spaces and Partial Differential Equations}, Springer, New York, 2011.

\bibitem{stuart2010}
{\sc A. Stuart}, {\em Inverse problems: A Bayesian perspective}, Acta Numer. 19 (2010) 451-559.

\bibitem{iglesias2010}
{\sc M. Iglesias, K. Law, A. Stuart}, {\em  Evaluation of Gaussian approximations for data
assimilation in reservoir models}, Comput. Geosci. 17 (2013) 851-885.

\bibitem{cotter2010}
{\sc S. Cotter, M. Dashti, A. Stuart}, {\em  Approximation of Bayesian inverse problems
for PDEs}, SIAM J. Numer. Anal. 48 (2010) 322-345.

\bibitem{wacher2017}
{\sc P. Wacker}, {\em Laplace's method in Bayesian inverse problems with Gaussian
priors}, arXiv preprint arXiv:1701.07989, (2017).

\bibitem{yan2015}
{\sc L. Yan, L. Guo}, {\em Stochastic collocation algorithms using l1-minimization for Bayesian
solution of inverse problems}, SIAM J. Sci. Comput. 37 (2015) A1410-A1435.

\bibitem{jiang2017}
{\sc L. Jiang, N. Ou}, {\em Multiscale model reduction method for Bayesian inverse problems of
subsurface flow}, J. Comput. Appl. Math. 319 (2017) 188-209.

\bibitem{Reynolds2008}
{\sc A. Reynolds, D. Oliver, N. Liu}, {\em Inverse Theory for
Petroleum Reservoir Characterization and History Matching}, 1st
edn. ISBN:9780521881517. Cambridge University Press, Cambridge,
2008.

\bibitem{draper1981}
{\sc N. Draper, H. Smith}, {\em Applied Regression Analysis}, John Wiley and Sons, New York, (1981)

\bibitem{Everitt2006}
{\sc B. Everitt}, {\em The Cambridge dictionary of statistics (Third
edition)}, Cambridge University Press, Cambridge, (2006)

\bibitem{vogel2002}
{\sc C. Vogel}, {\em Computational methods for inverse problems}, Society for Industrial and
Applied Mathematics, 2002.

\bibitem{jin2007}
{\sc B. Jin}, {\em Conjugate gradient method for the Robin inverse problem associated with the Laplace
equation}, Internat. J. Numer. Methods Engrg. 71 (2007) 433-453.

%
%
%

\bibitem{Lions 1972Non-homogeneous}
{\sc J. Lions, E. Magenes}, {\em  Non-homogeneous Boundary Value Problems and Applications}, Volume 1, Springer-Verlag, 1972.







\end{thebibliography}
\end{document}